\documentclass[runningheads, numbook, envcountsame]{svjour2}
\usepackage{geometry}
\usepackage{graphicx}
\usepackage{epstopdf}
 \usepackage{amsfonts}
 \usepackage{amsmath}
 \usepackage{amssymb}
 \usepackage{mathrsfs}
 \usepackage{txfonts}
 \usepackage{color}
\newcommand{\be}{\begin{equation} }
\newcommand{\ee}{\end{equation}}
\newcommand{\sech}{{\rm sech} }
\newcommand{\csch}{{\rm csch} }

\newenvironment{longremark}{\begin{remark}}{ \qed \end{remark}}
\newenvironment{indentenum}{\begin{enumerate} \item[] \begin{enumerate}}{\end{enumerate} \end{enumerate}}
\newenvironment{proof2}{\begin{proof}}{\qed \end{proof}}
\newenvironment{proofof}[1]{\paragraph{Proof of #1}}{\qed}

\smartqed  

\journalname{Archive for Rational Mechanics and Analysis}

\begin{document}

\title{Steady periodic gravity waves with surface tension
}

\titlerunning{Water waves with surface tension}        

\author{Samuel Walsh}


\institute{Brown University \at
               182 George St., Providence, RI 02912, USA\\
              Tel.: (401) 863-3113\\
              \email{samuel\_walsh@brown.edu}           
}

\date{\notused}

\maketitle

\begin{abstract}
In this paper we consider two-dimensional, stratified, steady water waves propagating over an impermeable flat bed and with a free surface.  The motion is assumed to be driven by capillarity (that is, surface tension) on the surface and a gravitational force acting on the body of the fluid.   We prove the existence of global continua of classical solutions that are periodic and traveling.   This is accomplished by first constructing a 1-parameter family of laminar flow solutions, $\mathcal{T}$, then applying bifurcation theory methods to obtain local curves of small amplitude solutions branching from $\mathcal{T}$ at an eigenvalue of the linearized problem.  Each solution curve is then continued globally by means of a degree theoretic theorem in the spirit of Rabinowitz. Finally, we complement the degree theoretic picture by proving an alternate global bifurcation theorem via the analytic continuation method of Dancer.
\end{abstract}

\section{Introduction}
\label{intro}

Initiated by the breakthrough work of Constantin and Strauss (cf. \cite{constantin2004exact}), there has recently been a great deal of interest in the theory of traveling water waves with vorticity (see, e.g., \cite{wahlen2006capgrav,wahlen2006capillary,wahlen2007onrotational,wahlenthesis,walsh2009stratified}.)  Historically, of course, most investigations of water waves have focused on the irrotational case, which has several extremely convenient mathematical features.  While appropriate in a wide variety of scenarios, however, the assumption of zero vorticity is far from universally valid.  For instance, ocean waves typically have a rotational layer near the surface due to the effects of wind shear (cf. \cite{mei1989applied}).

In fact, the wind is known to itself generate steady waves.  When it first begins to blow over a body of still water, the wind produces a fine, hexagonal meshwork of wavelets whose motion is governed entirely by capillarity (i.e. surface tension.)  If the breeze is maintained, then these waves will grow larger, eventually evolving into capillary-gravity waves and finally into gravity waves.  Over the course of this progression, the waves develop a remarkable robustness:  while pure capillary waves will quickly dissipate in the absence of external forcing, capillary-gravity and gravity waves can propagate long distances unaided (cf. \cite{kinsman}). 

Surface tension's central role in the formation of wind-driven waves argues strongly for the mathematical consideration of steady, capillary-gravity waves with vorticity, particularly since wind-waves are one of the most important classes of rotational water waves.  The goal of the present paper, therefore, is to develop a global theory for traveling, periodic, water waves with surface tension building upon the local theory derived by Wahl\'en (cf. \cite{wahlen2006capgrav,wahlen2007onrotational,wahlenthesis}).

Briefly, our new results come in five parts.  The first two of these relate to the phenomenon of double (local) bifurcation for capillary-gravity waves.  It is well known in the theory of irrotational capillary-gravity waves that, when the coefficient of surface tension resides in a distinguished set, the eigenvalue of the linearized problem at the point of bifurcation will be of multiplicity two and is otherwise simple.  This distinctive feature was first observed in the irrotational case by Wilton (cf. \cite{wilton1915ripples}) and has been the subject of numerous papers since. Wahl\'en observed the presence of double bifurcation points in the rotational case and proved the existence of local bifurcation curves in certain regimes (cf. \cite{wahlen2006capgrav}).  These results, however, are incomplete in the sense that they do not necessarily describe the entire solution set near the bifurcation point.   The first of our contributions is to describe in detail the local bifurcation diagram obtained in \cite{wahlen2006capgrav} at double bifurcation points.  In particular we confirm the existence of curves of mixed solutions conjectured by Wahl\'en (see the discussion following \textsc{Theorem \ref{mainresult1}}).  Secondly, using the bifurcation theory for analytic operators pioneered by Dancer, we extend the local curves to obtain global continua of solutions.  The application of these techniques to steady rotational waves is in and of itself novel and pays some additional dividends.  Most notably, it leads to our third result, which is proving that the solution continua are path-connected, a fact that was unknown even for the constant density, gravity wave case originally considered in \cite{constantin2004exact}.  Fourthly, for the case of simple local bifurcation, we also use the more standard degree theoretic method to obtain a global bifurcation theorem.  Finally, employing the machinery developed by the author in \cite{walsh2009stratified}, all of these results are obtained for fluids with general (stable) stratification, which is entirely new, even for the local bifurcation theory.  This improvement is quite important from a physical standpoint, since stratification phenomenon is both very common in ocean waves (due to salinity and temperature variations), and exerts a potentially significant effect on the dynamics of the flow (see, e.g., \cite{yih1965dynamics}, or the discussion in \cite{walsh2009stratified}.)   

That said, let us now recall the setup for traveling, stratified, capillary-gravity waves.  If we imagine a wave on the open ocean, past experience suggests that it may be \emph{regular} in the following sense.  First, it is essentially two-dimensional.  That is, the motion will be identical along any line that runs parallel to a crest.  Second, if we regard the wave in a coordinate system moving with some constant speed, it appears steady.  Finally, the profile is periodic in the direction of motion.
\begin{figure}[h] \setlength{\unitlength}{5cm} 
\centering
\includegraphics[scale=0.5]{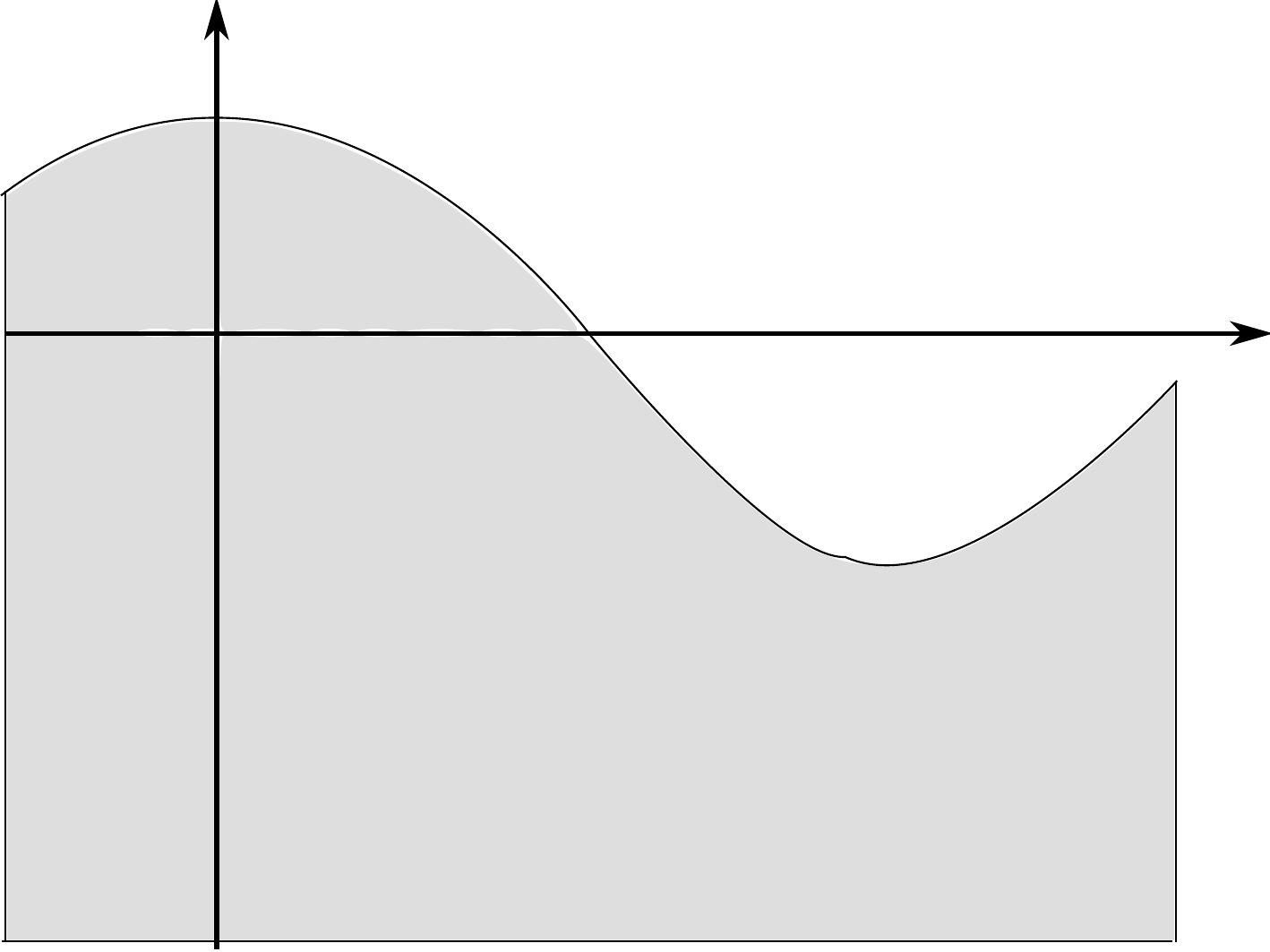}
\put(-0.85,0.86){$y = \eta(x,t)$}
\put(0.00,0.00){$y = -d$}
\put(0.00,0.60){$x$}
\put(-1.15,1.05){$y$}
\put(-0.75,0.30){$\overline{D_\eta}$}
\caption{The fluid domain $\overline{D_\eta}$.  The wave is assumed to propagate to the left with speed $c$.}
\end{figure}
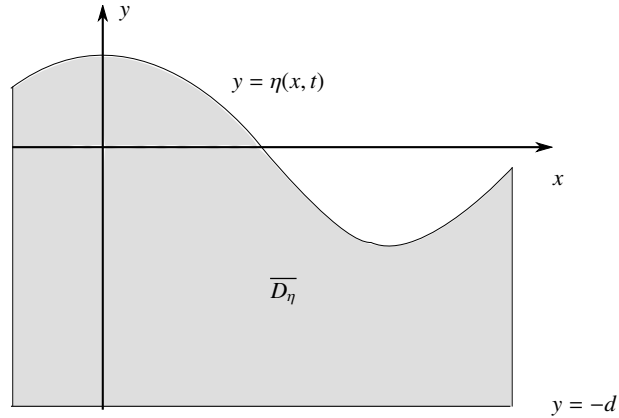

We now formulate governing equations for waves of this form.  Fix a Cartesian coordinate system so that the $x$-axis points in the direction of propagation, and the $y$-axis is vertical.  We assume that the floor of the ocean is flat and occurs at $y = -d$.  Let $y = \eta(x,t)$ be the free surface at the interface between the atmosphere and the fluid.  We shall normalize $\eta$ by choosing the axes so that the free surface is oscillating around the line $y = 0$.  As usual we let $u = u(x,y,t)$ and $v = v(x,y,t)$ denote the horizontal and vertical velocities, respectively.  Let $\rho = \rho(x,y,t) > 0$ be the density.  

Incompressibility of the fluid is represented mathematically by the requirement that the vector field $(u,v)$ be divergence free for all time:
\be u_x + v_y = 0. \label{incompress} \ee
Taking the fluid to be inviscid, conservation of mass implies that the density of a fluid particle remains constant as it follows the flow.  This is expressed by the continuity equation
\be \rho_t + u \rho_x + v \rho_y = 0. \label{consmass1} \ee
Next, the conservation of momentum is described by Euler's equations,
\be \left \{ \begin{array}{lll}
\rho u_t + \rho u u_x + \rho v u_y  & = & -{P_x}{}, \\
\rho v_t + \rho u v_x + \rho v v_y & = & -{P_y}- g\rho, \end{array} \right. \label{euler1} \ee
where $P=P(x,y,t)$ denotes the pressure and $g$ is the gravitational constant.  Here, of course, we assume that the only external force acting on the fluid is gravity.  

For capillary waves, we formulate the boundary condition on the free surface as folows.  At each point we say that the pressure of the fluid deviates from the atmospheric pressure of the air above, $P_{\textrm{atm}}$, in proportion to the (mean) curvature at that point.  This results in the dynamic boundary condition
\be P = P_{\textrm{atm}}-\sigma \frac{\eta_{xx}}{\left(1+\eta_x^2\right)^{\frac{3}{2}}}, \qquad \textrm{on } y = \eta(x,t), \label{presscond} \ee
where $\sigma > 0$ is the coefficient of surface tension, a material property that we take to be given and constant.

The corresponding boundary condition for the velocity is motivated by the fact that fluid particles that reside on the free surface continue to do so as the flow develops.  This observation is manifested in the kinematic condition:
\be v = \eta_t + u \eta_x, \qquad \textrm{on } y = \eta(x,t). \label{kincond1} \ee
Since we cannot have any fluid moving normal to the flat bed occurring at $y = -d$, we require
\be v = 0, \qquad \textrm{on } y = -d. \label{bedcond} \ee
Note that there is no accompanying condition on $u$ because in the inviscid case we allow for slip, that is, nonzero horizontal velocity along solid boundaries.  

We seek traveling periodic wave solutions $(u,v,\rho,P, \eta)$ to \eqref{incompress}-\eqref{bedcond}.  More precisely, we take this to mean that, for fixed $c > 0$, the solution appears steady in time and periodic in the $x$-direction when observed in a frame that moves with constant speed $c$ to the right.  The vector field will thus take the form $u = u(x-ct,y)$,  $v = v(x-ct, y)$, where each of these is $L$-periodic in the first coordinate.  Likewise for the scalar quantities: $\rho(x,y,t) = \rho(x-ct,y)$, $P = P(x-ct,y)$, and $\eta = \eta(x-ct)$, again with $L$-periodicity in the first coordinate.  We therefore take moving coordinates 
\[ (x-ct, y) \mapsto (x,y), \]
which eliminates time dependency from the problem.  In the moving frame \eqref{incompress}-\eqref{euler1} become 
\be \left \{ \begin{array}{lll} 
u_x + v_y & = & 0, \\
(u-c)\rho_x + v\rho_y & = & 0, \\
\rho (u-c)u_x + \rho v u_y  & = & -P_x ,\\
\rho(u-c) v_x + \rho v v_y & = & -P_y - g\rho, \end{array} \right. \label{euler2} \ee
throughout the fluid domain.  Meanwhile, the reformulated boundary conditions are
\be \left \{ \begin{array}{llll} 
v & = & (u-c) \eta_x, & \textrm{on } y = \eta(x), \\
v & = & 0, & \textrm{on } y = -d, \\
P & = & P_{\textrm{atm}} - \sigma \eta_{xx}\left(1+\eta_x^2\right)^{-3/2}, & \textrm{on } y = \eta(x), \end{array} \right. \label{boundcond} \ee
where $u,v,\rho, P$ are taken to be functions of $x$ and $y$, $\eta$ is a function of $x$, and all of them are $L$-periodic in $x$.

In the event that $u = c$ somewhere in the fluid we say that \emph{stagnation} has occurred.  This is a slight abuse of terminology, since in our definition we require only that the horizontal velocity vanishes.  In order to avoid roll-up and other instability phenomena, we shall restrict our attention to the case where $u < c$ throughout.  

Recall that we have chosen our axes so that $\eta$ oscillates around the line $y = 0$.  Mathematically this equates to
\be \fint_{0}^{L} \eta(x)dx = 0. \label{normalsurface} \ee 
In effect, this couples the depth to the problem, so we need not treat $d$ as a free parameter.  The trade-off, as we shall see, is that this normalization will have significant technical consequences later.  

Observe that by conservation of mass and incompressibility, $\rho$ is transported and the vector field is divergence free.  Therefore we may introduce a (relative) \emph{pseudo-stream function} $\psi = \psi(x,y)$, defined uniquely up to a constant by: 
\[ \psi_x = -\sqrt{\rho}v,\qquad \psi_y = \sqrt{\rho} (u-c). \]
Here we have the addition of a factor involving $\rho$ to the typical definition of the stream function for an incompressible fluid.  This neatly captures the inertial effects of the heterogeneity of the flow (see the treatment in \cite{yih1965dynamics}, for example).  The particular choice of $\sqrt{\rho}$ is merely to simplify algebraically what follows.

It is a straightforward calculation to check that $\psi$ is indeed a (relative) stream function in the usual sense, i.e. its gradient is orthogonal to the vector field in the moving frame at each point in the fluid domain:
\[ (u-c) \psi_x + v \psi_y = 0. \]
Moreover, \eqref{boundcond} implies that the free surface and flat bed are each level sets of $\psi$.  For definiteness we choose $\psi \equiv 0$ on the free boundary, so that $\psi \equiv -p_0$ on $y = -d$, where $p_0$ is the \emph{(relative) pseudo-voumetric mass flux}: 
\be p_0 := \int_{-d}^{\eta(x)} \sqrt{\rho(x,y)}\left[u(x,y) - c\right] dy. \label{defp0} \ee
It is an easy calculation to check that $p_0$ is well-defined (cf. \cite{walsh2009stratified}).  Physically, $p_0$ describes the rate of fluid moving through any vertical line in the fluid domain (and with respect to the transformed vector field $\sqrt{\rho}(u-c,v)$.)

Since $\rho$ is transported, it must be constant on the streamlines and hence, we may think of it as a function of $\psi$.  Abusing notation we may let $\rho:[p_0,0] \to \mathbb{R}^{+}$ be given such that 
\be \rho(x,y) = \rho(-\psi(x,y)) \label{defrho} \ee
throughout the fluid. When there is risk of confusion, we shall refer to the $\rho$ occurring on the right-hand side above as the \emph{streamline density function}.  We shall focus our attention on the case where the density is nondecreasing as depth increases.  This is entirely reasonable from a physical standpoint.  Indeed, hydrodynamic stability requires that the depth be monotonically increasing with depth, so that fluids with this type of density distribution are referred to in the literature as \emph{stably stratified}.  The level set  $-\psi = p_0$ corresponds to the flat bed, and the set where $-\psi = 0$ corresponds to the free surface.  Therefore, we require that the streamline density function is nonincreasing as a function of $-\psi$.  Note that this does not preclude the case of constant density.  

By Bernoulli's theorem (for gravity waves) the quantity
\[ E := P + \frac{\rho}{2}\left( (u-c)^2 + v^2\right) + g\rho y, \]
is constant along streamlines, that is 
\[ (u-c) E_x + v E_y = 0.\]
This can be verified directly from \eqref{euler2}. In the case of an inviscid fluid, $E$ represents the energy of the fluid particle at $(x,y)$.  The first term on the right-hand side, $P$, gives the energy due to internal pressure.  The second and third terms combined describe the kinetic energy, while the last is gravitational potential energy.  In particular, evaluating the above equation on the level set $\psi = 0$, we find
\be |\nabla \psi|^2 + 2g \rho(x, \eta(x)) \left(\eta(x)+d\right) = Q, \qquad \textrm{on } y = \eta(x) \label{capgravdefQ} \ee
where the constant $Q := 2(E|_\eta-P_{\textrm{atm}} + g\rho|_\eta d)$ gives roughly the energy density along the free surface of the fluid.  We shall view $Q$ as parameterizing the continuum of solutions.

Under the assumption that $u < c$ throughout the fluid, and given the fact that $E$ is constant along streamlines, there exists a function $\beta:[0,|p_0|] \to \mathbb{R}$ such that 
\be \frac{dE}{d\psi}(x,y) = -\beta(\psi(x,y)). \label{defbeta} \ee
For want of a better name we shall refer to $\beta$ as the \emph{Bernoulli function} corresponding to the flow.  Physically it describes the variation of specific energy as a function of the streamlines.   Define
\[ B(p) := \int_0^p \beta(-s)ds\]
for $p_0 \leq p \leq 0$ and let $B$ have minimum value $B_{\textrm{min}}$.  

We now briefly outline a few notational conventions.  Let $\overline{D_\eta}$ denote the closure of the fluid domain
\[ D_\eta := \left\{ (x,y) \in \mathbb{R}^2 : -d < y < \eta(x) \right\}.\]
For any integer $m \geq 1$ and $\alpha \in (0,1)$, we say that a bounded domain $D \subset \mathbb{R}^2$ is $C^{m+\alpha}$ provided that at each point in the boundary, denoted $\partial D$, is locally the graph of a $C^{m+\alpha}$ function.  That is, a function with H\"older continuous derivatives (of exponent $\alpha$) and up to order $m$.  Furthermore, for fixed $m \geq 1$ and $\alpha \in (0,1)$ we define the space $C_{\textrm{per}}^{m+\alpha}(\overline{D})$ to consist of those functions $f:\overline{D}\to\mathbb{R}$ with H\"older continuous derivatives (of exponent $\alpha$) up to order $m$ and that are $L$-periodic in the first coordinate.  Similarly, we shall take $C_{\textrm{per}}^{m}(\overline{D})$ to be the space of $m$-times continuously differentiable functions which are $L$-periodic in the first coordinate.  In this paper, we shall seek solutions to the water wave problem $(u,v,\rho, \eta, Q)$ in the space
\be \mathscr{S}:= C_{\textrm{per}}^{2+\alpha}(\overline{D_\eta}) \times C_{\textrm{per}}^{2+\alpha}(\overline{D_\eta}) \times C_{\textrm{per}}^{2+\alpha}(\overline{D_\eta}) \times C_{\textrm{per}}^{3+\alpha}(\mathbb{R}) \times \mathbb{R}. \label{defsolutionspace} \ee
\noindent
Occasionally, for brevity, we shall abuse notation and write $(u,v,\rho,\eta) \in \mathscr{S}$, or simply $(u,v,\eta) \in \mathscr{S}$, by which we mean $(u,v,\rho,\eta,Q) \in \mathscr{S}$, where $\rho$ is found from the (given) streamline density function using \eqref{defrho} and $Q$ is calculated via \eqref{capgravdefQ}.  

Finally, in describing the structure of the continua of solutions, we shall make use of the following definition.
\begin{definition} A subset $\mathcal{K} \subset \mathscr{S}$ is said to admit a \emph{global, locally injective, continuous parameterization} provided that $\mathcal{K}$ can be described as the image of a locally injective $\kappa \in C([0,\infty); \mathscr{S})$.  Furthermore, for $m \geq 1$, we say that $\mathcal{K}$ has a \emph{locally injective, $C^m$-reparameterization (resp. analytic reparameterization)} if, for each $s_* \in (0,\infty)$, there exists a continuous $\tau:(-1,1) \to \mathbb{R}^+$ such that $\tau(0) = s_*$ and the map $s \mapsto \kappa(\tau(s))$ is injective and of class $C^m$ (resp. analytic).  \label{defparameterization}
\end{definition}
There are several key features distinguishing the sets described above from the standard $C^m$-curves of differential geometry.  Observe that, while $\kappa$ may be locally injective, \textsc{Definition \ref{defparameterization}} does not require it to be \emph{globally} injective, i.e., the curve $\mathcal{K}$ may self-intersect.  This scenario is not permissible for $C^m$-curves since in the neighborhood of a point of self-intersection $\mathcal{K}$ fails to be a one-dimensional submanifold of $\mathscr{S}$.  Also, it is important to note that we do not require $\kappa^\prime$ to be non-vanishing.  As a consequence, the global chart may not be of class $C^m$ (resp. analytic), even though at each parameter value we can reparameterize to obtain that level of smoothness.  

Our notation established, we may now state our first theorem for stratified capillary-gravity waves.

\begin{theorem} \emph{(Global Continuum from Simple Bifurcation)} 
Fix a wave speed $c > 0$, wavelength $L > 0$, relative mass flux $p_0 < 0$ and coefficient of surface tension $\sigma > 0$.  Fix any $\alpha \in (0,1)$, and let the functions $\beta \in C^{3+\alpha}([0,|p_0|])$ and $\rho \in C^{1+\alpha}([p_0,0]; \mathbb{R}^+)$ be given such that the \eqref{lb} condition holds (see \textsc{Definition \ref{capgravlbc}}).   Also we assume the streamline density function $\rho$ is nonincreasing.  There exists a set $\Sigma_1 \subset (0,\infty)$ such that the following is true.  $\Sigma_1$ contains a half-line, $(\sigma_c, \infty)$, and, if $\sigma \in \Sigma_1$, then we have:  \\

\noindent
\emph{(a)}
There exists a connected set $\mathcal{C} \subset \mathscr{S}$ of solutions $(u,v,\rho, \eta,Q)$ of the traveling, stratified, water wave problem with surface tension \eqref{incompress}--\eqref{bedcond}.  The set $\mathcal{C}$, moreover, has the following properties.
\begin{itemize}
\item[$\bullet$] $\mathcal{C}$ contains a laminar flow (with a flat surface $\eta \equiv 0$ and all streamlines parallel to the bed); denote it by $(U_*, Q_*)$.
\item[$\bullet$] Either \emph{(i)} $\mathcal{C}$ contains more than one laminar flow, or it is unbounded in the sense that, along some sequence $\{(u_n, v_n, \rho_n, \eta_n)\} \subset \mathcal{C}$, we have 
\[ \mathrm{(ii)~} \max u_n \uparrow c \qquad \textrm{or} \qquad \mathrm{(iii)~} 
\min u_n \downarrow -\infty.\]
\end{itemize} 

\noindent
\emph{(b)} Furthermore, there exists a path-connected subset $\mathcal{K} \subset \mathcal{C}$ such that 
\begin{itemize}
\item[$\bullet$] $\mathcal{K}$ admits a global, locally injective, continuous parameterization with a locally injective $C^1$ reparameterization.    
\item[$\bullet$] Either $\mathcal{K}$ is a closed loop, or it unbounded in the sense that alternative \emph{(a)(ii)} or \emph{(a)(iii)} must occur along $\mathcal{K}$.  
\end{itemize}
\emph{(c)} All flows $(u,v,\rho,\eta) \in \mathcal{C}$ are regular in the sense that:
\begin{itemize}
\item[$\bullet$] $u,~v,~\rho$ and $\eta$ each have period $L$ in $x$;
\item[$\bullet$] $u$, $\rho$ and $\eta$ are symmetric, $v$ antisymmetric across the line $x =0$.
\end{itemize}

\label{mainresult1} \end{theorem}
\begin{longremark} The Local Bifurcation Condition \eqref{lb} essentially states that there exists an eigenvalue of the linearized problem with lowest mode, and is both necessary and sufficient for our result to hold.  Still, in \textsc{Lemma \ref{suffsizecondlemma}}, we derive an explicit sufficient condition that implies \eqref{lb}:
\be \begin{split}
 (g\rho(0)+\sigma)p_0^2 & > \int_{p_0}^0 \bigg\{ \frac{4\pi^2}{L^2} \left( 2B(p)-2B_{\mathrm{min}}+2\epsilon_0\right)^{3/2} \\
 & \qquad + (p-p_0)^2 \left( \left(2B(p)-2B_{\mathrm{min}}+2\epsilon_0\right)^{1/2} + g\rho^\prime(p) \right) \bigg\} dp, \end{split} \label{sizecond} \ee
where $\epsilon_0$ is defined by
\be \epsilon_0^{3/2} := \max\left\{2g \|\rho^\prime \|_{L^\infty} p_0^2 e^{|p_0|},~(2g\|\rho^\prime\|_{L^\infty})^{3},~(4\|\rho^\prime\|_{L^\infty})^3,~8g|p_0|\rho(0) \right\}. \label{defepsilon0}\ee
This is most certainly not optimal.  Indeed, obtaining a complete characterization of the set of data $(p_0, \rho, \beta)$ for which \eqref{lb} holds is an important and largely open problem, even in the special case when $\rho \equiv 1$ and $\sigma = 0$.  It is worth mentioning that \eqref{sizecond} can be satisfied by assuming $|p_0|$ and $\| \rho^\prime\|_{L^\infty}$ are sufficiently small, or $\sigma$ is large.  From this it is easy to see that the set $\Sigma := \{ \sigma > 0 : \eqref{lb} \textrm{ holds}\}$ contains a neighborhood of $+\infty$.  

We note also that much is known about the structure of the continuum $\mathcal{C}$ near the point of bifurcation $(U_*, Q_*)$.  There exists a neighborhood $\mathcal{U}$ of $(U_*,Q_*)$ in $\mathscr{S}$ in which $\mathcal{C} \cap \mathcal{U} = \mathcal{K} \cap \mathcal{U}$ comprises all solutions  to \eqref{incompress}--\eqref{bedcond}, and $(U_*,Q_*)$ is the only laminar solution in $\mathcal{U}$.
Moreover, all nonlaminar flows in $\mathcal{U} \cap \mathcal{C}$ exhibit the following monotonicity properties:
\begin{itemize}
\item[$\bullet$] within each period the wave profile $\eta$ has a single crest and trough; say the crest occurs at $x =0$; 
\item[$\bullet$] a water particle located at $(x,y)$ with $0 < x < L/2$ and $y > -d$ has positive vertical velocity $v > 0$.
\end{itemize}
For the proof of these results, see \textsc{Theorem \ref{capgravsimplelocalbifurcation}}.

Finally, let us note that in certain special cases it can be shown that $\Sigma_1$ is a countable set, see \textsc{Proposition \ref{genericproposition}}.  
\end{longremark}
As we shall see, \textsc{Theorem \ref{mainresult1}} is obtained by means of a bifurcation argument with the family of laminar flows taking on the role of the trivial solutions.  When $\sigma \in \Sigma_1$ (for example, if $\sigma$ is sufficiently large) we will show that local bifurcation occurs at a simple (generalized) eigenvalue.    

On the other hand, for $\sigma \notin \Sigma_1$, the eigenvalue will be not be simple.  In order to discuss this scenario more carefully, we must first make a few definitions.  Let $D$ be a smooth domain in $\mathbb{R}^2$ and let Banach spaces $X,Y \subset C_{\textrm{per}}^{m+\alpha}(\overline{D})$ be given for some $m \geq 1$, $\alpha \in (0,1)$.  For a smooth map $\mathcal{F}: \mathbb{R}\times X \to Y$  with $\mathcal{F}(\cdot, 0) \equiv 0$, the generic bifurcation problem is to find all nontrivial solutions $(\lambda, w)$ of  
\be \mathcal{F}(\lambda, w) = 0, \label{genericbifurcation} \ee
in a neighborhood of a given solution $(\lambda_*, 0)$.  Suppose that the null-space of $\mathcal{F}_w(\lambda_*,0)$ is two-dimensional, spanned by $\{ \phi_1, \phi_2\}$, and $\phi_i = \phi_i(x,y)$ is of minimal period $L/n_i$ in the first coordinate for $i = 1,2$, where $n_1$ and $n_2$ are distinct natural numbers.  Assume further that there exists a (continuously) parameterized curve of solutions, 
\[ \mathcal{C}_{\mathrm{loc}} := \{(\lambda(s), w(s)) : s \in (-\epsilon,\epsilon) \} \subset \mathbb{R} \times X,\]
to \eqref{genericbifurcation} with $(\lambda(0), w(0)) = (\lambda_*, 0)$.  We shall refer to these solutions as \emph{pure} provided that
\[ w(s) = s \phi_i + o(s) \qquad \textrm{in } X,\]
for either $i = 1$ or $i = 2$ but not both.  That is, pure solutions are those that locally lie near the interior of the span of either $\phi_1$ or $\phi_2$.  Conversely, the elements of $\mathcal{C}_{\textrm{loc}}$ are called \emph{mixed solutions} if there exists $\xi_1, \xi_2 \in \mathbb{R}^\times$ such that
\[ w(s) = s\left( \xi_1 \phi_1 + \xi_2 \phi_2 \right) + o(s) \qquad \textrm{in } X.\]
Clearly, mixed solutions locally lie near the interior of the span of $\phi_1$ and $\phi_2$.  

The relevance of these concepts to the case at hand stems from the fact that, for $\sigma \in \mathbb{R}^+ \setminus \Sigma_1$, local bifurcation occurs at eigenvalues of multiplicity two (see \textsc{Corollary \ref{Sigma12corollary}}.) That is, the null space of the linearized operator is spanned by two functions, $\{\phi_{1}, \phi_{2}\}$, of minimal period $L/n_1, L/n_2$ in $x$, respectively.   This feature of the eigenvalue problem is well-known for homogeneous capillary-gravity waves and has been very well analyzed in the irrotational case (cf. \cite{jones1986symmetry,okamoto1990on,toland1985bifurcation} for analytic results and, e.g. \cite{okamoto1996resonance,schwartz1979numerical,shoji1989new} and the references contained therein for some numerical results.)  For general vorticity, this scenario was analyzed by Wahl\'en \cite{wahlen2006capgrav}.  Briefly, his argument was the following.  So long as there is no resonance, i.e. $n_2$ is not an integer multiple of $n_1$, then the dimension of the null space can be reduced to one by considering the restriction to the subspaces of $L/n_1$- or $L/n_2$-periodic solutions.  Applying the Crandall-Rabinowitz arguments developed for the simple eigenvalue case, one can prove the existence of two $C^1$-curves bifurcating from the family of laminar flows, say $\mathcal{C}_{1,\textrm{loc}}$ and $\mathcal{C}_{2,\textrm{loc}}$.  These consist of small amplitude solutions of minimal period $L/n_1$ and $L/n_2$, respectively.  In fact, $\mathcal{C}_{i,\textrm{loc}}$ are both curves of pure solutions in the sense outlined above.  

The price payed for the relative ease of this approach is that, in restricting the solutions space, we are likely missing the mixed solutions.  To see why this may be the case, note that a linear combination $\xi_1 \phi_1 + \xi_2 \phi_2$ may have \emph{minimal} period $L$.  Such solutions would therefore not be included in either $\mathcal{C}_{1,\textrm{loc}}$ or $\mathcal{C}_{2,\textrm{loc}}$.  

Worse still, when resonance does occur (i.e. $n_2/n_1 \in \mathbb{N}$), this type of argument offers even less information.  Since restriction to the $L/n_1$-periodic solutions no longer renders the null space one-dimensional (as $L/n_2$--periodic solutions are now also $L/n_1$--periodic), the only conclusion we can draw via Crandall-Rabinowitz is that there is at least one local curve of pure $L/n_2$-minimal periodic solutions.  In the language we have been using thus far, this means that we can prove the existence of $\mathcal{C}_{2,\textrm{loc}}$, but not $\mathcal{C}_{1,\textrm{loc}}$.  This, of course, says nothing about the existence (or nonexistence) of curves of mixed solutions.

In the present paper we devote some effort to advancing the study of the bifurcation near double eigenvalues.  Rather than using Crandall-Rabinowitz, we instead carry out our analysis via the full Lyapunov-Schmidt procedure.   Owing to a great deal of degeneracy in the problem, a complete account is, unfortunately, beyond the scope of this paper.  Nonetheless, we are able to obtain some detailed results on the structure of the full solution set. 

\begin{theorem}  \emph{(Local Continua from Double Bifurcation)}
Let the wave speed $c >0$, the relative mass flux $p_0 < 0$, surface tension constant $\sigma$, streamline density function $\rho \in C^{2+\alpha}([0,|p_0|])$ and Bernoulli function $\beta \in C^{1+\alpha}([0,|p_0|])$, $0 < \alpha < 1$ be given. Suppose further that $p_0$, $\rho$, $\beta$ and $\sigma$ collectively satisfy \eqref{lb}.  

In contrast to \textsc{Theorem \ref{mainresult1}}, suppose $\sigma \notin \Sigma_1$ but non-degeneracy conditions \eqref{capgravnondegeneracycondition1} and \eqref{capgravnondegeneracycondition2} hold.  Then there exists four $C^1$-curves, $\{ \mathcal{C}_{\mathrm{i,loc}}\}_{i=1}^4$, of small-amplitude, traveling wave solutions of \eqref{incompress}-\eqref{bedcond} in the space $\mathscr{S}$ with the following properties.   
\begin{itemize}
\item[$\bullet$] Each $(u,v,\rho,\eta,Q) \in \bigcup_i \mathcal{C}_{\mathrm{i,loc}}$ has speed $c$, relative mass flux $p_0$ and satisfies $u < c$ throughout the fluid.
\item[$\bullet$] For $i=1,2$, each $(u,v,\rho, \eta,Q) \in \mathcal{C}_{\mathrm{i,loc}}$, $i = 1,2$, has minimal period $2\pi/n_i$ in $x$.
\item[$\bullet$] For $i=1,2,3,4$, each of the curves $\{\mathcal{C}_{\mathrm{i,loc}}\}$ contain precisely one laminar flow, while for each nontrivial solution $(u, v, \rho,\eta,Q) \in \mathcal{C}_{\textrm{i,loc}}$, with $i =1, 2$, we have \begin{enumerate}
\item[\emph{(i)}] the functions $u$ and $\eta$ are symmetric around the line $x = 0$ while $v$ is antisymmetric; 
\item[\emph{(ii)}] the function $\eta$ has precisely one maximum (crest) and one minimum (trough) per minimal period; and 
\item[\emph{(iii)}] the wave profile is strictly monotone between crest and trough.
\end{enumerate}
\end{itemize}
\label{mainresult2}
\end{theorem}
\begin{longremark}  Referring back to the preceding discussion, the first two curves, $\mathcal{C}_{\mathrm{1,loc}}$ and $\mathcal{C}_{\mathrm{2,loc}}$ are the pure solutions, while the remaining two are the mixed solutions.  

The nondegeneracy conditions \eqref{capgravnondegeneracycondition1}--\eqref{capgravnondegeneracycondition2} ensure that, upon applying the Lyapunov-Schmidt reduction, the solutions of the bifurcation equation are described totally by the expansion to cubic order.  It is possible to omit condition \eqref{capgravnondegeneracycondition2} but, in this case, rather than having the complete solution set near $(U_*,Q_*)$, we can only guarantee the existence of $\mathcal{K}_1$.  Condition \eqref{capgravnondegeneracycondition1} can be similarly relaxed (see the remarks in section \ref{doubleeigenvaluebifurcationsection}).  More detailed information on the local structure of the bifurcation diagram is outlined in \textsc{Lemma \ref{lineartermslemma}}--\textsc{Lemma \ref{regularvaluelemma}} and \textsc{Theorem \ref{bifurcationstructuretheorem}}.
\end{longremark}

Once the local bifurcation picture has been resolved, we are able to continue the curves $\{\mathcal{C}_{i,\textrm{loc}}\}_{i=1}^4$ to produce global continua of solutions.  All told, this gives the last of our main results:

\begin{theorem} \emph{(Global Continua from Double Bifurcation)} 
Assume the hypotheses of \textsc{Theorem \ref{mainresult2}}.  There exist global continua of solutions extending $\{\mathcal{C}_{i,\mathrm{loc}}\}_{i=1}^4$  in the sense that the following statements hold.  \\

\noindent
\emph{(a)} There exists a path-connected set $\mathcal{K} \subset \mathscr{S}$ of solutions $(u,v,\rho, \eta, Q)$ to \eqref{incompress}--\eqref{bedcond}.  The set $\mathcal{K}$ has the following structural properties.
\begin{itemize}
\item[$\bullet$] $\mathcal{K}$ contains a laminar flow, $(U_*, Q_*)$. 
\item[$\bullet$] There exists four path-connected sets $\{\mathcal{K}_i\}_{i=1}^4$ such that $\mathcal{K} = \bigcup_{i} \mathcal{K}_i$ and $(U_*, Q_*) \in \bigcap_{i} \mathcal{K}_i$.  In fact,  there exists a neighborhood $\mathcal{U}$ of $(U_*, Q_*)$ in $\mathscr{S}$ such that
\[ \mathcal{C}_{i, \mathrm{loc}} = \mathcal{K}_i \cap \mathcal{U}, \qquad \textrm{for i = 1, \ldots, 4.}\]
\item[$\bullet$] Each $\mathcal{K}_i$ admits a global, locally injective, continuous parameterization with a locally injective $C^1$-reparameterization.  
\item[$\bullet$] For each $i = 1,\ldots, 4$, we have the following alternatives:  Either \emph{(i)} $\mathcal{K}_i$ is a closed loop or, along some sequence $\{(u_n, v_n, \rho_n, \eta_n, Q_n)\} \subset \mathcal{K}_i$,
\[ \mathrm{(ii)~} \max u_n \uparrow c \qquad  \textrm{or} \qquad \mathrm{(iii)~}\min u_n \downarrow -\infty. \]
\end{itemize}

\noindent
\emph{(b)} All flows $(u,v,\rho,\eta) \in \mathcal{K}$ are regular in the sense that:
\begin{itemize}
\item[$\bullet$] $u,~v,~\rho$ and $\eta$ each have period $L$ in $x$;
\item[$\bullet$] $u$, $\rho$ and $\eta$ are symmetric, $v$ antisymmetric across the line $x =0$.
\end{itemize}

\label{mainresult3}
\end{theorem}

Steady capillary-gravity waves, as we remarked earlier, are an important, classical object of study in fluid mechanics and geophysical fluid dynamics.  Unsurprisingly, therefore, they have been the subject of a truly voluminous quantity of research, far more than we can hope to completely summarize here.    The vast majority of the work on steady, capillary-gravity waves considers the irrotational regime.  There are physical justifications for this choice, but perhaps the most compelling reason is mathematical:  when there is no vorticity, $\psi$ is a harmonic function.  All that remains is to determine the free surface $\eta$.  As the vector field has zero curl, moreover, we have the existence of a velocity potential $\phi$, which enables the use of conformal maps and other powerful techniques from complex analysis.  For instance, there are a number of equivalent reformulations of \eqref{incompress}--\eqref{euler2} as a scalar equation for $\theta$, the slope of $\eta$ (cf. \cite{levicivita1925determination,nekrasov1951exact} for the classical case where $\sigma = 0$, and, e.g. \cite{okamoto1990on} for the generalization to $\sigma > 0$.)  

Of course, with the introduction of vorticity, the usefulness of these tools is greatly diminished.  For that reason, the structure of our arguments in this paper follows the general program for steady, rotational waves established in \cite{constantin2004exact,wahlen2006capillary,wahlen2006capgrav,walsh2009stratified}.  We begin, in section \ref{reformulationsection}, by reformulating the problem in terms of scalar quantities and in a semi-Lagrangian coordinate system.   As a consequence, the transformed domain will be fixed, and the system of equations \eqref{incompress}-\eqref{bedcond} is replaced by a non-linear, integro-differential, boundary value problem for a scalar function, $h$, that describes the height above the flat bed.   The interior equation satisfied by $h$ has the following form:
\[ \mathcal{L}(D^2h, Dh) + \mathcal{I}(Dh, h) = 0,\] 
where $D^2 h$ denotes the Hessian of $h$, $Dh$ the gradient, $\mathcal{L}$ is an elliptic operator and $\mathcal{I}$ is an integro-differential operator.  As in \cite{walsh2009stratified}, we will show that $\mathcal{I}$ totally accounts for the presence of stratification (indeed, if $\rho$ is constant, then $\mathcal{I} \equiv 0$).  On the other hand, the technical challenges due to surface tension present themselves, naturally, on the boundary.  The reformulation of the dynamic condition (in particular, the mean curvature term in \eqref{presscond}) results in a non-degenerate Venttsel-type condition on the image of the free surface.  

In section \ref{capgravlocalbifurcationsection}, we prove the existence of a 1-parameter family, $\mathcal{T}$, of laminar flows.  Since these flows are, by definition, functions of the streamlines, the integro-differential equation reduces to an ODE with operator coefficients (the operator coming from $\mathcal{I}$).   We proceed to show that, if $\sigma \in \Sigma_1$ --- e.g., it is sufficiently large --- then  there exists a point on $\mathcal{T}$ at which the problem has a simple (generalized) eigenvalue.  Using the theory of Crandall-Rabinowitz, we prove that there is a branch of small-amplitude, non-laminar solutions bifurcating from $\mathcal{T}$ at that point.  This is, essentially, a generalization of \cite{wahlen2006capgrav} to the stratified case.  

Next, in section \ref{doubleeigenvaluebifurcationsection}, we consider the case when $\sigma \notin \Sigma_1$.  As mentioned above, this means that the bifurcation from $\mathcal{T}$ is not simple.  In fact, we show that, for all but at most one value of $\sigma$, the bifurcation is precisely double, i.e. the null space of the linearized operator has algebraic multiplicity two.  Using the classical method of Lyapunov-Schmidt, we give several results on the bifurcation diagram for such $\sigma$.  This is done by explicit computation of the coefficients in the bifurcation equation.  In particular, we show that, under the non-degeneracy condition, there are four $C^1$-curves of small-amplitude, non-laminar solutions emanating from the bifurcation point, as described in \textsc{Theorem \ref{mainresult2}}.  Let us note that this result, though preliminary in the sense that it requires some additional non-degeneracy assumption, is new even in the case of constant density.  

The remainder of the paper is devoted to the global analysis.  In section \ref{globalbifurcationsection}, we use a degree theoretic argument to prove an alternative theorem in the style of Rabinowitz \cite{rabinowitz1971some}.  This allows us to continue globally the curve found via simple local bifurcation in section \ref{capgravlocalbifurcationsection}.  Of course, this involves first proving some strong compactness properties about the reformulated problem --- a task made difficult both by the nonlocal operator term in the interior equation, and the nonstandard, second-order boundary conditions.  

To continue the local curves obtained in the case of double bifurcation, we employ a global bifurcation theory due to Dancer (cf. \cite{buffoni2003analytic,dancer1973bifurcation,dancer1973globalsolution,dancer1973globalstructure}) that exploits the analyticity of the integro-differential operator in the reformulated problem.  As we make clear, this method enjoys several advantages with respect to the more standard degree theoretic approach of section \ref{globalbifurcationsection}.  First, it enables us to prove an alternative theorem that continues each of the local curves found in section \ref{doubleeigenvaluebifurcationsection}. Since the crossing at the bifurcation point is no longer odd, this would be non-trivial with a degree theoretic argument.  Second, we are able to prove higher regularity of the global continuum.  In fact, we show that the continuum $\mathcal{K}$ of solutions to the reformulated problem is a curve with a locally injective, \emph{analytic} reparameterization.  

 In section \ref{uniformregularitysection}, we derive uniform bounds on the $C^{3+\alpha}$-norm along the continuum $\mathcal{C}$ in terms of the first-order derivatives, the mean curvature on the surface and the energy $Q$.  In section \ref{finalsection}, these bounds, along with the alternative theorems of sections \ref{globalbifurcationsection} and \ref{analyticitysection}, are combined to produce \textsc{Theorem \ref{mainresult1}} and \textsc{Theorem \ref{mainresult3}}.  Finally, section \ref{examplessection} analyzes some important special cases:  homogeneous, irrotational waves and homogeneous rotational waves.
 
\section{Reformulation of the problem} \label{reformulationsection}

In this section we set ourselves to the task of finding an equivalent formulation of \eqref{euler2}--\eqref{boundcond} that eliminates the free boundary problem.   We shall accomplish this by making a change of variables to transform the fluid domain $\overline{D_\eta}$ in a fixed domain, whose closure we shall denote $\overline{R}$.  As usual for incompressible problems, this will be done by means of the pseudo-stream function $\psi$ and will come at the cost of some additional nonlinearity in the governing equation.  

To see why this is a natural choice of coordinates, recall that we have $\psi \equiv 0$ on the free surface, $\psi \equiv -p_0$ on the flat bed, where $p_0$ is the relative pseudo-mass flux and $\psi$ is defined uniquely by the requirement that
\be \psi_x = -\sqrt{\rho} v, \qquad \psi_y = \sqrt{\rho}(u-c). \label{defpsi} \ee
In light of \eqref{defpsi} and \eqref{incompress}-\eqref{consmass1}, the governing equations inside the fluid become
\be \left\{ \begin{array}{lll}
\psi_y \psi_{xy} - \psi_x \psi_{yy} & = & -P_x \\
-\psi_y \psi_{xx} + \psi_x\psi_{xy} & = & -P_y - \rho g \end{array} \right. \qquad \textrm{in } \overline{D_\eta}, \label{equivform} \ee
whereas the boundary conditions \eqref{boundcond} are
\be \left\{ \begin{array}{llll}
\psi_x & = & -\psi_y \eta_x, & \qquad \textrm{on } y = \eta(x), \\
P & = & P_{\textrm{atm}} - \sigma {\eta_{xx}}{\left(1+\eta_x^2\right)^{-3/2}}, & \qquad \textrm{on } y = \eta(x), \\
\psi_x & = & 0, & \qquad \textrm{on } y = -d. \end{array} \right. \label{equivboundcond} \ee

Also, by Bernoulli's theorem, the quantity 
\be E = P + \frac{\rho}{2}\left( (u-c)^2 + v^2\right) + g\rho y, \label{defE1} \ee
is constant along streamlines.  In particular, then, evaluating the above equation on the free surface we find
\be |\nabla \psi|^2 + 2g \rho(x, \eta(x)) \left(\eta(x)+d\right) = Q, \qquad \textrm{on } y = \eta(x), \label{defQ} \ee
where the constant $Q := 2(E|_\eta-P_{\textrm{atm}} + g\rho|_\eta d)$.  Note that this $Q$ gives roughly the energy density along the free surface of the fluid.  Critically, however, we shall see that as a consequence of our normalization of $\eta$, $d$ is not a parameter for the problem.  On the contrary, in all but the most trivial cases, $d$ varies along $\mathcal{C}$.  In our analysis, therefore, we shall instead be viewing $Q$ as parameterizing the continuum.

With that in mind, we begin by deriving a scalar PDE satisfied by $\psi$.  In the case where the fluid is homogeneous, it is not hard to check that
\be -\Delta \psi = \omega, \label{homgoverneq} \ee
where $\omega := v_x - u_y$ is the vorticity.  Moreover, again exploiting the homogeneity, one can readily verify that $\omega$ is transported, i.e. there exists some \emph{vorticity function}, $\gamma:[0,|p_0|] \to \mathbb{R}$, such that $\omega = \gamma(\psi)$.  It is therefore possible to view \eqref{homgoverneq} as a semilinear elliptic equation, where $\gamma$ is prescribed.  This is the approach taken, for instance, by Constantin and Strauss (\cite{constantin2004exact,constantin2007rotational}) in the (homogeneous) gravity wave case, and Wahl\'en (\cite{wahlen2006capillary,wahlen2006capgrav}) for (homogeneous) capillary and (homogeneous) capillary-gravity waves.
    
In the presence of stratification, however, we cannot generally expect the vorticity to be a function of $\psi$.  Instead, we must look back to the Bernoulli theorem to find an appropriate quantity conserved along streamlines.  In particular, one can prove that
\[ \frac{d E}{d \psi} = \Delta \psi + gy \frac{d \rho}{d \psi}. \]
Recall that in the previous section we introduced the Bernoulli function $\beta$ and streamline density function $\rho$.  Rewriting the above expression we arrive at an equation of roughly the same form as \eqref{homgoverneq}:
\be -\beta(\psi) = \Delta \psi - gy \rho^\prime(-\psi). \label{yiheq} \ee
This is known as Yih's equation or the Yih--Long equation (see, e.g., \cite{yih1965dynamics}).  In this paper, it shall serve as our governing equation for $\psi$.  Were it not for the presence of $y$ on the right-hand side, this would reduce to the homogeneous case --- a fact that agrees with our intuitive notion that density stratification should reintroduce the depth into the problem as a serious consideration. One can easily check, moreover, that, if $\rho \equiv 1$, then $\beta$ and $\gamma$ coincide exactly. 

With \eqref{yiheq} in hand, we now make a change of variables to eliminate the free boundary.  The new coordinates we will denote $(q,p)$ where 
\[ q = x, \qquad p = -\psi(x,y).\]
This scheme is sometimes referred to as \emph{semi-Lagrangian} coordinates, in recognition of the fact that we are working, in some sense, halfway between Lagrangian (or streamline) coordinates and the usual Eulerian system (cf. \cite{turner2002traveling}).  

By means of a scaling argument, we may take $L := 2\pi$.  Then, under the transformation
\[ (x,y) \mapsto (q,p),\]
the closed fluid domain $\overline{D_\eta}$ is mapped to the rectangle
\[ \overline{R} := \{ (q,p) \in \mathbb{R}^2 : 0 \leq q \leq 2\pi,~p_0 \leq p \leq 0 \}.\]
Given this, it will be convenient to put
\[ T := \{ (q,p) \in R : p = 0 \}, \qquad B := \{ (q,p) \in R : p = p_0 \}. \]
  
Note that, in light of \eqref{defrho} and \eqref{defbeta}, we have that
\[ \beta = \beta(-p), \qquad \rho = \rho(p).\]
Moreover, the assumption that the streamline density function is nonincreasing becomes
\be \rho_p \leq 0. \label{rhoincdepth} \ee
  
Next, following the ideas of Dubreil--Jacotin (cf. \cite{dubreil1934determination,dubreil1937theoremes}), define 
\be h(q,p) := y + d \label{defh} \ee
which gives the height above the flat bottom on the streamline corresponding to $p$ and at $x = q$.  We calculate:
\be \psi_y = -\frac{1}{h_p}, \qquad \psi_x = \frac{h_q}{h_p} \label{hqhpeq}. \ee
Note that this implies $h_p > 0$, because we have stipulated that $u < c$ throughout the fluid.  The change of variables then gives
\be u  =  c - \frac{1}{\sqrt{\rho} h_p},~\qquad v = -\frac{h_q}{\sqrt{\rho} h_p}. \label{uvhphq} \ee

Performing the Dubreil--Jacotin transformation on \eqref{yiheq}, we find that, in the interior of the fluid, 
\[ (1+h_q^2)h_{pp} + h_{qq}h_p^2 - 2h_q h_p h_{pq}-g(h-d(h))\rho_p h_p^3  = -h_p^3 \beta(-p), \]
where we have used \eqref{defh} to write $y = h-d$.  Recall, however, that we have normalized $\eta$ so that it has mean zero.  Taking the mean of \eqref{defh} along $T$, we obtain
\be d = d(h) = \fint_0^{2\pi} h(q,0)dq. \label{defd} \ee
That is, the average depth $d$ must be viewed as a linear operator acting on $h$.  Namely, it is the average value of $h$ over $T$.  Where there is no risk of confusion, we shall suppress this dependency and simply write $d$.  As we shall see, the addition of this integral term to the governing equations will be the single most significant departure from the homogeneous case, both technically and qualitatively.  In recognition of this fact, we shall refer to the mapping $h \mapsto -g(h-d(h))h_p^3 \rho_p$ as the \emph{stratification operator}.  See \cite{walsh2009stratified} for details.

On the bottom, we have by construction that $h \equiv 0$.  The more interesting boundary condition is found on the top of $R$, which corresponds to the free surface in the original formulation. There, we have from \eqref{equivboundcond} that
\[ 1+h_q^2 + h_p^2( 2\sigma\kappa[h]+2g\rho h- Q)  = 0, \qquad p = 0,\]
where the curvature, $\kappa$, is defined by
\[ \kappa[h] := - \frac{ h_{qq}}{\left(1+h_q^2\right)^{3/2}}.\]
This calculation is carried out in \cite{wahlen2006capillary}, for example.  

Combining these observations together we see that the completely reformulated problem is the following.  For a given Bernoulli function $\beta$, find  $(h,Q) \in C^{3+\alpha}_{\textrm{per}}(\overline{R}) \times \mathbb{R}$ with $h_p > 0$ and satisfying the height equation
\be \left \{ \begin{array}{lll}
(1+h_q^2)h_{pp} + h_{qq}h_p^2 - 2h_q h_p h_{pq}  
\\ \qquad -g(h-d(h))\rho_p h_p^3  = -h_p^3 \beta(-p), & p_0 < p < 0, \\
& & \\
1+h_q^2 + h_p^2( 2\sigma\kappa[h]+2g\rho h- Q)  = 0, & p = 0, \\
h = 0, & p = p_0. \end{array} \right. \label{heighteq} \ee

\begin{lemma} \emph{(Equivalency)} Problem \eqref{heighteq} is equivalent to problem \eqref{incompress}-\eqref{bedcond}
\label{equivalencelemma} \end{lemma}

\begin{proof2} This is proven for $\sigma = 0$ in \cite{walsh2009stratified}.  A straightforward modification of that argument --- which is itself a simple generalization of a similar result in \cite{constantin2004exact} --- proves the lemma.  For brevity, we omit the details here.  The reader is also referred to the discussion in Remark \ref{eulernotanalyticremark}. \end{proof2}

\section{Local bifurcation at simple eigenvalues } \label{capgravlocalbifurcationsection}

\subsection{Overview}  \label{overviewlocalbifurcationsection}

The goal of this section is to prove the existence of local curves of small amplitude solutions to \eqref{euler2}--\eqref{boundcond}.  In continuation from the previous section, we shall make use of the equivalent formulation \eqref{heighteq} to prove the existence of a 1-parameter family of laminar flow solutions, $\mathcal{T}$.  These will serve as the trivial solution curve for the standard Crandall-Rabinowitz machinery.  Linearizing the problem around the laminar flows, we are thereby able to prove the existence of a local curve of non-laminar solutions bifurcating from $\mathcal{T}$ at a simple eigenvalue under certain size conditions.  The analysis culminates in the following result.

\begin{theorem} \emph{(Local Bifurcation at Simple Eigenvalues)} 
Fix $\alpha \in (0,1)$ and let the wave speed $c >0$, the relative mass flux $p_0 < 0$, surface tension constant $\sigma$, streamline density function $\rho \in C^{2+\alpha}([p_0,0]; \mathbb{R}^+)$ and Bernoulli function $\beta \in C^{1+\alpha}([0,|p_0|])$ be given. Suppose further that $p_0$, $\rho$, $\beta$ and $\sigma$ collectively satisfy \eqref{lb}.  If $\sigma \in \Sigma_1$ (as defined in \eqref{defSigma1}), then the following holds.

There exists a $C^1$-curve $\mathcal{C}_{\textrm{loc}} \subset \mathscr{S}$ of small amplitude, traveling wave solutions $(u, v, \rho, \eta)$ of \eqref{incompress}-\eqref{bedcond} with period $2\pi$, speed $c$ and relative mass flux $p_0$, satisfying $u<c$ throughout the fluid. The curve $\mathcal{C}_{\mathrm{loc}}$ has the following properties.  
\begin{enumerate}
\item[\emph{(i)}] $\mathcal{C}_{\textrm{loc}}$ contains precisely one laminar flow. denote it $(Q_*, U_*)$.
\item[\emph{(ii)}] $\mathcal{C}_{\textrm{loc}} \setminus \{(Q_*, U_*)\}$ comprises all nonlaminar solutions in a sufficiently small neighborhood of $(Q_*, U_*)$ in $\mathbb{R} \times C^{2+\alpha}_{\mathrm{per}}(\overline{D_\eta})$.  
\item[\emph{(iii)}] For each nontrivial solution $(u, v, \rho, \eta) \in \mathcal{C}_{\textrm{loc}}$, 
\begin{itemize}
\item[$\bullet$] the functions $u$, $\rho$ and $\eta$ are symmetric around the line $x = 0$ while $v$ is antisymmetric; 
\item[$\bullet$] the function $\eta$ has precisely one maximum (crest) and one minimum (trough) per period; 
\item[$\bullet$] the wave profile is strictly monotone between crest and trough.
\end{itemize}
\end{enumerate} \label{capgravsimplelocalbifurcation}
\end{theorem}

\subsection{Eigenvalue problem} \label{eigenvalueproblemsection} 

We begin by establishing the existence of a 1-parameter family of laminar flows.  By laminar flow we mean that the streamlines (which, we recall, include the free surface) are parallel to the bed.  Waves of this type will accordingly have the ansatz $h = H(p)$, which reduces the height equation \eqref{heighteq} to the following (operator coefficient) ODE boundary value problem:
\be \left\{ \begin{array}{lll}
H_{pp} -gH_p^3(H-d(H)) \rho_p  =  - H_p^3 \beta(-p), & p_0  < p < 0, \\
1+H_p^2(2g\rho H -Q) = 0, & p = 0, \\
H = 0, & p = p_0. \end{array} \right. \label{laminar2} \ee

The basic existence theory for laminar solutions comes in the form of the following lemma.

\begin{lemma}
\emph{(Laminar Flow)} 
Suppose that the streamline density function $\rho$ satisfies $\rho_p \leq 0$.  Then there exists a 1-parameter family of solutions $(H(\cdot;\lambda), Q(\lambda))$ to the laminar flow equation \eqref{laminar2} with $H_p > 0$, where $0 \leq  -2B_{\min{}} < \lambda < Q$. \label{laminarflowlemma} \end{lemma} 

\begin{proof2}  Because the curvature of the laminar flows is zero, \eqref{laminar2} is identical to the problem satisfied by laminar stratified gravity waves in the absence of surface tension.  The lemma therefore follows from \textsc{Lemma 3.2} of \cite{walsh2009stratified}. \end{proof2}

\begin{longremark} For later reference, we note it was also shown in \cite{walsh2009stratified} that
\[ H(0) = \frac{Q-\lambda}{2g\rho(0)}, \qquad H_p(0) = \lambda^{-1/2}.\]

In contrast to the theory for homogeneous gravity and capillary-gravity waves, there is no explicit formula for these laminar flows.  However, $H_p$ satisfies the following implicit equation 
\[ H_p(p;\lambda) = \frac{1}{\sqrt{\lambda+G(p;\lambda)}}, \qquad p_0 < p < 0, \]
where
\[ G(p;\lambda) := 2B(p) + 2\int_p^0 g\left(H(r;\lambda)-d(H(r;\lambda))\right)\rho_p(r) dr, \qquad p_0 < p < 0.\]
It will occasionally be convenient to make use of the shorthand 
\[ a(\cdot;\lambda) := H_p(\cdot;\lambda)^{-1},\qquad Y(\cdot; \lambda) := H(\cdot;\lambda) -d(H).\]  
Finally, we shall denote $\dot{H} := \partial H/ \partial \lambda$, $\dot{Y} := \partial Y/\partial \lambda$, $\dot{G} := \partial G/\partial\lambda$ and $\dot{Q} := dQ/d\lambda$.  In fact, we have the following useful identity
\be \dot{Y}(p) = \frac{1}{2} \int_{p_0}^0 \frac{1+\dot{G}(r)}{(\lambda+G(r))^{3/2}}dr = \frac{1}{2} \int_{p}^0 \left(1+\dot{G}(r)\right) a(r)^{-3} dr, \qquad p_0 \leq p \leq 0. \label{Ydotidentity} \ee
This is found by differentiating the relation $Y_p = (\lambda+G)^{1/2}$ in $\lambda$, followed by integrating in $p$.   \label{futurereferenceremark} \end{longremark}

We now linearize the full height equation \eqref{heighteq} along the curve of laminar solutions, which we shall denote $\mathcal{T}$.  As in \cite{walsh2009stratified}, when there is heterogeneity in the fluid, it will be necessary to ensure that $\lambda$ is bounded strictly away from $-2B_{\textrm{min}}$.  Recall that we have defined  
\[ \epsilon_0^{3/2} :=  \max\left\{2g \|\rho^\prime \|_{L^\infty} p_0^2 e^{|p_0|},~(2g\|\rho^\prime\|_{L^\infty})^{3},~(4\|\rho^\prime\|_{L^\infty})^3,~8g|p_0|\rho(0) \right\}.\] 
Note that having $\epsilon_0^{3/2}$ greater than the second quantity on the right-hand side above implies, for $\lambda > -2B_{\textrm{min}} +\epsilon_0$,  
\[a+g\rho_p \geq \epsilon_0^{1/2} - g\|\rho_p\|_{L^\infty} \geq \epsilon_0/2 > 0\]  
if $\rho_p \nequiv 0$.  The same holds true, of course, when $\rho_p \equiv 0$, since in this case the left-hand side reduces to $a$, which is strictly positive when $\lambda > -2B_{\textrm{min}}$.  It was shown in \cite{walsh2009stratified} that, when $\epsilon_0^{3/2} > 2g\|\rho_p\|_{L^\infty} p_0^2 e^{|p_0|}$, then $-1/2 \leq \dot{G}(\cdot; \lambda) \leq 0$.  In particular, this implies that, $a+g\rho_p > 0$, and $1+\dot{G} > 1/2 > 0$.  Both of these facts are critical for our analysis and so from this section onward we shall always be assuming that $\lambda > -2B_{\textrm{min}}+\epsilon_0.$  The remaining two quantities on the right-hand side above are, in fact, not strictly necessary.  In \cite{walsh2009stratified} they were used to avoid a certain pathology that might occur if the point of bifurcation $\lambda_*$ occurred to the right of a particular value $\lambda_0$.  That will not be an issue for the present paper.  Nonetheless, for simplicity of exposition, we keep the same definition of $\epsilon_0$.  

Now, fixing any such $\lambda$, we calculate that the linearized problem about $H(\cdot; \lambda)$ is: 
\be \left \{ \begin{array}{lll}
m_{pp} + m_{qq}H_p^2 -g \rho_p (3m_p (H-d(H)) H_p^2 + (m-d(m))H_p^3), & \\
\qquad =  -3H_p^2 m_p \beta(-p) & p_0<p<0, \\
m =  0, & p = p_0, \\
g\rho m  - \sigma m_{qq} - \lambda^{3/2} m_p = 0, & p = 0. \end{array} \right. \label{linearheighteq} \ee

Since we seek solutions that are $2\pi$-periodic and even in $q$, it is natural to search for $m$ with the ansatz $m(q,p) = M(p)\cos{(nq)}$, for some $n \geq 0$.  Inserting this into \eqref{linearheighteq} we see immediately that $M$ must satisfy the following:
\be \left \{ \begin{array}{lll}
m_{pp} -n^2 M H_p^2 -g \rho_p (3M_p (H-d(H)) H_p^2 + (M-\delta_{n0}M(0))H_p^3 ) & \\
\qquad =  -3H_p^2 M_p \beta(-p), & p_0<p<0, \\
M =  0, & p = p_0, \\
g\rho M  +n^2 \sigma M - \lambda^{3/2} M_p = 0 & p = 0, \end{array} \right. \label{linearheighteq2} \ee
where $\delta_{n0} = 1$ if $n = 0$ and vanishes otherwise.  This is a simple consequence of the fact that
\[ d(m) = \fint_0^{2\pi} m(q,0)dq = M(0) \fint_{0}^{2\pi} \cos{(nq)} dq = \delta_{n0} M(0).\]
Using the equation satisfied by $H$, \eqref{laminar2}, we may put \eqref{linearheighteq2} into self-adjoint form:
\be\left\{ \begin{array}{ll}
\{H_p^{-3} M_p \}_p = (n^2 H_p^{-1} +g\rho_p)M - g\rho_p \delta_{n0} M(0),  & p_0 < p < 0,\\ 
M  = 0, & p = p_0, \\
\lambda^{3/2} M_p = (n^2\sigma + g\rho)M, & p = 0. \end{array}\right. \label{sturmliouvilleode}  \ee
Solving the above problem for any value of $n$ will produce a $2\pi$-periodic solution to the linearized problem, but, as we are primarily interested in solutions of minimal period $2\pi$, we shall focus on the case where $n = 1$. 

With that in mind, we make the following definition.
\be \begin{split} \Sigma_1 &:= \{ \sigma \in \mathbb{R}^+ \colon \textrm{there exists a unique } \lambda > -2B_{\mathrm{min}} + \epsilon_0 \textrm{ such that} \\ 
& \qquad \textrm{problem \eqref{sturmliouvilleode} has a solution $(M, n)$ if and only if } n = 1  \}, \end{split} \label{defSigma1} \ee
\be \begin{split} \Sigma_2 &:= \{ \sigma \in \mathbb{R}^+ \colon \textrm{there exists a } \lambda > -2B_{\mathrm{min}} + \epsilon_0 \textrm{ such that} \\ 
& \qquad \textrm{problem \eqref{sturmliouvilleode} has precisely two linearly independent} \\
& \qquad \textrm{nontrivial solutions $(M_1, n_1)$ and $(M_2, n_2)$,} \\
& \qquad \textrm{where $n_1 = 1$, $n_2 \geq 2$} \}. \end{split} \label{defSigma2} \ee
and
\be \begin{split} \Sigma_3 & := \{ \sigma \in \mathbb{R}^+ : \textrm{there exists a $\lambda > -2B_{\mathrm{min}}+\epsilon_0$ such that} \\
& \qquad \textrm{problem \eqref{sturmliouvilleode} has at least two linearly independent solutions} \\
& \qquad \textrm{$(M_1,n_1)$ and $(M_2,n_2)$ where $n_1 = 1$, $n_2 = 0$}\} \end{split}\label{defSigma3} \ee

For fixed $p_0$, $\rho$ and $\beta$, we shall prove that the set $\Sigma_1$ consists of those values of $\sigma$ for which simple bifurcation occurs at some point along $\mathcal{T}$.  The curve of nonlaminar solutions will, locally, be of the form $(Q(\lambda), H(p; Q(\lambda)) + \epsilon M(p) \cos{q} + O(\epsilon^2) )$, and hence of minimal period $2\pi$. On the other hand, if $\sigma$ lies in the set $\Sigma_2$, then we expect double bifurcation to occur  and there will be several curves of $2\pi$-periodic solutions that, nevertheless, have a minimal period strictly less than $2\pi$.  The last set, $\Sigma_3$, consists of those values of $\sigma$ where there is a zero eigenvalue of the linearized problem.  We shall prove that $\Sigma_3$ is either a singleton or the empty set (see \textsc{Lemma \ref{zeroeigenvaluelemma}}).  Furthermore, in \textsc{Corollary \ref{Sigma12corollary}}, we show that every $\sigma$ for which \eqref{lb} holds is in precisely one of $\Sigma_1$,  $\Sigma_2$ and $\Sigma_3$.   We strongly suspect that $\Sigma_1$ is generic.  In fact, in section \ref{examplessection} we provide an elementary proof of the (previously known) result that, for irrotational, constant density waves, $\Sigma_2$ is countable.  

This section is concerned with the simple bifurcation case, so our first task is to produce sufficient conditions under which $\sigma \in \Sigma_1$.  In order to apply the linear theory developed in \cite{walsh2009stratified}, we introduce a modified density function, $\tilde{\rho}$, defined by 
\[\tilde{\rho} := \rho+\frac{\sigma}{g}.\]
It follows that, for each $\lambda \geq -2B_{\textrm{min}} + \epsilon_0$, $(H(\cdot; \lambda), \tilde{Q}(\lambda))$ is a laminar solution with streamline density function $\tilde{\rho}$, where 
\be \tilde{Q}(\lambda) := 2\sigma H(0;\lambda)+Q(\lambda) = \sigma \left( \frac{\left(\sigma^{-1}g\rho(0)+1\right)Q(\lambda)-\lambda}{g\rho(0)} \right) =  \frac{  g \tilde{\rho}(0)Q(\lambda)-\sigma \lambda }{g\rho(0)}. \label{defQtilde} \ee

It is important to note that $\tilde{Q}$ inherits all of the important properties of $Q$.  These we summarize in the following lemma.  

\begin{lemma}  Let $\lambda > -2B_{\textrm{min}} + \epsilon_0$ be given and let $\tilde{Q}$ be defined as in \eqref{defQtilde}.  Then $\tilde{Q}$ is strictly convex as a function of $\lambda$ with a unique minimum occurring at $\tilde{\lambda}_0$.  Moreover, if  the unique minimum of $Q$ occurs at $\lambda_0$, then $\lambda_0 < \tilde{\lambda}_0$.  
\label{propQtilde}
\end{lemma}
\begin{proof2}
We know by \textsc{Corollary 3.7} of \cite{walsh2009stratified} that $Q$ is a strictly convex function of $\lambda$.  Differentiating \eqref{defQtilde} twice in $\lambda$ confirms that the same is true of $\tilde{Q}$.  To see $\tilde{Q}$ attains its minimum, note that evaluating \eqref{defQtilde} at $\lambda = \lambda_0$, we have $d{\tilde{Q}}/d\lambda < 0$.  But, since $dQ/d\lambda \nearrow 1$ as $\lambda \to \infty$, while $g\tilde{\rho} > \sigma$, taking $\lambda$ sufficiently large we must have $d{\tilde{Q}}/d\lambda > 0$.  Hence, there exists some $\tilde{\lambda}_0 > \lambda_0$ where $d{\tilde{Q}}/d\lambda = 0$. \end{proof2}

Now, since we are searching for a solution to \eqref{sturmliouvilleode} with $n^2 = 1$, we can simply replace $\rho$ with $\tilde{\rho}$ to rewrite the problem as something reminiscent of the stratified gravity wave case.  In particular, we see that it is equivalent to finding a solution $M$ of
\be\left\{ \begin{array}{ll}
H_p \{H_p^{-3} M_p \}_p -g \tilde{\rho}_p H_p M =  M,  & p_0 < p < 0\\ 
M  = 0, & p = p_0 \\
\lambda^{3/2} M_p = g\tilde{\rho} M, & p = 0. \end{array}\right. \label{sturmliouvilleode2}  \ee

We have therefore reduced the problem to that considered in \cite{walsh2009stratified}.  Recalling the methodology of that paper, let us now consider the minimization problem:
\be \mu = \mu(\lambda) = \inf_{\phi \in \mathcal{S}} \mathcal{R}(\phi; \lambda) \label{minproblem} \ee
where $\mathcal{S} := \{ \phi \in H^1((p_0,0)) \colon \phi \nequiv 0, ~\phi(p_0) = 0\}$ and the Rayleigh quotient $\mathcal{R}$ is given by
\be \mathcal{R}(\phi; \lambda) := \frac{-g \tilde{\rho}(0) \phi(0)^2 + \int_{p_0}^0 H_p(p;\lambda)^{-3} \phi_p(p)^2 dp }{ \int_{p_0}^0 \phi(p)^2 \left(H_p(p;\lambda)^{-1} + g \tilde{\rho}_p(p) \right)dp}. \label{defcalR} \ee
By the standard variational theory for Sturm-Liouville problems, the existence of a nontrivial solution to \eqref{sturmliouvilleode2} for some $\lambda_* > -2B_{\textrm{min}}+\epsilon_0$ is equivalent to the statement $\mu(\lambda_*) = -1$.   Given that fact, we make the following definition.

\begin{definition} We say the pseudo-volumetric mass flux $p_0$, coefficient of surface tension $\sigma$, streamline density function $\rho$ and Bernoulli function $\beta$ collectively satisfy the \emph{capillary-gravity local bifurcation condition} provided that
\be \textrm{there exists $\lambda > -2B_{\mathrm{min}}+\epsilon_0$ and $M \nequiv 0$ solving the linearized problem \eqref{sturmliouvilleode2}}. \tag{L-B} \label{lb} \ee
Equivalently,
\[ \inf_{\lambda \geq -2B_{\mathrm{min}}+\epsilon_0} \mu(\lambda) \leq -1.\]
This last statement can also be written 
\be  \inf_{\lambda \geq -2B_{\mathrm{min}} + \epsilon_0} ~\inf_{\phi \in \mathcal{S}} \left\{ \frac{-g\rho(0)\phi(0)^2-\sigma\phi(0)^2 + \int_{p_0}^0 H_p(p;\lambda)^{-3} \phi_p(p)^2 dp }{ \int_{p_0}^0 \phi(p)^2 \left(H_p(p;\lambda)^{-1} + g \rho_p(p)\right)dp} \right\} < -1. \label{lbexplicit} \ee
where $\mathcal{S} := \{ \phi \in H^1((p_0,0)) : \phi \nequiv 0,~\phi(p_0) = 0 \}$ and $H(\cdot;\lambda)$ is the laminar flow solution given by \textsc{Lemma \ref{laminarflowlemma}}.  \label{capgravlbc} \end{definition}

\eqref{lb} is both necessary and sufficient, but clearly suffers from being highly non-explicit.  However, as was alluded to in the remarks following \textsc{Theorem \ref{mainresult1}}, if we are willing to sacrifice necessity,  there is an entirely explicit sufficient condition under which the capillary-gravity local bifurcation holds, namely \eqref{sizecond}.  

\begin{lemma} \emph{(Sufficiency of Size Condition)} If $p_0$, $\rho$, $\beta$ and $\sigma$ satisfy the size condition \eqref{sizecond}, then they satisfy \eqref{lb}. \label{suffsizecondlemma}
\end{lemma}
\begin{proof2} This lemma follows from \textsc{Lemma 3.5} of \cite{walsh2009stratified} with $\tilde{\rho}$ in place of $\rho$.  \end{proof2}

We also remark that from \eqref{lbexplicit} it is evident that, if \eqref{lb} holds for $(p_0, \sigma_1, \rho, \beta)$, then it hold for $(p_0, \sigma_2, \rho, \beta)$ whenever $\sigma_2 \geq \sigma_1$. In particular, this implies that if the Local Bifurcation Condition is satisfied for the problem without surface tension (i.e., with $\sigma = 0$), then it is satisfied for all values of $\sigma > 0$ as well.  

For future reference we now quote another result from \cite{walsh2009stratified}.  

\begin{lemma} \emph{(Monotonicity)} $\mu$ is a strictly increasing function of $\lambda$ when $\mu(\lambda) < 0$.  \end{lemma}
\noindent
This result implies, for instance, that any $\lambda > -2B_{\textrm{min}}+\epsilon_0$ satisfying $\mu(\lambda) = -1$ is necessarily unique.  The next lemma gives the existence of such a $\lambda$.  

\begin{lemma} \emph{(Eigenvalue Problem)} Assuming \eqref{lb}, there exists $\lambda_* > -2B_{\mathrm{min}}+\epsilon_0$ and a nontrivial solution $M = M(p)$ to \eqref{sturmliouvilleode} with $n = 1$.  
\end{lemma}
\begin{proof2}  In the previous paragraphs we have shown that the producing such an $M$ is equivalent to finding $\lambda_*$ such that $\mu(\lambda_*) = -1$.  Condition \eqref{lb} ensures that $\mu(\lambda) < -1$, for some $\lambda > -2B_{\textrm{min}}+\epsilon_0$, so it only remains to show that $\mu(\lambda) \geq -1$ for $\lambda$ sufficiently large.  With that in mind, let $\lambda$ be given such that
\[ \lambda > -2B_{\mathrm{min}} + \left( g\|\rho_p\|_{L^\infty} + \sqrt{g\rho(0) + \sigma}\right)^2.\]
Recalling that $G(\cdot; \lambda) \geq 2B_{\mathrm{min}}$, we find that, for any such $\lambda$, 
\[ H_p^{-1} \geq H_p^{-1} + g\rho_p = (\lambda+G)^{1/2} + g\rho_p \geq (\lambda + 2B_{\textrm{min}})^{1/2} - g\|\rho_p\|_{L^\infty} \geq \sqrt{g\rho(0)+\sigma}.\]
Now let $\phi \in \mathcal{S}$ be given and take $\lambda$ as above.  Then
\begin{align*}
\int_{p_0}^0 \left(H_p^{-3} \phi_p^2 +(H_p^{-1} + g\rho_p) \phi^2 \right) dp & \geq \sqrt{g\rho(0)+\sigma} \int_{p_0}^0 \left((g\rho(0)+\sigma)\phi_p^2+\phi^2 \right)dp \\
& \geq 2\left( g\rho(0)+\sigma \right) \int_{p_0}^0 \phi \phi_p dp =  \left(g\rho(0)+\sigma\right)\phi(0)^2.
\end{align*}
From the last inequality we conclude $\mathcal{R}(\phi) \geq -1$.  Since the selection of $\phi$ above was arbitrary, this proves $\mu(\lambda) \geq -1$, which completes the lemma.  \end{proof2}

From Remark \ref{futurereferenceremark}, we have that 
\[ \dot{Q} = 1-2g\rho(0)\dot{Y}(p_0).\]
As we stated in the proof of \textsc{Lemma \ref{propQtilde}}, the right-hand side above vanishes for $\lambda = \lambda_0$, and approaches 1 monotonically from below as $\lambda \to \infty$.  Therefore, for some $\lambda_c$, we have 
\be 4 \dot{Y}(p_0;\lambda) < \frac{1}{g\rho(0)+g\|\rho_p\|_{L^\infty}|p_0|}, \qquad \textrm{for all } \lambda > \lambda_c. \label{deflambdac} \ee
  As $2g\rho(0)\dot{Y}(p_0;\lambda_0) = 1$, and $Q$ obtains its unique minimum at $\lambda_0$, we must have, in particular, that $\lambda_c$ occurs to the right of ${\lambda}_0$.  The significance of this quantity will be made clear later.  Briefly, the point is that, when $\lambda_* > \lambda_c$, we will show that the null space of the linearized problem is one-dimensional, whereas in general it may be two-dimensional.  In fact, when $\rho$ is taken to be a constant, one can instead take $\lambda_c = \lambda_0$ (cf. \cite{wahlen2006capgrav}), in which case it satisfies the implicit relation 
 \be \int_{p_0}^0 H_p(p;\lambda_c)^{3} dp = \frac{1}{g\rho(0)}. \label{deflambdac2} \ee
 Since $H_p(\cdot; \lambda)$ is strictly monotonic as a function of $\lambda$, \eqref{deflambdac2} serves as an equivalent definition of $\lambda_c$ when $\rho$ is constant. 

Next, put
\be \sigma_c := (g\rho(0))^2 \int_{p_0}^0 \left(H_p(p;\lambda_c)^{-1} + g\rho_p(p)\right)\left( \int_{p_0}^p H_p(s;\lambda_c)^3 ds \right)^2 dp.\label{defsigmac} \ee
The importance of this quantity is found in the following lemma.
\begin{lemma} Suppose that $p_0$, $\sigma$, $\rho$, and $\beta$ are given satisfying \eqref{lb} (or size condition \eqref{sizecond}).  Then, for a unique $\lambda_*$, there exists a solution of the form $m = M(p)\cos{q}$ to \eqref{linearheighteq}.  Moreover, if we assume that 
\be \sigma \geq \sigma_c, \label{sigmasizecond} \ee
then $\lambda_c < \lambda_*$.  (Note that this is an explicit statement as neither $H$ nor $\rho$ depend on $\sigma$.) \label{eigenvaluelemma} \end{lemma} 

\begin{proof2} 
Writing \eqref{lb} in terms of $\tilde{\rho}$ we find that $p_0$, $\tilde{\rho}$ and $\beta$ satisfy the Local Bifurcation Condition of \cite{walsh2009stratified}.  Thus applying \textsc{Lemma 3.4} through \textsc{Lemma 3.9}  of that paper gives the existence of a solution $M$ of \eqref{sturmliouvilleode2} for a unique $\lambda_*$.  As we have commented before, this is equivalent to solving \eqref{sturmliouvilleode}, so $m(q,p) := M(p) \cos{q}$ is indeed the solution to \eqref{linearheighteq} we seek.  

Next, we the  minimization problem \eqref{minproblem}.  For each $\lambda$, the function $M$ that attains the minimum, $\mu(\lambda)$, will satisfy the boundary conditions of the linearized problem \eqref{sturmliouvilleode2} as well as the ODE:
\[ \{H_p^{-3} M_p\}_p = -\mu(\lambda)(H_p^{-1} +g \rho_p)M. \]
By construction, therefore, $\mu(\lambda_*) = -1$. Fix $\lambda = \lambda_c$, and let $\phi \in \mathcal{S}$ be defined as 
\[ \phi(p) = \int_{p_0}^p a(s)^{-3} ds, \qquad p_0 < p < 0.\] 
Note that $\phi$ is therefore the solution of the following boundary value problem: 
\[ (a^3 \phi_p)_p = 0,~ \textrm{ for } p_0 < p < 0, \qquad \phi(p_0) = 0, \qquad \lambda^{3/2} \phi_p(0) = g\rho(0) \phi(0).\]
Evaluating $\mathcal{R}(\phi; \lambda)$ yields
\begin{align*}
{\mathcal{R}}(\phi; \lambda_c) & =  \frac{-g\tilde{\rho}(0)\phi(0)^2 + \int_{p_0}^0 a(p)^3 \phi_p(p)^2 dp }{ \int_{p_0}^0 \phi(p)^2 \left(a(p) + g \tilde{\rho}_p(p) \right)dp} \\ 
& =  \frac{-\sigma \phi(0)^2}{ \int_{p_0}^0 \phi(p)^2 \left(a(p;\lambda) + g \rho_p(p) \right)dp} \\
& =-\sigma \sigma_c^{-1}. \end{align*}
For $\sigma \geq \sigma_c$, the above identity implies $\mu(\lambda_c) \leq -1$. By the monotonicity of $\mu$, we conclude that $\lambda_c$ lies to the left of $\lambda_*$ as desired. \end{proof2}

We end this subsection with a lemma characterizing the set of $\sigma$ for which zero is an eigenvalue of the linearized problem.

\begin{lemma} \emph{(Zero Eigenvalue)} Fix $p_0$, $\rho$ and $\beta$.  There exists at most one value of  $\sigma$ for which both \eqref{lb} holds and there exists a nontrivial $(M, n)$ solving \eqref{sturmliouvilleode} with $n = 0$ and $\lambda = \lambda_*$.   In particular, $\Sigma_3$ is either a singleton or the empty set.\label{zeroeigenvaluelemma}
\end{lemma}

\begin{proof2}  Let $\sigma_0$ be given as in the statement of the lemma.  Then by \eqref{sturmliouvilleode} with $n=0$, there exists $M$ solving 
\[ (a^3 M_p)_p - g\rho_p (M-M(0)) = 0, \qquad M(p_0) = 0, \qquad \lambda^{3/2} M_p(0) = g\rho(0)M(0),\]
for $\lambda = \lambda_*$.    Note that all the effects of surface tension have been removed from the equation, as $\sigma$ occurs only as a coefficient of $n$ in \eqref{sturmliouvilleode}.  Indeed, this problem is identical to that studied in \cite{walsh2009stratified}, where it was shown that a nontrivial solution exists if and only if $\lambda_* = \lambda_0$.  

In order to prove the lemma, therefore, we must show that there is at most one value of $\sigma$ for which $\lambda_*$ can coincide with $\lambda_0$.  To make the dependence of $\lambda_*$ on $\sigma$ explicit, define
\[ \mathcal{R}(\phi; \lambda, \sigma) := \frac{-g\rho(0)\phi(0)^2-\sigma\phi(0)^2 +\int_{p_0}^0 a^3 \phi_p^2 dp}{\int_{p_0}^0 \phi^2 \left( a+g\rho_p\right)} dp,\]  
for $\phi \in \mathcal{S}$. Then, for each $\sigma$ for which \eqref{lb} is satisfied, we may let $\lambda_* = \lambda_*(\sigma)$ to be the smallest value of $\lambda > -2B_{\textrm{min}}+\epsilon_0$ such that
\[ \inf_{\phi \in \mathcal{S}} \mathcal{R}(\phi; \lambda_*(\sigma), \sigma) = -1.\]
Clearly, $\lambda_*$ is continuous in $\sigma$ on its domain.  

Now let $\sigma_1,\sigma_2$ be given with $\sigma_1 < \sigma_2$, \eqref{lb} hold for $\sigma_1$ (hence also for $\sigma_2$.)  By the definition of $\lambda_*(\cdot)$, there exists $\phi_1, \phi_2 \in \mathcal{S}$ such that 
\[  \mathcal{R}(\phi_1; \lambda_*(\sigma_1), \sigma_1) = -1 = \mathcal{R}(\phi_2; \lambda_*(\sigma_2), \sigma_2).\]
But, recalling the meaning of $\mathcal{R}(\cdot;\cdot,\cdot)$, we have
\begin{align*} \mathcal{R}(\phi_2; \lambda_*(\sigma_2), \sigma_2) & \leq \mathcal{R}(\phi_1; \lambda_*(\sigma_2), \sigma_2)\\
&= \mathcal{R}(\phi_1; \lambda_*(\sigma_2), \sigma_1) + \frac{\left( \sigma_1 - \sigma_2\right) \phi_1(0)^2}{\int_{p_0}^0 \phi_1^2 \left(a+g \rho_p\right)dp} \\
& <  \mathcal{R}(\phi_1; \lambda_*(\sigma_2), \sigma_1). \end{align*}
The strict inequality derives from the fact that $\phi_1(0) \neq 0$, in light of the equation \eqref{sturmliouvilleode} satisfied by $\phi_1$.  We have therefore shown $\mathcal{R}(\phi_1; \lambda_*(\sigma_2), \sigma_1) > -1$.  By the strict monotonicity of $\mu$, we conclude that $\lambda_*(\sigma_2)$ lies to the right of $\lambda_*(\sigma_1)$.  It follows that $\lambda_*$ is a strictly increasing function of $\sigma$. Hence there is at most one value of $\sigma$ for which $\lambda_*(\sigma) = \lambda_0$.
\end{proof2}

\subsection{Proof of local bifurcation from simple eigenvalues} \label{simpleigenvaluebifurcationsection}

In short, the previous section assures us that, if \eqref{lb} and size condition \eqref{sigmasizecond}  each hold, then, at some $\lambda_* > -2B_{\textrm{min}} + \epsilon_0$, there exists a solution of the linearized problem \eqref{linearheighteq}.  In order to infer the existence of curve of non-laminar flows local to this point, we shall make use of the classical theorem of Crandall and Rabinowitz.  In this section and onward, we shall use $\mathcal{N}(\mathcal{L})$ to denote the null space of a linear operator $\mathcal{L}$ and $\mathcal{R}(\mathcal{L})$ to denote its range.

\begin{theorem} \emph{(Crandall-Rabinowitz, \cite{crandall1971bifurcation})} Let $X$ and $Y$ be Banach spaces, $I \subset \mathbb{R}$ an open interval with $\lambda_* \in I$.  Suppose that $\mathcal{F} : I \times X \to Y$ is a continuous map with the following properties:
\begin{itemize}
\item[\emph{(i)}] $\mathcal{F}(\lambda, 0) = 0$, for all $\lambda \in I$;
\item[\emph{(ii)}] $D_1 \mathcal{F}$, $D_2 \mathcal{F}$ and $D_1 D_2 \mathcal{F}$ exist and are continuous, where $D_i$ denotes the Fr\'echet derivative with respect to the $i$-th coordinate;
\item[\emph{(iii)}]  $D_2 \mathcal{F}(\lambda_*,0)$ is a Fredholm operator of index 0, in particular, the null space is one-dimensional and spanned by some element $w_*$;
\item[\emph{(iv)}] $D_1 D_2 \mathcal{F}(\lambda_*, 0)w_* \notin \mathcal{R}(D_2 \mathcal{F}(\lambda_*,0))$.
\end{itemize}
 Then there exists a continuous local bifurcation curve $\{(\lambda(s), w(s)) \in \mathbb{R} \times X : |s| < \epsilon \}$ with $\epsilon > 0$ sufficiently small such that $(\lambda(0), w(0)) = (\lambda_*, w_*)$, and
 \[ \{ (\lambda, w) \in \mathcal{U} : w \neq 0, \mathcal{F}(\lambda, w) = 0 \} = \{(\lambda(s), w(s)) \in \mathbb{R} \times X : |s| < \epsilon \}\]
 for some neighborhood $\mathcal{U}$ of $(\lambda_*,0)$ in $\mathbb{R} \times X$.  Moreover, we have 
 \[ w(s) = sw_* +o(s), \qquad \textrm{in } X,~|s| < \epsilon. \]
 If $D_2^2 \mathcal{F}$ exists and is continuous, then the curve is of class $C^1$.  
\label{crandallrabinowitz}
\end{theorem}

All that is left for us now is to verify the hypotheses of the preceding theorem. As in the previous sections let the transformed fluid domain be
\[ R := \{(q,p) : 0 < q < 2\pi, ~p_0 < p < 0 \}, \]
with boundaries
\[ \qquad T := \{ (q,p) \in R : p = 0 \}, \qquad B := \{ (q,p) \in R : p = p_0 \},\]
and define
\[ X := \{ h \in C_{\textrm{per}}^{3+\alpha}(\overline{R}) : h = 0 \textrm{ on } B\}, \qquad
 Y = Y_1 \times Y_2 := C_{\textrm{per}}^{1+\alpha}(\overline{R}) \times C_{\textrm{per}}^{1+\alpha}(T).\]
Let $h(q,p) := H(p) + w(q,p)$.  Then by the full height equation, $w$ must satisfy the following PDE
\be \begin{split} (1+w_{q}^2)(H_{pp}+w_{pp})+w_{qq}(H_p+w_p)^2-2w_q w_{pq} (H_p+w_p) &  \\ 
\qquad -g(H+w-d(H)-d(w))(H_p+w_p)^3 \rho_p + (H_p+w_p)^3\beta(-p) &= 0 \qquad \textrm{ in } R,  \end{split} \label{interHwheighteq} \ee
\be 1+w_q^2 + (H_p+w_p)^2(2\sigma \kappa[w]+2g\rho(H+w)-Q) = 0 \qquad  \textrm{ on } T, \label{Hwheighteqboundcond} \ee
together with periodicity in $q$ and vanishing on $B$.  

Conforming to the framework of \cite{crandall1971bifurcation}, we introduce a nonlinear operator
\[\mathcal{F} = (\mathcal{F}_1, \mathcal{F}_2): (-2B_{\textrm{min}} + \epsilon_0, \infty) \times X \to Y\] 
defined, for $w \in X$, $\lambda > -2B_{\textrm{min}} + \epsilon_0$, by
\be \begin{split}
\mathcal{F}_1(\lambda,w) & :=  (1+w_{q}^2)(H_{pp}+w_{pp})+w_{qq}(H_p+w_p)^2-2w_q w_{pq} (H_p+w_p) \\ 
&  \qquad -g(H+w-d(H)-d(w))(H_p+w_p)^3 \rho_p + (H_p+w_p)^3\beta(-p) 
\end{split} \label{defF1} \ee
\be \mathcal{F}_2(\lambda,w)  :=  1+w_q^2 + (H_p+w_p)^2(2\sigma \kappa[w]+2g\rho(H+w)-Q). \label{defF2} \ee
Note that by definition of the laminar solution $H(\cdot; \lambda)$,  we have $\mathcal{F}(\lambda,0) \equiv 0$.  For future reference, we evaluate the Fr\'echet derivatives $\mathcal{F}_{1w}$, $\mathcal{F}_{2w}$ at $w = 0$:
\begin{align}
\mathcal{F}_{1w} &= \partial_p^2 + H_p^2 \partial_q^2 +3H_p^2 \beta(-p) \partial_p - 3g(H-d(H))H_p^2\rho_p \partial_p - gH_p^3 \rho_p (1-d) \label{F1w} \\ 
\mathcal{F}_{2w} &= \bigg(2g\rho H_p^2 +2H_p (-2\sigma \partial_q^2 + 2g\rho H-Q)\partial_p  \bigg)\bigg|_T \nonumber \\
& =  \bigg(2g\rho \lambda^{-1}-2\lambda^{1/2}\partial_p -2\sigma \lambda^{-1} \partial_q^2 \bigg)\bigg|_T. \label{F2w} \end{align}

In order to show that the eigenvalue at $\lambda_*$ is simple we must characterize the null space and range of the linearized operator $\mathcal{F}_w (\lambda_*, 0)$. This is accomplished in the following three lemmas

\begin{theorem} \emph{(Null Space)} For $\lambda_* \neq \lambda_0$, the null space of $\mathcal{F}_w(\lambda_*,0)$ is at most two-dimensional.  In particular, it is one-dimensional provided $\lambda_* > \lambda_c$.  Zero is an eigenvalue if and only if $\lambda_* = \lambda_0$.  \label{nullspacetheorem} \end{theorem}
\begin{proof2} In the previous section we showed that $M(p) \cos q \in \mathcal{N}(\mathcal{F}_w(\lambda_*,0))$, hence the null space is at least one-dimensional.  Using a generalization of the arguments in \cite{wahlen2006capgrav}, we now demonstrate that, when $\lambda_* \neq \lambda_0$, then $\mathcal{N}(\mathcal{F}_w(\lambda_*,0))$ can be at most two-dimensional.  Furthermore, when $\lambda_* > \lambda_c$ it is at most one-dimensional.

Let $m \in \mathcal{N}(\mathcal{F}_w(\lambda_*,0))$ be given.  Then, by the evenness of $m$ in $q$, we may decompose it via a cosine series:
\[ m(q,p) = \sum_{n=0}^\infty m_n(p) \cos{(nq)}.\]
By the  definition of $\lambda_*$, we have that $m_1$ is in the span of $M$, the minimizer of the Rayleigh quotient $\mathcal{R}$.   Because $m$ is in the null space of the linearized operator, each $m_{n_i}$ must satisfy \eqref{linearheighteq2} for $n = n_i$. 

The case when $m_0 \nequiv 0$ was dealt with in \textsc{Lemma \ref{zeroeigenvaluelemma}}.  Recall that it was argued that this occurs precisely when $\lambda_*  = \lambda_0$.  

Let us now consider the case $\lambda_* \neq \lambda_0$. Suppose that  $n_1, n_2$ and $n_3$ are distinct positive integers such that $1 \in \{ n_1, n_2, n_3\}$ and $m_{n_1}, m_{n_2}, m_{n_3} \nequiv 0$.  We shall prove that (i) for $\lambda > \lambda_c$, each $m_{n_i}$ is in the span of $m_1$, and (ii) at most two of $\{m_{n_1}, m_{n_2}, m_{n_3}\}$ are linearly independent.  For simplicity let us denote $m_{n_i} := \phi_i$ for $i = 1,2,3$.

As $n_i \geq 1$, from \eqref{sturmliouvilleode} we conclude that the functions $\phi_i$ solve the following ODE:
\[ -(a^3 \phi_i^\prime)^\prime + g\rho^\prime \phi_i = -n_i^2 a \phi, \qquad p_0 < p < 0,\]
along with the boundary conditions
\[ \phi_i(p_0) = 0, \qquad \lambda^{3/2} \phi_i^\prime(0) = (n_i^2 \sigma + g\rho(0)) \phi(0), \qquad i = 1,2,3. \label{nullspaceidentity1} \]
Multiplying the first equation by $\phi_j$ and integrating, we obtain the following identity: for any $p_0 \leq p \leq 0$, 
\be \left(a^3\phi_i\phi_j^\prime\right)\bigg|^p = \int_{p_0}^p a^3 \phi_i^\prime \phi_j^\prime dr + \int_{p_0}^p (n_i^2 a + g\rho^\prime) \phi_i \phi_j dr > 0.\label{nullspaceidentity2}\ee
When $p = 0$ we find from the boundary conditions that
\be -n_i^2 \left( \int_{p_0}^0 a\phi_i \phi_j dp - \sigma \phi_i(0) \phi_j(0) \right) = -g\rho(0)\phi_i(0)\phi_j(0) + \int_{p_0}^0 a^3 \phi_i^\prime \phi_j^\prime dp + \int_{p_0}^0 g\rho^\prime \phi_i \phi_j dp. \label{nullspaceidentity1a}\ee
In particular, when $i = j$ we get the useful identity:
\be -n_i^2 \left( \int_{p_0}^0 a\phi_i^2 dp - \sigma \phi_i(0)^2 \right) = -g\rho(0)\phi_i(0)^2 + \int_{p_0}^0 a^3 (\phi_i^\prime)^2 dp + \int_{p_0}^0 g\rho^\prime \phi_i^2 dp  . \label{nullspaceidentity1}\ee
Note that the quantity in parenthesis on the left-hand side need not have a sign.  In some sense, this is what is responsible for the appearance of double bifurcation.

Importing some notation from \cite{wahlen2006capgrav}, let us define the space $\mathbb{H} := L^2([p_0, 0]) \times \mathbb{C}$ and endow it with additional structure by prescribing the indefinite inner product 
\[ [ \tilde{u}_1, \tilde{u}_2 ]_{\mathbb{H}} := \int_{p_0}^0 a u_1 \overline{u_2} dp  - \sigma b_1 \overline{b_2}, \qquad \textrm{for each $\tilde{u}_i = (u_i,b_i) \in \mathbb{H}$, $i = 1,2$.}\]
Note that $[\cdot,\cdot]_{\mathbb{H}}$ consists of an infinite-dimensional positive definite term and a one-dimensional negative definite term. As a trivial consequence, we have that any maximal negative definite subspace of $\mathbb{H}$ (with respect to $[\cdot,\cdot]_{\mathbb{H}}$) is one-dimensional.  

We also observe that, taking \eqref{nullspaceidentity1a} and exchanging the roles of $\phi_i$ and $\phi_j$, reveals
\[ -n_i^2 \left( \int_{p_0}^0 a \phi_i \phi_j dp - \sigma\phi_i(0)\phi_j(0) \right) = -n_j^2 \left( \int_{p_0}^0 a \phi_i \phi_j dp - \sigma\phi_i(0)\phi_j(0) \right).\]
Recalling the definition of $[\cdot,\cdot]_{\mathbb{H}}$, this last identity implies immediately that,
\[ [\tilde{\phi}_i, \tilde{\phi}_j ]_{\mathbb{H}} = 0, \qquad \textrm{for } i,j = 1,2,3, i \neq j,\]
where $\tilde{\phi}_i := (\phi_i, \phi_i(0))$.  In other words, the $\tilde{\phi}_i$ are mutually orthogonal with respect to this indefinite inner product.  

Let $D \subset \mathbb{H}$ be the dense set $D := \{ (u,b) \in \mathbb{H} : u \in H^2([p_0,0]), u(p_0) = 0, b = u(0) \}$.  In particular, we have by the boundary conditions imposed in \eqref{linearheighteq2} that $\tilde{\phi}_i$ is an element of $D$.    

Now let $\lambda > \lambda_c$ be given.  Then, for any $\tilde{u} \in D \setminus \{0\}$, we estimate
\begin{align} 
 |b|^2 & =   |u(0)|^2 =  \left| \int_{p_0}^0 u^\prime(p) dp \right|^2 \nonumber \\
 & \leq  \int_{p_0}^0  \frac{1+\dot{G}(p;\lambda)}{a(p;\lambda)^{3}} dp  \int_{p_0}^0 \frac{a(p; \lambda)^3}{1+\dot{G}(p;\lambda)} |u^\prime(p)|^2 dp \nonumber \\
 & < 4\dot{Y}(p_0) \int_{p_0}^0 a(p;\lambda)^3 u^\prime(p)^2 dp. \nonumber \end{align}
 The last inequality come from the estimate $1+\dot{G} \geq 1/2$ and the identity \eqref{Ydotidentity} evaluated at $p = p_0$.   Since $\lambda_* > \lambda_c$, it follows from above that 
 \be |b|^2 <  \frac{1}{g\rho(0)+g\|\rho_p\|_{L^\infty}|p_0|} \int_{p_0}^0 a(p;\lambda)^3 u^\prime(p)^2 dp. \label{uinequality} \ee
Using \eqref{uinequality} with $\tilde{u} = \tilde{\phi}_i$, we find 
\[ \int_{p_0}^0 a^3 (\phi_i^\prime)^2 dp - g\rho(0) \phi_i(0)^2 - g\|\rho^\prime\|_{L^\infty} |p_0| > 0.\]
From \eqref{nullspaceidentity2} with $i = j$, however, we conclude that $\phi_i(\cdot)^2$ is a strictly increasing function attaining its maximum value at p = 0.  It follows that
\[ 0 < - g\rho(0) \phi_i(0)^2+ \int_{p_0}^0 a^3 (\phi_i^\prime)^2 dp  + \int_{p_0}^0 g\rho^\prime \phi_i^2 dp  = -n_i^2 [\tilde{\phi}_i, \tilde{\phi}_i]_{\mathbb{H}}.\]
Here the last equality is simply a restatement of \eqref{nullspaceidentity1}.  But now we have shown that each nontrivial $\tilde{\phi}_i$ is negative definite with respect to $[\cdot,\cdot]_{\mathbb{H}}$.  Since they are all mutually orthogonal, this is impossible unless they lie within the same one-dimensional subspace of $\mathbb{H}$.  Thus, for $\lambda > \lambda_c$, $\mathcal{N}(\mathcal{F}_w(\lambda,0))$ is given by the linear span of $m_1$.

For $\lambda \leq \lambda_c$, the quantity appearing on the right-hand side of \eqref{nullspaceidentity1} may be negative and so we must approach the analysis with greater care.  Let us consider the indefinite inner product
\[ [\tilde{u}_1, \tilde{u}_2]_D := \int_{p_0}^0 a^3 u_1^\prime \overline{u_2^\prime} dp + \int_{p_0}^0 g\rho^\prime u_1 \overline{u_2} dp - g\rho(0) u_1(0) \overline{u_2(0)}, \qquad \textrm{for every } \tilde{u}_1, \tilde{u}_2 \in D.\] 
We claim that every subspace of $D$ that is maximal and negative definite with respect to $[\cdot,\cdot]_D$ must be one-dimensional.  Since there is an infinite dimensional negative semi-definite term in the definition of $[\cdot,\cdot]_D$, proving this claim is not automatic.

Let $\tilde{u} \in D$ be given.  Then we estimate
\begin{align*} 
\int_{p_0}^0 g\rho^\prime(p)u(p)^2 dp & = \int_{p_0}^0 g\rho^\prime(p) \left( \int_{p_0}^p u^\prime(r)dr\right)^2 dp \\
& \geq \int_{p_0}^0 g\rho^\prime(p) (p-p_0) \int_{p_0}^p u^\prime(r)^2 dr dp \\
& = \int_{p_0}^0 u^\prime(r)^2 \int_r^0 g\rho^\prime(p) (p-p_0) dp dr \\
& \geq - \frac{1}{2}g\|\rho^\prime\|_{L^\infty} p_0^2 \int_{p_0}^0 u^\prime(p)^2 dp. \end{align*}
Recall, however, that by the definition of $\epsilon_0$, 
\[ a^3 \geq \epsilon_0^{3/2} \geq 2g\|\rho^\prime\|_{L^\infty} p_0^2 e^{|p_0|} > 2 g \|\rho^\prime\|_{L^\infty}p_0^2.\]
With this inequality we are more-or-less done, for combining it with the last inequality we get
\[ [\tilde{u},\tilde{u}]_D \geq \frac{1}{2} \int_{p_0}^0 a^3 (u^\prime)^2 dp - g\rho(0) u(0)^2.\]
The right-hand side above clearly defines an indefinite inner product on $D$ for which each maximal negative semi-definite subspace is of dimension one.  The same must therefore be true of $(D, [\cdot,\cdot]_D)$.  This proves the claim.  

Now, using the notation we have developed, \eqref{nullspaceidentity1a} can be rewritten: 
\be -n_i^2 [\tilde{\phi}_i, \tilde{\phi}_j]_{\mathbb{H}} = [ \tilde{\phi}_i, \tilde{\phi}_j]_D, \qquad i,j  = 1,2,3. \label{nullspaceidentity3} \ee
Since the $\tilde{\phi}_i$ are orthogonal with respect to $[\cdot,\cdot]_\mathbb{H}$, \eqref{nullspaceidentity3} says they must also be orthogonal with respect to $[\cdot,\cdot]_D$.  Appealing to the claim, we conclude that at most one of them is negative semi-definite with respect to $[\cdot,\cdot]_D$. Without loss of generality, we may therefore assume that $[\tilde{\phi}_2, \tilde{\phi}_2]_D,~[\tilde{\phi}_3,\tilde{\phi}_3]_D > 0$.  But then \eqref{nullspaceidentity3} implies $\tilde{\phi}_2$ and $\tilde{\phi}_3$ are each negative definite with respect to $[\cdot,\cdot]_{\mathbb{H}}$, which is impossible unless they are in the same one-dimensional subspace.  For general $\lambda \neq \lambda_0$, therefore, we have proved $\mathcal{N}(\mathcal{F}_w(\lambda,0))$ is at most two-dimensional.   \end{proof2}

\begin{longremark} From the arguments above it is clear that, when $\lambda_* = \lambda_0$, the dimension of $\mathcal{N}(\mathcal{F}_w(\lambda_*,0))$ is at most three.  This is simply a consequence of the fact that there can be at most two nonzero eigenvalues --- one whose eigenfunction is negative semi-definite with respect to $[\cdot,\cdot]_\mathbb{H}$ and one whose eigenfunction is positive definite.  When $\rho$ is a constant, one can leverage this observation to prove that the null space will be at most two-dimensional.  Unfortunately, for general $\rho$, $\phi_0$, the eigenfunction for $n = 0$, need not be orthogonal to the other eigenfunctions (with respect to $[\cdot,\cdot]_{\mathbb{H}}$.)
\end{longremark}

The previous lemma motivates the choice of notation in \eqref{defSigma1}--\eqref{defSigma3}:  $\Sigma_1$ consists of those values of $\sigma$ for which simple bifurcation occurs; $\Sigma_2$ consists of those for which we have double bifurcation; and $\Sigma_3$ (if it is nonempty) contains the unique value of $\sigma$ for which zero is an eigenvalue and we cannot rule out triple bifurcation.  The next corollary makes this statement more precise.

\begin{corollary}  Fix $p_0, \rho$ and $\beta$.  Denote 
\[ \Sigma := \{ \sigma \in \mathbb{R}^+ : (p_0, \sigma,\rho,\beta) \textrm{ collectively satisfy \eqref{lb}}\}.\]  Then 
\[ (\sigma_c, \infty] \subset \Sigma_1 \qquad \textrm{and} \qquad \Sigma_1 \cup \Sigma_2 \cup \Sigma_3 = \Sigma.\] 
\label{Sigma12corollary} \end{corollary}
\begin{proof2} \textsc{Theorem \ref{nullspacetheorem}} implies that, if $\lambda_* \neq \lambda_0$, then $\mathcal{N}(\mathcal{F}_w(\lambda_*, 0))$ is either one-dimensional or two-dimensional.   First consider the one-dimensional case, i.e. assume there exists  $\phi \in X$ generating the null space.      Since $\phi$ is even in $q$, we may write it as a cosine series:
\[ \phi(q,p) = \sum_{n=0}^\infty M_n(p) \cos{(nq)}. \]
Observe that, because $\mathcal{F}_w(\lambda_*,0)\phi = 0$, we have additionally 
\[ \mathcal{F}_w(\lambda_*,0) \left(M_n(p) \cos{(nq)}\right)= 0, \qquad \textrm{for all } n \geq 0.\]
The local bifurcation condition guarantees that $M_1 \nequiv 0$.  But, since $\mathcal{N}(\mathcal{F}_w(\lambda_*,0))$ is one-dimensional, it follows that $M_n \equiv 0$, for $n \neq 1$.  Hence there is a unique solution $(n_1,M_1)$ to \eqref{sturmliouvilleode} with $\lambda = \lambda_*$, $n_1 = 1$.  We infer that, if $\sigma > \sigma_c$, then $\sigma \in \Sigma_1$.

In particular, note that by \textsc{Lemma \ref{eigenvaluelemma}}, when $\sigma > \sigma_c$, then the corresponding value of $\lambda_*$ will lie to the right of $\lambda_c$ and thus the null space is one-dimensional.  In other words, $(\sigma_c, \infty] \subset \Sigma_1$.

Suppose now that $\lambda_* \neq \lambda_0$ and the null space of $\mathcal{F}_w(\lambda_*,0)$ is two-dimensional with generators $\phi_1$ and $\phi_2$.  Then, decomposing each of these into cosine series we find that 
\[ \phi_1(q,p) = M_1(p) \cos{(n_1 q)}, \qquad \phi_2(q,p) = M_2(p) \cos{(n_2 q)},\]
where $(n_1, M_1)$ and $(n_2, M_2)$ are distinct solutions to \eqref{sturmliouvilleode} with $\lambda = \lambda_*$, $n_1 = 1$ and $n_2 \neq 1$.  In fact, by \textsc{Lemma \ref{zeroeigenvaluelemma}}, we actually have $n_2 \geq 2$.  By definition this means $\sigma \in \Sigma_2$. 

Finally, if $\lambda_* = \lambda_0$, then in view of \textsc{Lemma \ref{zeroeigenvaluelemma}}, we have that  $\sigma \in \Sigma_3$  --- indeed, it must be the case that $\Sigma_3 = \{\sigma\}$.   Since these possibilities are exhaustive, we conclude $\Sigma = \Sigma_1 \cup \Sigma_2 \cup \Sigma_3$.  \end{proof2}

\noindent
\begin{lemma} \emph{(Local Fredholm Map)} For each $\lambda$, the map $\mathcal{F}_w(\lambda,0): X \to Y$ is Fredholm with index 0. \label{localfredholmlemma} \end{lemma} 

\begin{proof2}  This follows trivially from combining \textsc{Lemma 3.5} of \cite{wahlen2006capgrav} and \textsc{Lemma 4.2} of \cite{walsh2009stratified}. \end{proof2}

\begin{lemma} \emph{(Range)} The pair $(\mathcal{A},\mathcal{B})$ belongs to the range of the linearized operator $\mathcal{F}_w(\lambda_*,0)$ if and only if it satisfies the orthogonality condition:
\[ \int \!\!\! \int_{R} \mathcal{A} a^3 \phi^* dqdp + \frac{1}{2} \int_T \mathcal{B} a^2 \phi^* dq = 0, \]
where $\phi^*$ generates the null space of $\mathcal{F}_w(\lambda_*,0)$. \label{rangelemma} \end{lemma}
\begin{proof2} Suppose first that $(\mathcal{A},\mathcal{B}) \in Y$ is in the range of $\mathcal{F}_w(\lambda_*,0)$.  Then there exists some $v \in X$ such that $\mathcal{A} = \mathcal{F}_{1w}(\lambda_*,0)v$, and $\mathcal{B} = \mathcal{F}_{2w}(\lambda_*,0)v$.  It follows that
\begin{align*}
\iint_R \mathcal{A} a^3 \phi^* dqdq & =  \iint_R \big((a^3 v_p)_p + av_{qq}- a^3 g \rho_p v\big) \phi^* dqdp \\
& =  \iint_R \big( (a^3 \phi_p^*)_p + a\phi_{qq}^* - a^3 g \rho_p \phi^*\big) v dqdp + \int_T (a^3 v_p \phi^*-a^3v\phi_p^*) dq \\
& =   \int_T (a^3 v_p \phi^*-a^3v\phi_p^*) dq, \end{align*}
as on $B$, both $v$ and $\phi^*$ vanish.  Now, on $T$ we have $a^3 = \lambda^{3/2}$.  Moreover, from the definition of $\phi^*$ we have on $T$ that $g\rho \phi^* - \sigma \phi^*_{qq} = \lambda^{3/2}\phi_p^*$.  Thus,
\begin{align*}
\frac{1}{2} \int_T \mathcal{B} a^2 \phi^* dq & =  \int_T \big( g\rho \lambda^{-1} v - \lambda^{1/2} v_p -  \sigma \lambda^{-1} v_{qq}\big) \lambda \phi^* dq \\
& =  \int_T \big((g\rho\phi^* -\sigma \phi_{qq}^*)v - \lambda^{3/2} \phi^* v_p \big)dq \\ 
& =  \int_T \big(\lambda^{3/2} \phi_p ^* v - \lambda^{3/2} \phi^* v_p \big)dq. \end{align*}
Combining this expression with the last gives necessity of the orthogonality condition.

Sufficiency follows in the same way as \textsc{Lemma 3.6} of \cite{wahlen2006capgrav}.  The main point of that argument is that, once we have established $\mathcal{F}_w(\lambda,0)$ is Fredholm of index 0, then dimension counting and necessity imply sufficiency. For brevity we omit the details.  \end{proof2}
   
Finally, we must ensure that the so-called transversality or crossing condition holds.

\begin{lemma} \emph{(Transversality Condition)} Suppose that $\phi^*$ generates the null space of $\mathcal{F}(\lambda_*,0)$. Then $\mathcal{F}_{w\lambda}(\lambda_*,0)\phi^* \notin \mathcal{R}(\mathcal{F}_w(\lambda_*,0))$. \label{technicallemma} \end{lemma}

\begin{proof2} First we calculate the mixed Fr\'ecfhet derivatives of $\mathcal{F}$ at $(\lambda_*,0)$:
\begin{align*}
\mathcal{F}_{1\lambda w}(\lambda_*,0) & =  -(1+\dot{G})a^{-4} \partial_q^2 - 3(1+\dot{G})a^{-4} \beta(-p) \partial_p - 3g \dot{Y}a^{-2}\rho_p \partial_p \\
& \qquad +3gY (1+\dot{G})a^{-4} \rho_p \partial_p + \frac{3}{2} g(1+\dot{G})a^{-5} \rho_p(1-d)\\
\mathcal{F}_{2\lambda w}(\lambda_*,0) & =  \bigg(-2g\rho \lambda^{-2} - \lambda^{-1/2} \partial_p +2\sigma \lambda^{-2} \partial_q^2 \bigg)\bigg|_T .\end{align*}
By the previous lemma, it suffices to show that the pair $(\mathcal{F}_{1\lambda w}(\lambda_*,0)\phi^*, \mathcal{F}_{2\lambda w}(\lambda_*,0) \phi^*)$ does not satisfy the orthogonality condition.  Equivalently, if we put 
\begin{align*}
\Xi & :=  \iint_{R}\bigg( (1+\dot{G})a^{-1} (\phi^*)^2 -3(1+\dot{G})a^{-1}\beta(-p)\phi^*\phi_p^* \\
& \qquad -3g\dot{Y} a \rho_p \phi^* \phi_p^* + 3gY(1+\dot{G})a^{-1} \rho_p \phi^* \phi_p^* + \frac{3}{2} g(1+\dot{G})a^{-2} \rho_p (\phi^*)^2 \bigg)dqdp \\
& \qquad\qquad +\int_T \bigg(-g\rho a^{-2} (\phi^*)^2 - \frac{1}{2}a \phi^* \phi_p^* +\sigma a^{-2} \phi_{qq}^* \phi^* \bigg) dq, \end{align*}
our lemma will be proven if we can show $\Xi \neq 0$ (again, because $d(\phi^*) = 0$).  We will demonstrate that, in fact, $\Xi < 0$.  To keep our notation concise, let $\Xi = \Xi_1+ \ldots + \Xi_8$, where $\Xi_i$ denotes the $i$-th term in the sum above.  Note that the only additional term from the $\Xi$ considered in \textsc{Lemma 3.12} of \cite{walsh2009stratified} is $\Xi_8$.  By simply integrating by parts once, and noting the boundary terms vanish by periodicity, we find
\[ \Xi_8 = -\sigma \lambda_*^{-1} \int_{T} \phi_q^2 dq \leq 0.\]
The lemma follows.  \end{proof2}

Combining these lemmas, we now prove the main theorem for simple eigenvalues.

\begin{proofof}{Theorem \ref{capgravsimplelocalbifurcation}.}  That $\mathcal{F}$ satisfies conditions (i)-(iii) of \textsc{Theorem \ref{crandallrabinowitz}} is elementary.   Likewise, $\mathcal{F}_{ww}(\lambda_*,0)$ exists and is continuous.  Finally, conditions (iii)-(iv) follow from \textsc{Lemma \ref{localfredholmlemma}} and \textsc{Lemma \ref{technicallemma}}, respectively.  

We are therefore able to conclude the existence of a $C^1$-curve $\mathcal{C}_{\textrm{loc}}^{\prime}$ of non-laminar solutions to the height equation bifurcating from $\mathcal{T}$ at $(H(\cdot; \lambda_*), Q(\lambda_*))$.  In a neighborhood $\mathcal{U}$ of the bifurcation point, moreover, this curve and $\mathcal{T}$ together comprise the complete solution set and each $h \in \mathcal{C}_{\textrm{loc}}^{\prime} \cap \mathcal{U}$ can be written
\be h(q,p) = H(p; \lambda_*)+\epsilon M(p) \cos{q} + o(\epsilon), \qquad \textrm{in } X. \label{lineardecomposition} \ee
Since $H_p(\cdot; \lambda) > 0$, \eqref{lineardecomposition} implies  that, possibly by shrinking $\mathcal{U}$, we may assume that $h_p > 0$ on $\mathcal{C}_\textrm{loc}^\prime$.  Then, in light of \textsc{Lemma \ref{equivalencelemma}}, there is a corresponding $C^1$-curve $\mathcal{C}_{\textrm{loc}}$ of solutions to the original problem, \eqref{incompress}--\eqref{boundcond}.  Since $\mathcal{C}_\textrm{loc}^\prime \cap \mathcal{T} = \{(\lambda_*,0)\}$, we are guaranteed that $\mathcal{C}_\textrm{loc}$ contains only one laminar solution: that corresponding to $H(\cdot; \lambda_*)$.  Finally, part (iii) of the Theorem statement follows from \eqref{lineardecomposition}, possibly by further restricting to a subneighborhood of $(H_*, Q_*)$ in $\mathcal{U}$.  \end{proofof}

\section{Local bifurcation from double eigenvalues} \label{doubleeigenvaluebifurcationsection}

If we do not assume that $\sigma \in \Sigma_1$, then we cannot rule out the possibility that null space of $\mathcal{F}_{w}( \lambda_*,0)$ is two-dimensional (or, indeed, that zero is an eigenvalue.)  Were this to be the case, the Crandall-Rabinowitz arguments of the previous section break down and we are forced to revert to a more general methodology.  To successfully prosecute this program requires that we impose some additional nondegeneracy conditions. Suppose that $\sigma \in \Sigma_2 \cup \Sigma_3$.  We may therefore let $n_2 \neq 1$ be given so that there exists a nontrivial solution to \eqref{sturmliouvilleode} with $\lambda = \lambda_*$, $n = n_2$.  We require
\be n_2 \geq 3 \label{capgravnondegeneracycondition1} \ee
and 
\be 0 \neq \Theta_{1111}\Theta_{2222} \neq \Theta_{1122}\Theta_{2211}, \qquad 
\Theta_{1111}\Psi_{22} \neq \Theta_{2211}\Psi_{11} \qquad \Theta_{2222} \Psi_{11} \neq \Theta_{1122}\Psi_{22}, \label{capgravnondegeneracycondition2} \ee
where $\Theta_{iijj}$, $\Psi_{ii}$ are given by \eqref{Thetaiiii}--\eqref{Thetaiijj} and \eqref{Psiii}, respectively, for $i,j = 1,2$.   In particular, we note that \eqref{capgravnondegeneracycondition1} implies $\sigma \notin \Sigma_3$.  The purpose of dictating $n_2 \neq 0$ is merely to avoid considering the zero eigenvalue case, which require some additional work.  In fact, Wahl\'en was able to give a fairly complete analysis of the local bifurcation when this does occur, though the argument is quite subtle (cf. \cite{wahlen2006capgrav}).  The reason for avoiding $n_2 = 2$ is less obvious, but results from certain symmetries in the underlying problem that shall become clear later (in particular, see \textsc{Lemma \ref{quadratictermslemma}}.)

Under these assumptions, the eventual product of our efforts in this section will be the proof of \textsc{Theorem \ref{mainresult2}}.  
We shall accomplish this in several stages, according to the Lyapunov-Schmidt procedure.  The techniques here are not new, though, as we shall see, there is a certain degree of degeneracy in the problem that will require some care. \\

\paragraph{Step I:  Reduction to finite-dimensional problem.}  We begin by recalling a weighted inner-product on $Y$ introduced in the previous section.  For all pairs $(\mathcal{A}_i, \mathcal{B}_i) \in Y$, $i = 1,2$, let
\be
\begin{split} \left( \left(\mathcal{A}_1, \mathcal{B}_1\right),\left(\mathcal{A}_2,\mathcal{B}_2\right)\right)_Y &:= \iint_R a^3(\cdot; \lambda_*)\mathcal{A}_1(\cdot)\mathcal{A}_2(\cdot) dqdp \\
& \qquad + \frac{1}{2} \int_T a^2(\cdot; \lambda_*) \mathcal{B}_1(\cdot) \mathcal{B}_2(\cdot) dq.  \end{split} \label{defYinnerprod} \ee
Denote $X_A := \mathcal{N}(\mathcal{F}_w(\lambda_*,0))$.  Then, as we have seen in \textsc{Theorem \ref{nullspacetheorem}}, $X_A$ will be at most dimension two.  Suppose that it is spanned by $\phi_1$ and $\phi_2$, where $\phi_i(q,p) = M_i(p) \cos{n_i q}$, $i  = 1,2$.  Of course, the most interesting case given our previous analysis, is when $n_1 = 1$. But there is little to gain at this stage in making such a restriction.

It follows directly from \textsc{Lemma \ref{rangelemma}} that the range of $\mathcal{F}_w(\lambda_*, 0)$ is given by the orthogonal complement (with respect to $(\cdot, \cdot)_Y$) of the linear span of $\{\phi_1, \phi_2\}$.  So, if we now define 
\begin{align*}
X_B &:= \left\{ h \in X : \left((h,h|_T),(\phi_i,\phi_i|_T)\right)_Y = 0,~i=1,2 \right\}, \\
Y_A &:= \left\{ (\mathcal{A},\mathcal{B}) \in Y : \mathcal{A}|_T = \mathcal{B},~(\mathcal{A},\mathcal{B}) \in \textrm{span}\{(\phi_1, \phi_1|_T), (\phi_2, \phi_2|_T)\}\right\}, \\
Y_B &:= \left\{ (\mathcal{A},\mathcal{B}) \in Y : \mathcal{A}|_T = \mathcal{B},~\left((\mathcal{A},\mathcal{B}),(\phi_i,\phi_i|_T)\right)_Y = 0,~i=1,2 \right\}, \end{align*}
then $X = X_A \oplus X_B$, $Y = Y_A \oplus Y_B$.  Finally, let $A_X: X\to X_A$, $B_X : X \to X_B$ denote the projections onto $X_A$ and $X_B$, respectively, and similarly for $A_Y:Y \to Y_A$, $B_Y:Y \to Y_B$.

As always in bifurcation theory, we are trying to find all solutions to the problem
\be \mathcal{F}(\lambda, w) = 0, \label{lyapunovproblem} \ee
in a neighborhood of $(\lambda_*,0) \in \mathbb{R}\times X$.  Applying the projections above, we see that this is equivalent to solving
\begin{align}
A_Y \mathcal{F}(\lambda,A_X w + B_X w) &= 0, \label{finiteequation} \\
 B_Y \mathcal{F}(\lambda,A_Xw + B_X w) &=  0,\label{infiniteequation} \end{align}
locally near $(\lambda_*,0)$.  

Consider \eqref{infiniteequation}.   Merely by unravelling the definitions above, we have $A_Y \mathcal{F}_w(\lambda_*,0) B_X \equiv 0$ and therefore, for any $w \in X$, 
\begin{align*}
B_Y \mathcal{F}_w(\lambda_*,0) w - \mathcal{F}_w(\lambda_*,0) B_X w  &=  B_Y\mathcal{F}_w(\lambda_*,0) (1-B_X) w =  B_Y\mathcal{F}_w(\lambda_*,0) A_X w = 0.\end{align*}  
This proves that $\mathcal{F}_w(\lambda_*, 0)$ commutes with projections onto the infinite dimension spaces in the sense that 
\be \mathcal{F}_w(\lambda_*, 0) B_X = B_Y \mathcal{F}_w(\lambda_*,0). \label{commutes} \ee

Now we note that, by construction, $\mathcal{F}_w(\lambda_*,0) B_X$ is an isomorphism from $X$ onto $Y_B = \mathcal{R}(\mathcal{F}_w(\lambda_*,0))$.  Thus, in view of \eqref{commutes}, the Implicit Function Theorem allows us to solve for $B_X w$ as a function of $A_X w$ and $\lambda$ near $(\lambda_*,0)$.   In particular, if we write $A_X w = \xi_1 \phi_1 + \xi_2 \phi_2$, for some $\xi_1, \xi_2 \in \mathbb{R}$, then there exists a neighborhood $U \subset \mathbb{R}^2 \times \mathbb{R}$ of $(0,0,\lambda_*)$ and a smooth $\zeta: U \to Y_B$ such that, in a sufficiently small neighorhood of $(\lambda_*,0)$ in $\mathbb{R} \times X$,  all solutions of \eqref{infiniteequation} are given by $\{ (\lambda,\xi_1 \phi_1 + \xi_2 \phi_2 + \zeta(\xi_1,\xi_2; \lambda)) : (\xi_1, \xi_2, \lambda) \in U\}$.  Moreover, it is easy to see that $\zeta(\cdot,\cdot; \lambda)  = O(|\xi|^2)$, in this neighborhood.

Taking these solutions and inserting them into \eqref{finiteequation}, we infer that \eqref{lyapunovproblem} reduces to the following finite-dimensional problem:  Find all $(\xi_1, \xi_2) \in \mathbb{R}^2,~\lambda \in \mathbb{R}$ in a sufficiently small neighborhood $U$ of $(0,0,\lambda_*)$ such that  
\be \mathfrak{B}(\xi_1, \xi_2; \lambda) = 0, \label{bifurcationequation} \ee
where the bifurcation function $\mathfrak{B}: U \to \mathbb{R}^2$ is defined by
\be \mathfrak{B}_i(\xi_1,\xi_2; \lambda) := \left( \phi_i, \mathcal{F}(\lambda, \xi_1 \phi_1 + \xi_2\phi_2+\zeta(\xi_1, \xi_2; \lambda)) \right)_Y, \qquad i = 1,2. \label{defbifurcationfunction} \ee

To simplify this expression, we repeatedly Taylor expand $\mathcal{F}$.  In doing so, we shall let  
\[\langle \mathcal{F}_{ww}(\lambda, 0) \phi_i, \phi_j \rangle \in \mathcal{L}(X,Y)\] 
denote the second Fr\'echet derivative of $\mathcal{F}$ in the $\phi_i$ then $\phi_j$ direction, $i, j = 1,2$.  With this notation, expanding \eqref{defbifurcationfunction} around $(\lambda_*,0)$ yields, for $i,j,k,\ell = 1,2$, 
\be \mathfrak{B}_i = (\lambda- \lambda_*) \xi_j \Psi_{ij} + \xi_j \xi_k \Phi_{ijk} + \xi_j \xi_k \xi_{\ell} \Theta_{ijk\ell} + O(|\lambda- \lambda_*|^2 |\xi| + |\lambda-\lambda_*||\xi|^4 + |\xi|^5 ),\label{expandedbifurcationfunction} \ee
where we are using summation convention and
\begin{align}
\Psi_{ij} &:= \frac{1}{2}\left( \phi_i,~\mathcal{F}_{\lambda w}(\lambda_*, 0) \phi_j \right)_Y, \qquad i,j = 1,2, \label{defPsiij} \\
\Phi_{ijk} &:= \frac{1}{2}\left( \phi_i, ~\left \langle \mathcal{F}_{ww}(\lambda_*, 0)\phi_j, \phi_k \right \rangle \right)_Y , \qquad i,j,k = 1,2.\label{defPhiijk} \\
\Theta_{ijk\ell} & :=  \frac{1}{6}\left( \phi_i,~\left\langle \left\langle \mathcal{F}_{www}(\lambda_*,0) \phi_j, \phi_k \right\rangle, \phi_{\ell} \right\rangle \right)_Y, \qquad i,j,k,\ell = 1,2. \label{defThetaijkl} \end{align}  \\

\paragraph{Step II: Computation of coefficients.}  Clearly, before we can begin to characterize the solutions of \eqref{bifurcationequation}, we must first determine $\Psi_{ij}$, $\Phi_{ijk}$ and $\Theta_{ijk\ell}$ more precisely.  As one might expect, this will be an elementary but involved computation.  What will simplify things immensely, however, is the special structure of $\phi_1$ and $\phi_2$.  That is, we know $\phi_i(q,p) = M_i(p) \cos{(n_i q)}$, where $M_i$ is the solution of \eqref{sturmliouvilleode} with $n = n_i$, for $i = 1,2$.  This will aid us in evaluating the integrals in $(\cdot, \cdot)_Y$.

Also, observe that multiplying \eqref{sturmliouvilleode} by $M_i$ and integrating by parts one obtains, for $p \in (p_0, 0]$, $i = 1,2$,
\[ M_i M_i^\prime = a^{-3} \int_{p_0}^p a^3 (M_i^\prime)^2 dr + a^{-3} \int_{p_0}^p (an_i^2 + g\rho_p) M_i^2 dr .\]
As the right-hand side above is strictly positive, we surmise that $M_i$ and $M_i^\prime$ do not vanish on $(p_0, 0)$ and are of the same sign. By linearity of the equation, we may simply take them both to be positive.  
\begin{longremark} In general, for a Sturm-Liouville problem one has that the only eigenfunction without nodes is the ground state.  We have argued, however, that $M_i, M_i^\prime$ are both nonvanishing for $i = 1,2$.  This seeming contradiction is precisely due to the presence of the surface tension term on the boundary, which alters the orthogonality relation between the eigenfunctions.  
\end{longremark}

\begin{lemma} \emph{(Linear Terms)}  The matrix $\Psi_{ij}$ defined by \eqref{defPsiij} is diagonal and negative definite.  \label{lineartermslemma} \end{lemma}
\begin{proof2} We shall start by explicitly computing $\Psi_{ij}$.  Using the expressions for $\mathcal{F}_{\lambda w}(\lambda_*,0)$ found in the previous section, we see that, for $i,j = 1,2$,
\be \begin{split} \Psi_{ij} & =  \bigg[\int_{p_0}^0 n_j^2 a^{-1}(1+\dot{G}) M_i M_j dp -3 \int_{p_0}^0 (1+\dot{G})a^{-1}\beta M_i M_j^\prime dp -3g \int_{p_0}^0 \dot{Y}a \rho_p M_i M_j^\prime dp  \\
&  \qquad +3g \int_{p_0}^0 Y(1+\dot{G}) a^{-1} \rho_p M_i M_j^\prime dp +\frac{3}{2}g\int_{p_0}^0(1+\dot{G})a^{-2} \rho_p M_i M_j \\
& \qquad\qquad -\frac{1}{2}\left(2g\rho \lambda^{-1}M_i M_j + \lambda^{1/2}M_i M_j^\prime + 2\sigma \lambda^{-1} n_j^2 \right)\bigg|_T\bigg]  \int_0^{2\pi} \cos{(n_i q)}\cos{(n_j q)} dq.   \end{split} \label{Psiii} \ee
Thus, $\Psi_{ij} = 0$ for $i \neq j$.  Arguing as in \textsc{Lemma \ref{technicallemma}}, moreover, we can show that $\Psi_{ii} < 0$, $i = 1,2$.  \end{proof2}

Let us now attempt to characterize $\Phi_{ijk}$.  Unfortunately, the conclusions we will be able to draw will be far less satisfying than for the previous result.

\begin{lemma} \emph{(Quadratic Terms)} Let $\Phi_{ijk}$ be defined as in \eqref{defPhiijk} then, if $n_2 \neq 2 n_1$, we have $\Phi_{ijk} = 0$, for $i,j,k = 1,2$.  On the other hand, if $n_2 = 2n_1$, then necessarily,  $\Phi_{111},\Phi_{222},\Phi_{221},\Phi_{212},\Phi_{122} = 0$.   In general, $\Phi_{112}$ and $\Phi_{122}$ ( $= \Phi_{211}$) are given by equations \eqref{Phi112} and \eqref{Phi211}, respectively.\label{quadratictermslemma} \end{lemma}

\begin{proof2} Recycling notation slightly, for the purpose of the next computation we let $\phi$ and $\psi$ be given elements of $X$ with $d(\phi) = d(\psi) = 0$.  This restriction is natural as we shall eventually be applying this to $\phi_1, \phi_2$, each of which have this property.  Then, 
\begin{eqnarray*}
\left\langle \mathcal{F}_{1ww}(\lambda, w) \phi, \psi \right\rangle & = &  2w_{q}\psi_{q}\phi_{pp}  + 2(H_{pp}+w_{pp}) \phi_{q} \psi_{q} +2w_{q}\psi_{pp} \phi_{q} \\
& & +2(H_p+w_p) \psi_p \phi_{qq} +2(H_p +w_p) \psi_{qq} \phi_p + 2w_{qq} \psi_p \phi_p \\
& & -2w_{q} \psi_{pq} \phi_p -2w_{qp} \psi_q \phi_p -2(H_p+w_p)\psi_q \phi_{pq} -2w_q \psi_p \phi_{pq} \\
& &  -2w_{pq} \psi_p \phi_q - 2(H_p+w_p)\psi_{pq} \phi_q  -3g\rho_p(H_p+w_p)^2 \psi_p \phi \\
& & -3g\rho_p (H_p+w_p)^2 \psi \phi_p +6 (H_p+w_p) \psi_p \phi_p \beta(-p) \\
& &  -6 g\rho_p (H_p+w_p)(H+w-d(H)-d(w)) \psi_p \phi_p ,
\end{eqnarray*}
and
\begin{eqnarray*}
\left\langle \mathcal{F}_{2ww}(\lambda, w)\phi,\psi \right\rangle & = & 2\phi_q \psi_q + 2\phi_p \psi_p \left(2\sigma \kappa[w] + 2g\rho(H+w)-Q \right) \\
& & +2(H_p+w_p)\psi_p \left(2\sigma \kappa^\prime[w] \phi + 2g\rho\phi\right) \\
& & + 2(H_p+w_p)\phi_p \left(2\sigma \kappa^\prime[w] \psi + 2g\rho\psi\right) \\
&  &+ 2\sigma(H_p+w_p)^2 \langle \kappa^{\prime\prime}[w] \phi, \psi \rangle, \end{eqnarray*}
where 
\[ \langle \kappa^{\prime\prime}[w] \phi, \psi \rangle = \frac{3(\phi_{qq}\psi_q+\phi_q \psi_{qq})}{(1+w_q^2)^{5/2}}   - \frac{15w_q w_{qq}}{(1+ w_q^2)^{7/2}} \phi_q \psi_q. \]
Evaluating each of these at $(\lambda_*, 0)$, we arrive at the following
\be \begin{split}
 \left\langle \mathcal{F}_{1ww}(\lambda_*, 0) \phi, \psi \right\rangle & =  2H_{pp} \phi_{q} \psi_{q}+2H_p \psi_p \phi_{qq} +2H_p \psi_{qq} \phi_p -2H_p \psi_q \phi_{pq} \\
 &  \qquad-2H_p \psi_{qp}\phi_q -3g\rho_p H_p^2 \psi_p \phi -3g\rho_p H_p^2 \psi \phi_p \\
 &  \qquad\qquad-6g\rho_p H_p(H-d(H)) \psi_p \phi_p + 6H_p \psi_p \phi_p \beta(-p) \end{split} \label{comptedF1ww} \ee
 \be \begin{split}
 \left\langle \mathcal{F}_{2ww}(\lambda_*,0) \phi, \psi \right\rangle & =  2\phi_q \psi_q + 2\phi_p \psi_p(2g\rho H - Q)  +2H_p \psi_p (2g \rho \phi - 2\sigma \phi_{qq}) \\ 
 &  \qquad+ 2H_p \phi_p( 2g\rho\psi -2\sigma \psi_{qq})+ 6\sigma H_p^2 (\phi_{qq} \psi_q + \phi_q \psi_{qq}). \end{split} \label{computedF2ww} \ee
Note that, as the above expressions make clear, the regularity of $\mathcal{F}_{ww}$ guarantees that $\langle \mathcal{F}_{ww}(\lambda,0) \phi, \psi \rangle = \langle\mathcal{F}_{ww}(\lambda,0) \psi, \phi \rangle$.  So, recalling the definition of $\Phi_{ijk}$ in \eqref{defPhiijk}, we see immediately that $\Phi_{ijk} = \Phi_{ikj}$, for $i,j,k = 1,2$.

To keep things manageable we will evaluate the inner-product in \eqref{defPhiijk} in two steps.  
 First we compute weighted the $L^2$ inner-product over $R$:
\begin{align*} 
\left( a^3 \phi_i, \left\langle \mathcal{F}_{1ww}(\lambda_*,0)\phi_j, \phi_k \right\rangle \right)_{L^2} & = C^{(1)}_{ijk} \left[n_j n_k \left(a^3 H_{pp} M_i, M_j M_k \right)_{L^2} - 2n_j n_k \left(a^2 M_i, (M_j M_k)^\prime\right)_{L^2} \right]\\
&  \qquad-C_{ijk}^{(2)} \bigg[2n_j^2 \left(a^2 M_i, M_j M_k^{\prime} \right)_{L^2} + 2n_k^2 \left(a^2 M_i, M_j^\prime M_k \right)_{L^2} \\
&  \qquad\qquad+3g \left(a \rho_p M_i, (M_jM_k)^\prime \right)_{L^2} -6 \left( a^2 \beta M_i, M_j^\prime M_k^\prime \right)_{L^2}   \\
&\qquad\qquad\qquad + 6g \left(a^2 \rho_p (H-d(H)) M_i, M_j^\prime M_k^\prime\right)_{L^2}\bigg], \end{align*}
where, for $i,j,k = 1,2$, we put
\begin{align*} C_{ijk}^{(1)} &:= \int_0^{2\pi} \cos{(n_i q)} \sin{(n_j q)} \sin{(n_k q)} dq, 
\\ C_{ijk}^{(2)} &:= \int_0^{2\pi} \cos{(n_i q)} \cos{(n_j q)} \cos{(n_k q)} dq.  \end{align*}
Similarly, evaluating the weighted $L^2$ inner-product on $T$ yields,
\begin{align*}
\frac{1}{2}\left(a^2 \phi_i, \left \langle \mathcal{F}_{2ww}(\lambda_*,0) \phi_j, \phi_k \right\rangle \right)_{L^2(T)} & =  C_{ijk}^{(3)} \left[n_j n_k \lambda_* M_i M_j M_k + 3n_j^2 n_k \sigma M_i M_j M_k \right] \\
&   \qquad + C_{ijk}^{(2)} \bigg[\lambda_*^{1/2}(2g\rho+2\sigma n_k^2) M_i M_j^\prime M_k \\
& \qquad\qquad + \lambda_*^{1/2} (2g\rho+2\sigma n_j^2) M_iM_j  M_k^\prime \\ 
&  \qquad\qquad\qquad-\lambda_*^2 M_i M_j^\prime M_k^\prime \bigg] \\
& \qquad +3C_{ijk}^{(4)} \sigma n_j n_k^2 M_i M_j M_k,
\end{align*}
where
\begin{align*} C_{ijk}^{(3)} &:= \int_0^{2\pi} \cos{(n_i q)} \cos{(n_j q)} \sin{(n_k q)} dq, \\
 C_{ijk}^{(4)} & :=  \int_0^{2\pi} \cos{(n_i q)} \sin{(n_j q)} \cos{(n_k q)} dq.\end{align*}

On both the boundary and interior, therefore, we will have massive cancellation due to the integrals of the various trigonometric terms.  For instance, when $i = j = k$, all of these integrals vanish, which implies $\Phi_{111} = \Phi_{222} = 0$.  More strikingly, unless $n_2 = 2 n_1$, \emph{every} term will vanish and $\Phi_{ijk} = 0$, for $i,j,k = 1, 2$. Suppose that we do have this relationship between $n_2$ and $n_1$, then, after the dust clears, we are left finally with the following:
\be \Phi_{111},~\Phi_{222},~\Phi_{221},~\Phi_{212},~\Phi_{122} = 0, \label{Phiijkcompute1} \ee
\be\begin{split} 
\frac{2}{\pi}\Phi_{211} & =  -n_1^2 \left( a^3 H_{pp} M_2, M_1^2 \right)_{L^2}  -3g \left( a \rho_p M_2,(M_1^2)^\prime \right)_{L^2}   \\
&  \qquad- 6g \left(a^2 \rho_p(H-d(H)) M_2, (M_1^\prime)^2 \right)_{L^2} + 6\left(a^2 \beta M_2, (M_1^\prime)^2 \right)_{L^2}\\ 
&  \qquad+ \bigg(4\lambda_*^{1/2}(2g\rho + 2\sigma n_1^2)M_2 M_1 M_1^\prime -\lambda_* n_1^2 M_2 M_1^2 -  \lambda_*^2 M_2 (M_1^\prime)^2 \bigg)\bigg|_T, \end{split}  \label{Phi211} \ee
and
\be \begin{split}
\frac{2}{\pi}\Phi_{112},\,\frac{2}{\pi} \Phi_{121} & =  n_1 n_2 \left(a^3 H_{pp} M_2, M_1^2 \right)_{L^2} - 2n_1^2 \left( a^2 M_1^2, M_2^\prime \right)_{L^2} -n_2^2 \left( a^2 M_2, (M_1^2)^\prime \right)_{L^2}  \\
&\qquad -3g \left( a \rho_p M_1,(M_1 M_2)^\prime \right)_{L^2} - 3g \left( a^2 \rho_p(H-d(H)) M_2^\prime, (M_1^2)^\prime\right)_{L^2}\\
& \qquad\qquad  + 3\left( a^2 \beta M_2^\prime, (M_1^2)^\prime \right)_{L^2} -2n_1 n_2 \left( a^2 M_1, (M_1 M_2)^\prime \right)_{L^2} \\ 
&  \qquad\qquad\qquad+ \bigg(\lambda_* n_1 n_2 M_1^2 M_2+\lambda_*^{1/2}(2g\rho + 2\sigma n_2^2)M_1 M_1^\prime M_2 \\
&  \qquad\qquad\qquad\qquad + \lambda_*^{1/2}(2g\rho+2\sigma n_1^2)M_1^2 M_2^\prime -  \lambda_*^2 M_1 M_1^\prime M_2^\prime \bigg)\bigg|_T. \end{split} \label{Phi112} \ee 
This completes the lemma.
 \end{proof2}

Given the high degree of degeneracy at the quadratic level, in order to characterize the complete solution set of the bifurcation equation we forced to consider the cubic terms.  These we compute in the next lemma.

\begin{lemma}\emph{(Cubic Terms)} Let $\Theta_{ijk\ell}$ be given as in \eqref{defThetaijkl}.  Then 
\[\Theta_{1112},\, \Theta_{1121},\, \Theta_{1211},\, \Theta_{2111},\, \Theta_{2221},\, \Theta_{2212},\, \Theta_{2122},\, \Theta_{1222} = 0.\]
The remaining coefficients are generally nonvanishing; $\Theta_{iiii}$ is given by equation \eqref{Thetaiiii}, while $\Theta_{iijj}$ is given by \eqref{Thetaiijj}, for $i,j = 1,2$ and $i \neq j$.  \label{cubictermslemma} \end{lemma}

\begin{proof2} As in the previous lemma, we begin with an explicit calculation.  Differentiating once more in $w$ the expressions for $\mathcal{F}_{ww}(\lambda, w)$ found in \textsc{Lemma \ref{quadratictermslemma}}, we find, for $\phi, \psi, \theta \in X$ with $d(\phi) = d(\psi) = d(\theta) = 0$, 
\begin{eqnarray*}
\left\langle \left \langle \mathcal{F}_{1www}(\lambda, w)\phi, \psi \right\rangle, \theta \right \rangle & = & 2 \phi_{pp} \psi_q \theta_q + 2 \phi_q \psi_q \theta_{pp} + 2 \phi_q \psi_{pp} \theta_q  +2 \phi_{qq} \psi_p \theta_p \\
& & + 2\phi_p \psi_{qq} \theta_p + 2\phi_p \psi_p \theta_{qq} -2\phi_{pq} \psi_{p} \theta_q  -2 \phi_{pq} \psi_q \theta_p   \\
& &  -2 \phi_p \psi_{pq} \theta_q -2\phi_q -2 \phi_p \psi_{pq} \theta_q -2 \phi_p \psi_q \theta_{pq} - 2\phi_q \psi_p \theta_{pq} \\
& & -6g \rho_p (H_p + w_p)(\phi_p \psi_p \theta + \phi_p \psi \theta_p + \phi \psi_p \theta_p) \\
& & -6g \rho_p (H+w -d(H) -d(w)) \phi_p \psi_p \theta_p + 6 \beta(-p) \phi_p \psi_p \theta_p,\end{eqnarray*} 
and, on the boundary, 
\begin{eqnarray*}
\left\langle \left \langle \mathcal{F}_{2www}(\lambda, w)\phi, \psi \right\rangle, \theta \right \rangle & = & 2\phi_p \psi_p (2\sigma \kappa^\prime[w] \theta + 2g\rho \theta ) + 2\psi_p \theta_p ( 2\sigma \kappa^\prime[w] \phi + 2g\rho \phi) \\
& & + 2\phi_p \theta_p (2\sigma \kappa^\prime[w] \psi + 2g\rho \psi) +4\sigma (H_p+w_p) \phi_p \langle \kappa^{\prime\prime}[w] \theta, \psi \rangle \\
& & + 4\sigma(H_p+w_p) \psi_p \langle \kappa^{\prime\prime}[w] \theta, \phi \rangle +4\sigma(H_p+w_p) \theta_p \langle \kappa^{\prime\prime}[w] \phi, \psi \rangle \\
& & + 2\sigma(H_p+w_p)^2 \langle \langle \kappa^{\prime\prime\prime}[w] \phi, \psi \rangle, \theta \rangle,\end{eqnarray*}
where 
\[\langle \langle \kappa^{\prime\prime\prime}[w] \phi, \psi \rangle, \theta \rangle = -15\frac{ (w_q\phi_q \psi_q \theta_q)_q }{(1+w_q^2)^{7/2}} + 105 \frac{w_q^2 w_{qq} \phi_q \psi_q \theta_q}{(1+w_q^2)^{9/2}}.\]
Note that this last line implies $\kappa^{\prime\prime\prime}[0] = 0$.  

Evaluating each of these at $(\lambda_*,0)$, we compute  
\be \begin{split}
\left\langle \left\langle \mathcal{F}_{1www}(\lambda_*,0) \phi, \psi \right\rangle, \theta \right\rangle & =  
2 \phi_{pp} \psi_q \theta_q + 2 \phi_q \psi_q \theta_{pp} + 2 \phi_q \psi_{pp} \theta_q  +2 \phi_{qq} \psi_p \theta_p \\
& \qquad + 2\phi_p \psi_{qq} \theta_p + 2\phi_p \psi_p \theta_{qq} -2\phi_{pq} \psi_{p} \theta_q  -2 \phi_{pq} \psi_q \theta_p   \\
&  \qquad \qquad -2 \phi_p \psi_{pq} \theta_q -2\phi_q -2 \phi_p \psi_{pq} \theta_q -2 \phi_p \psi_q \theta_{pq}  \\
& \qquad \qquad \qquad -6g \rho_p H_p (\phi_p \psi_p \theta + \phi_p \psi \theta_p + \phi \psi_p \theta_p) \\
& \qquad \qquad \qquad \qquad -6g \rho_p (H-d(H)) \phi_p \psi_p \theta_p  \\
& \qquad\qquad\qquad\qquad\qquad+ 6 \beta(-p) \phi_p \psi_p \theta_p  - 2\phi_q \psi_p \theta_{pq}, \label{F1www} \end{split} \ee
\be \begin{split} 
\left\langle \left\langle \mathcal{F}_{2www}(\lambda_*,0) \phi, \psi \right\rangle, \theta \right\rangle & =  2\phi_p \psi_p (-2\sigma \theta_{qq} + 2g\rho \theta ) + 2\psi_p \theta_p ( -2\sigma \phi_{qq} + 2g\rho \phi) \\
&  \qquad + 2\phi_p \theta_p (-2\sigma \psi_{qq} + 2g\rho \psi) \\
& \qquad \qquad +12\sigma H_p \phi_p (\psi_{qq} \theta_q + \psi_q \theta_{qq})  \\
& \qquad \qquad\qquad + 12\sigma H_p \psi_p (\phi_{qq} \theta_q + \phi_q \theta_{qq}) \\
& \qquad\qquad\qquad\qquad+ 12 \sigma H_p \theta_p (\phi_{qq} \psi_q + \phi_q \psi_{qq} ). \end{split}
\label{F2www} \ee
We remark again that, owing to the regularity of $\mathcal{F}$ --- and as \eqref{F1www}-\eqref{F2www} make explicit --- we have $\Theta_{ijk\ell} = \Theta_{i \varpi(j)\varpi(k)\varpi(\ell)}$, where $\varpi$ is any permutation on $\{1,2\}$, and $i,j,k,\ell =1,2$.    

To compute the inner product $(\cdot, \cdot)_Y$, we first calculate that, by \eqref{F1www}, 
\be \begin{split}
\frac{1}{2} \left( a^3 \phi_i , \left \langle \left \langle \mathcal{F}_{1www}(\lambda_*,0) \phi_j, \phi_k \right \rangle, \phi_\ell \right\rangle \right)_{L^2(R)} & =  n_k n_\ell C_{ijkl}^{(1)}\bigg[ \left( a^3 M_i M_j^{\prime\prime}, M_k M_\ell\right)_{L^2}   \\ 
& \qquad - \left(a^3 M_i M_j^\prime, \left(M_k M_\ell \right)^\prime \right)_{L^2} \bigg]  \\ 
& + n_j n_\ell C_{ijk\ell}^{(2)} \bigg[ \left( a^3 M_i M_k^{\prime\prime},M_j, M_\ell \right)_{L^2} \\
&  \qquad-  \left(a^3 M_i  M_k^\prime, \left(M_j, M_\ell\right)^\prime \right)_{L^2} \bigg] \\
&  + n_j n_k C_{ijk\ell}^{(3)} \bigg[ \left( a^3 M_i M_\ell^{\prime\prime}, M_j M_k \right)_{L^2} \\
& \qquad - \left(a^3 M_i M_\ell^\prime, \left(M_j M_k \right)^\prime \right)_{L^2}  \bigg]  \\
&  - C_{ijk\ell}^{(4)} \bigg[ n_j^2 \left( a^3 M_i M_j, M_k^\prime M_\ell^\prime \right)_{L^2} \\
&  \qquad +n_k^2 \left( a^3 M_i M_j^\prime, M_k M_\ell^\prime \right)_{L^2} \\
& \qquad +n_\ell^2 \left( a^3 M_i M_j^\prime, M_k^\prime M_\ell \right)_{L^2} \\
&  \qquad +3g \left( \rho_p a^2 M_i M_j^\prime, \left( M_k M_\ell \right)^\prime \right)_{L^2} \\
& \qquad +3g\left( \rho_p a^2 M_i M_j, M_k^\prime M_\ell^\prime \right)_{L^2} \\
& \qquad +3 \left( a^6 H_{pp} M_i M_j^\prime, M_k^\prime M_\ell^\prime \right)_{L^2} \bigg],   \end{split} \label{interThetaijkl} \ee
where the constants above are defined by 
\be  
\left\{ \begin{split}
C_{ijk\ell}^{(1)} & :=  \int_0^{2\pi} \cos{(n_i q)} \cos{(n_j q)} \sin{(n_k q)} \sin{(n_\ell q)} dq, \\
C_{ijk\ell}^{(2)} & :=  \int_0^{2\pi} \cos{(n_i q)} \sin{(n_j q)} \cos{(n_k q)} \sin{(n_\ell q)} dq, \\
C_{ijk\ell}^{(3)} & :=  \int_0^{2\pi} \cos{(n_i q)} \sin{(n_j q)} \sin{(n_k q)} \cos{(n_\ell q)} dq, \\
C_{ijk\ell}^{(4)} & :=  \int_0^{2\pi} \cos{(n_i q)} \cos{(n_j q)} \cos{(n_k q)} \cos{(n_\ell q)} dq.
\end{split} \right. \label{Thetaijklcoefficients} \ee
Likewise, on the top we find
\be
\frac{1}{12} \left( a^2 \phi_i , \left \langle \left \langle \mathcal{F}_{2www}(\lambda_*,0) \phi_j, \phi_k \right \rangle, \phi_\ell \right\rangle \right)_{L^2(T)}  =  C_{ijk\ell}^{(4)} \left(a^2 M_i M_j^\prime M_k^\prime M_\ell^\prime \right)\bigg|_T\label{boundaryThetaijkl}.
\ee
Here we have used the boundary condition on $T$ to simplify.  

As we saw in the analysis of the quadratic terms, the symmetries of the underlying equation will be expressed through the quantities in \eqref{Thetaijklcoefficients}.  In particular, we can check that, for $i,j,k,\ell$ with an odd number of ones and twos, all of the coefficients vanish.   Hence $\Theta_{ijk\ell} = 0$ in this case.  If there are an even number, then by the permutation property observed above,  it suffices to consider two types of indices: (i) $\Theta_{iiii}$ and (ii) $\Theta_{iijj}$, for $i,j = 1,2$ with $i \neq j$.  A truly elementary calculation tells us
\[ C_{iiii}^{(1)} = \frac{\pi}{4},~C_{iijj}^{(1)} = \frac{\pi}{2},~C_{iiii}^{(4)} = \frac{3\pi}{4},~C_{iijj}^{(4)} = \frac{\pi}{2},\]
while $C_{iiii}^{(2)}, C_{iijj}^{(2)}, C_{iiii}^{(3)}, C_{iijj}^{(3)} = 0$, for $i,j = 1,2$ with $i \neq j$.     Using the information we obtained from \eqref{interThetaijkl}--\eqref{boundaryThetaijkl} we find that, for such $i,\,j$, 
\be \begin{split}
\frac{1}{\pi} \Theta_{iiii} & =  \frac{1}{12} n_i^2 \left[ \left(a^3 M_i^3, M_i^{\prime\prime}\right)_{L^2} - 2\left(a^3 M_i^2, \left(M_i^\prime\right)^2\right)_{L^2} \right] \\
&\qquad -\frac{3}{4}  \bigg[n_i^2\left(a^3 M_i^2, \left(M_i^\prime\right)^2 \right)_{L^2} + 3g\left(\rho_p a^2 M_i^2, (M_i^\prime)^2\right)_{L^2} \\
& \qquad \qquad+ \left( a^6 H_{pp} M_i, \left(M_i^\prime\right)^3\right)_{L^2} \bigg] + \frac{3}{4} \left( a^2 M_i \left(M_i^\prime\right)^3 \right)\bigg|_T
\end{split} \label{Thetaiiii} \ee 
and
\be \begin{split}
 \frac{1}{\pi} \Theta_{iijj} &= \frac{1}{6} n_j^2 \left[ \left(a^3 M_i M_i^{\prime\prime}, M_j^2\right)_{L^2} -2 \left(a^3 M_i M_i^\prime, M_j M_j^\prime\right)_{L^2} \right] \\
 & \qquad -\frac{1}{6} \bigg[ n_i^2 \left( a^3 M_i^2, \left(M_j^\prime\right)^2\right)_{L^2} +2n_j^2 \left(a^3 M_i M_i^\prime, M_j M_j^\prime \right)_{L^2} \\
 & \qquad\qquad + 6g \left(\rho_p a^2 M_i M_i^\prime, M_j M_j^\prime \right)_{L^2} + 3g\left(\rho_p a^2 M_i^2, \left(M_j^\prime\right)^2 \right)_{L^2} \\ 
 & \qquad\qquad\qquad + 3\left(a^6 H_{pp} M_i M_i^\prime, \left(M_j^\prime\right)^2 \right)_{L^2} \bigg] + \frac{1}{2} \left(a^2 M_i M_i^\prime \left(M_j^\prime\right)^2 \right)\bigg|_T. \end{split}
  \label{Thetaiijj}  \ee
  This completes the proof.  
\end{proof2}

\paragraph{Step III: Control of higher order terms.}  The next step is to use the scaling in \eqref{expandedbifurcationfunction} to show that all solutions (locally) must lie in a particular subregion of $U$.  In the case of mixed solutions, this effort is complicated by the fact that $\Phi_{ijk}$ is degenerate in the $\xi_1$-direction.  Consequently, in order to apply the standard theory we must restrict ourselves to working in sets that lie in complements of cusps surrounding the $\xi_1$-axis.

\begin{lemma} \emph{(Scaling of Solutions)} \emph{(a)} Suppose that 
\be \Phi_{112},\Phi_{121}, \Phi_{211} \neq 0. \label{nondegeneracyconda} \ee  
For each $ C > 0$, there is a constant $K > 0$ and open subregion $V \subset U$ with $(0,0,\lambda_*) \in V$ such that, if $\mathfrak{B}(\xi, \lambda) = 0$ for some $(\xi, \lambda) \in V \setminus \{ (\xi,\lambda) : |\xi_1| < C |\xi_2| \} $, then $(\xi,\lambda)$ satisfies the scaling relation $|\xi| \leq K |\lambda-\lambda_*|$.  \\

\noindent
\emph{(b)} Suppose instead that
\be  \Phi_{112}, \Phi_{121}, \Phi_{211} = 0, \qquad  0 \neq \Theta_{1111}\Theta_{2222} \neq \Theta_{1122}\Theta_{2211}. \label{nondegeneracycond} \ee
Then there exists a constant $K > 0$ and and open subset $V$ of $U$ with $(0, 0,\lambda_*) \in V$ such that, if $\mathfrak{B}(\xi, \lambda) = 0$ for some $(\xi, \lambda) \in V$, then $|\xi| \leq K |\lambda-\lambda_*|^{1/2}$. \label{scalinglemma} \end{lemma}

\begin{proof2} First consider the situation in part (a).  By contradiction suppose that, for some $C > 0$, there exists a sequence $\{(\xi^{(n)}, \lambda_n)\} \subset V \setminus \{ (\xi, \lambda) : |\xi_1| < C |\xi_2| \}$ such that  $\lambda_n \searrow \lambda_*$, $\xi^{(n)} \to 0$ as $n \to \infty$, and 
\[ \frac{|\xi^{(n)}|}{|\lambda_n - \lambda_*|} \to \infty, \qquad \textrm{as } n \to \infty. \]
We may therefore write $\lambda_n - \lambda_* = \rho_n |\xi^{(n)}| $, where $\rho_n \to 0$ as $n \to \infty$.  Then the equation for $\mathfrak{B}_1$ in \eqref{expandedbifurcationfunction} becomes, for all $n \geq 1$, 
\[ 0 = \xi_1^{(n)} \left|\xi^{(n)}\right| \rho_n \Psi_{11} + 2\xi_1^{(n)} \xi_2^{(n)} \Phi_{112} + O\left(\left|\xi^{(n)}\right|^3\right). \]
Diving through by $|\xi^{(n)}|^2$ and taking $n \to \infty$, each term above but the second clearly vanishes.  So, in order to produce a contradiction we will prove that the second term is bounded uniformly away from zero.   By construction we know $|\xi_1^{(n)}| > C|\xi_2^{(n)}|$, for $n \geq 1$. Suppose that there exists a $C^\prime > 0$ such that $C^\prime |\xi_1^{(n)}|^2 < |\xi_2^{(n)}|^2$.  Then, 
\[ \left| \frac{\xi_1^{(n)} \xi_2^{(n)}}{|\xi^{(n)}|^2} \Phi_{112} \right| \geq \frac{C}{1+C^\prime} |\Phi_{112}| > 0.\]
Thus, we are we done if we can produce such a constant $C^\prime$.

Suppose no such $C^\prime$ exists.  Then, possibly by passing to a subsequence, it must be the case that  $|\xi_2^{(n)}|/ |\xi_1^{(n)}| \to 0$ as $n \to \infty$.  By the equation for $\mathfrak{B}_2$ we have, $\forall n \geq 1$,
\[ 0 = \xi_1^{(n)} \left|\xi^{(n)}\right| \rho_n \Psi_{22} + (\xi_1^{(n)} )^2 \Phi_{211} + O\left(\left|\xi^{(n)}\right|^3\right). \] 
As before, we may divide through by $|\xi^{(n)}|^2$ and take $n \to \infty$ to find that the first term and the higher order terms vanish in the limit.  But,
\[ \lim_{n\to\infty} \left| \frac{|\xi_1^{(n)}|^2}{|\xi^{(n)}|^2} \Phi_{211} \right| = |\Phi_{211}| > 0,\]
so we have a contradiction.  This completes part (a).  

To prove (b), we shall employ a similar blow-up argument.  Now, however, as the coefficients $\Theta_{ijk\ell}$ are assumed to be non-degenerate, we will not be restricted to working in the complements of cusps.  

By contradiction suppose that there exists a sequence $\{(\xi^{(n)}, \lambda_n)\} \subset U$ with $\lambda_n \searrow \lambda_*$, $\xi^{(n)} \to 0$, and 
\[ \frac{\left|\xi^{(n)} \right|}{\left| \lambda_n - \lambda_* \right|^{1/2}} \to \infty, \qquad \textrm{as } n \to \infty.\]
We may therefore write $\lambda_n - \lambda_* = \rho_n |\xi^{(n)}|^2$, where $\rho_n \to 0$ as $n \to \infty$.  Then the bifurcation equation becomes, for $n \geq 1$, $i = 1,2$,
\be 0 = \xi_i^{(n)} \rho_n \left| \xi^{(n)} \right|^2 \Psi_{ii} + \xi_j \xi_k \xi_\ell \Theta_{ijk\ell} + O(|\xi^{(n)}|^4), \label{contrascaledbifurcationfunction} \ee
Dividing through by $|\xi^{(n)}|^3$ and taking $n \to \infty$ we see that the first term and the higher order terms vanish in the limit.  On the other hand, by hypothesis we have, for $i = 1,2$, 
\[  \left|\frac{\xi_j \xi_k \xi_\ell}{\left| \xi^{(n)} \right|^3} \Theta_{ijk\ell} \right| \geq  \inf_{|\xi| = 1}  \left| \xi_j \xi_k \xi_\ell \Theta_{ijk\ell}\right| > 0 . \]
Hence the third term is bounded above and strictly away from zero.  Having arrived at a contradiction, we conclude the lemma holds.  \end{proof2}

\begin{longremark}  Recall that, when we prove \textsc{Theorem \ref{mainresult2}}, we will be assuming \eqref{capgravnondegeneracycondition1}--\eqref{capgravnondegeneracycondition2} with $n_1 = 1$, which together imply \eqref{nondegeneracycond}. Hence, we expect the conclusions of part (b) to hold.    
\end{longremark}

Fix $C > 0$ and define the set $W$ as follows
\be W := \left \{ 
\begin{array}{ll}  V \setminus \{ (\xi,\lambda) \in U : |\xi_1| < C |\xi_2|\}, & \textrm{if \eqref{nondegeneracyconda} holds,} \\
V, & \textrm{if \eqref{nondegeneracycond} holds,}
\end{array} \right. \label{defW} \ee 
where $V$ is as in \textsc{Lemma \ref{scalinglemma}}.  
Let $W^{\pm} \subset W$ denote the sets $\{ (\xi, \lambda) \in W : \pm(\lambda - \lambda_*) > 0 \}$.   Then, working in $W^{\pm}$, we are justified in letting $\epsilon = \epsilon(\xi,\lambda) > 0$, $\theta = \theta(\xi,\lambda) \in \mathbb{R}^2$ be defined by the relations
\[  \xi = \pm \epsilon \theta, \qquad  \epsilon = \left\{ \begin{array}{ll} |\lambda - \lambda_*|, & \textrm{if \eqref{nondegeneracyconda} holds,} \\
|\lambda-\lambda_*|^{1/2}, & \textrm{if \eqref{nondegeneracycond} holds.} \end{array} \right.\]
Note that \textsc{Lemma \ref{scalinglemma}} guarantees $\|\theta\|_{L^\infty(V)}$ is bounded.  Rewriting \eqref{expandedbifurcationfunction} in terms of $\epsilon$ and $\theta$ we obtain the following system of equations for $(\theta, \epsilon)$:  
\be 0 = \widetilde{\mathfrak{B}}_i^{\pm}(\theta, \epsilon) := \left\{\begin{array}{ll} 
\pm \theta_j \Psi_{ij} + \theta_j \theta_k \Phi_{ijk} + O(\epsilon), & \textrm{if \eqref{nondegeneracyconda} holds}, \\
\pm \theta_j \Psi_{ij} + \theta_j \theta_k \theta_\ell \Theta_{ijk\ell} + O(\epsilon), & \textrm{if \eqref{nondegeneracycond} holds,} \end{array} \right.
 \label{reducedbifurcationequation} \ee
for $i, j, k, \ell = 1,2$.

We wish to say that solutions to \eqref{reducedbifurcationequation} with $\epsilon = 0$ are, in some sense, the only solutions to \eqref{expandedbifurcationfunction} in the set $W$.  To reach this conclusion requires taming the higher order terms, represented by the $O(\epsilon)$ in \eqref{reducedbifurcationequation}.  We shall do this by proving the following lemma.

\begin{lemma} \emph{(Regular Value)} Let $W \subset U$ be as in \eqref{defW}.  Then zero is a regular value of $\theta \mapsto \widetilde{\mathfrak{B}}^{\pm}(\theta,0)$ in $W$, provided that either (a) \eqref{nondegeneracyconda} holds, or (b) \eqref{nondegeneracycond} holds and we have, additionally, that 
\be \Theta_{1111}\Psi_{22} \neq \Theta_{2211}\Psi_{11}, \qquad \Theta_{2222}\Psi_{11}  \neq \Theta_{1122} \Psi_{22}. \label{regularvalueconditions} \ee \label{regularvaluelemma} \end{lemma}

\begin{proof2} Assume first that hypothesis (a) holds and say that, for some $\theta \neq 0$ in $W^\pm$, we have $\widetilde{\mathfrak{B}}^\pm(\theta, 0) = 0$.  Then, using our knowledge of the specific forms of  $\Phi_{ijk}$ and $\Psi_{ij}$ garnered from \textsc{Lemma \ref{quadratictermslemma}}--\textsc{Lemma \ref{cubictermslemma}}, we see that, from the equation for $\widetilde{\mathfrak{B}}_1^{\pm}$, 
\[ \pm \Psi_{11} + 2 \theta_2 \Phi_{112} = 0.\]
Here we have used the fact that $\theta_1 \neq 0$, since $\theta \in W$ and is nonzero.  But, the Jacobian of $\theta \mapsto \widetilde{\mathfrak{B}}^{\pm}(\theta,0)$ is easily computed to be
\be \det \left(\begin{array}{cc}
\pm\Psi_{11} + 2\theta_2 \Phi_{112} & 2\theta_1 \Phi_{112} \\
2\theta_1\Phi_{211} & \pm\Psi_{22} \end{array} \right) =  \pm(\pm\Psi_{11} + 2 \theta_2 \Phi_{112}) \Psi_{22} - 4 \theta_1^2 \Phi_{112} \Phi_{211} . \label{jacobiancalculation} \ee
We have already shown that the first term on the right-hand side vanishes, and, by hypothesis, the second term is strictly nonzero.   On the other hand, in the case where $\theta = 0$, \eqref{jacobiancalculation} shows us  the Jacobian is given by $\Psi_{11} \Psi_{22} > 0$. Hence zero is a regular value in $W$. 

Next consider what happens when the hypotheses of (b) are satisfied.  Then, using the information on the cubic terms obtained in \textsc{Lemma \ref{cubictermslemma}}, we see that
\[ \widetilde{\mathfrak{B}}^\pm(\theta, 0) = \left( \begin{array}{l}\pm \theta_1 \Psi_{11} + \theta_1 \theta_2^2 \Theta_{1122} + \theta_1^3 \Theta_{1111} \\
\pm \theta_2 \Psi_{22} + \theta_1^2 \theta_2 \Theta_{2211} + \theta_2^3 \Theta_{2222} \end{array}\right).\]  
Taking $\theta_1 = 0$ we have automatically that $\widetilde{\mathfrak{B}}_1^\pm = 0$.  To get nontrivial solutions, therefore, requires that $\theta_2$ satisfies
\be \pm \Psi_{22} + \theta_2^2 \Theta_{2222} = 0. \label{theta2eq1} \ee
Note this implies that, depending on the sign of $\Theta_{2222}$, either $\widetilde{\mathfrak{B}}^+(0,\theta_2,0)=0$ has two solutions and $\widetilde{\mathfrak{B}}^-(0,\theta_2,0)=0$ has none, or the reverse is true.  In any event, evaluating the Jacobian of $\widetilde{\mathfrak{B}}^\pm$ at these points yields
\[ \det \left( \begin{array}{cc} \pm\Psi_{11} + \theta_2^2 \Theta_{1122} & 0 \\
0 & \pm\Psi_{22} + 3\theta_2^2 \Theta_{2222}  \end{array} \right) = 2\theta_2^2 \Theta_{2222}(\pm \Psi_{11}+\theta_2^2\Theta_{1122}). \]
Thus the Jacobian vanishes if and only if $\pm \Psi_{11} + \theta_2^2 \Theta_{1122} = 0$.  Of course, if $\Theta_{1122} = 0$, this never occurs and we are done.  Otherwise, combining this with \eqref{theta2eq1} leads to  
\[  \mp \frac{\Psi_{11}}{\Theta_{1122}} = \theta_2^2 = \mp \frac{\Psi_{22}}{\Theta_{2222}},\]
which contradicts \eqref{regularvalueconditions}.  (Note that $\Theta_{1111}, \Theta_{2222} \neq 0$ by the nondegeneracy condition \eqref{nondegeneracycond}.)  Thus all zeros of $\theta \mapsto\widetilde{\mathfrak{B}}^\pm(\theta,0)$ in $W^\pm$ with $\theta_1 = 0$ are simple.  Due to  the symmetry of $\widetilde{\mathfrak{B}}$, we can simply reverse the roles of one and two above to get a similar result for the case when $\theta_2 = 0$, $\theta_1 \neq 0$.  

Finally let us consider solutions for which neither $\theta_1$ nor $\theta_2$ vanish.  Then simplifying we find
\be 0 = \left( \begin{array}{c} \pm \Psi_{11} + \theta_2^2 \Theta_{1122} + \theta_1^2 \Theta_{1111} \\
\pm \Psi_{22} + \theta_1^2 \Theta_{2211} + \theta_2^2 \Theta_{2222} 
\end{array} \right), \label{Btildeatmixedsolution} \ee
so that the Jacobian evaluated at such points comes to
\be \begin{split}
 \det \left(\frac{\partial\widetilde{\mathfrak{B}}^\pm}{\partial\theta} \right)\bigg|_{(\theta,0)} &= (\pm \Psi_{11} + \theta_2^2 \Theta_{1122} +3\theta_1^2 \Theta_{1111})( \pm \Psi_{22} + \theta_1^2 \Theta_{2211} + 3\theta_2^1 \Theta_{2222}) \\
 & \qquad -4\theta_1^2\theta_2^2 \Theta_{1122} \Theta_{2211} \\
 & =  4\theta_1^2\theta_2^2\left(\Theta_{1111}\Theta_{2222} - \Theta_{1122}\Theta_{2211} \right).
 \end{split}
 \label{jacobianmixedsolutions} \ee
But this quantity is nonvanishing for $\theta \neq 0$ by \eqref{nondegeneracycond}.  This completes the lemma. \end{proof2}
 
 Examining the proof of \textsc{Lemma \ref{regularvaluelemma}}, we can count the number of nontrivial roots of the map $\theta \mapsto \widetilde{\mathfrak{B}}^\pm(\theta,0)$ in $W$ under various assumptions.  First suppose hypothesis (b) (i.e.,  \eqref{nondegeneracycond} and \eqref{regularvalueconditions} hold.) Then, by \eqref{theta2eq1}, taking $\theta_1 = 0$ and $\Phi_{2222} > 0$, there are no nontrivial solutions to $\widetilde{\mathfrak{B}}^-(0,\theta_2, 0) = 0$, and two nontrivial solutions to $\widetilde{\mathfrak{B}}^+(0,\theta_2, 0) = 0$.  If $\Phi_{2222} < 0$, then there are no nontrivial solutions to $\widetilde{\mathfrak{B}}^+(0,\theta_2, 0) = 0$ and two for $\widetilde{\mathfrak{B}}^-(0,\theta_2, 0) = 0$.  The same holds true when we take $\theta_2 = 0$ and search for nontrivial solutions.  Thus there are a total of four possible nontrivial solutions where precisely one of $\theta_1$ and $\theta_2$ is zero.  
 
 On the other hand, \eqref{Btildeatmixedsolution} shows us that if neither $\theta_1$ nor $\theta_2$ vanish, then $(\theta_1^2, \theta_2^2)$ solve the linear system
 \[ A \left( \begin{array}{c} \theta_1^2 \\ \theta_2^2 \end{array}\right) = \mp \left(\begin{array}{c} \Psi_{11} \\ \Psi_{22} \end{array}\right), \qquad \textrm{where } A := \left( \begin{array}{cc} \Theta_{1111} & \Theta_{1122} \\ \Theta_{2211} & \Theta_{2222} \end{array}\right).\]
 Note that $A$ is invertible by assumption \eqref{nondegeneracycond}.  Therefore, one of $\widetilde{\mathfrak{B}}^+, \widetilde{\mathfrak{B}}^-$ has no such solutions with $\theta_1, \theta_2 \neq 0$, while the other has four.  In total, we find that there are eight nontrivial solutions of $\widetilde{\mathfrak{B}}^\pm(\theta, 0) = 0$ in $W$:
 \be \begin{split} \left(\theta_1 = 0,~\theta_2 = \pm |\Psi_{22}/\Theta_{2222}|^{1/2}\right), \qquad \left( \theta_1= \pm |\Psi_{11}/\Theta_{1111}|^{1/2},~ \theta_2 = 0\right), & \\ 
 \left( \theta_1=\pm  |A|^{-1}|\Theta_{2222} \Psi_{11} - \Theta_{2211} \Psi_{22}|^{1/2},~ \theta_2 = \pm  |A|^{-1}|\Theta_{1111}\Psi_{22}-\Theta_{1122} \Psi_{11}|^{1/2} \right)&. \label{characterizationsolutions} \end{split} \ee

Suppose now that the situation in (a) occurs, that is,  $\Phi_{112}, \Phi_{121}, \Phi_{211} \neq 0$.  Then if we let $\theta_1 = 0$, we find that $\theta_2$ must be zero as well (this follows both from the bifurcation equation and, at a more basic level, from the definition of $W$.)  Similarly, there are no nontrivial solutions with $\theta_2 = 0$.  However, if we suppose that $\theta_1, \theta_2 \neq 0$, then we can solve $\widetilde{\mathfrak{B}}_1^\pm(\theta, \epsilon)$ for $\theta_2$ to find
\[ \theta_2 = \mp \frac{\Psi_{11}}{2\Phi_{112}},\]
hence, $\widetilde{\mathfrak{B}}_2^\pm(\theta,0) = 0$ if and only if
\[ \theta_1^2 = \mp \frac{\theta_2 \Psi_{22}}{\Phi_{211}} = \frac{\Psi_{11}\Psi_{22}}{2\Phi_{112}\Phi_{221}}.\]  
Since $\Psi_{11}$ and $\Psi_{22}$ are of the same sign, we see that the above equation has two solutions or zero, depending on whether or not $\Phi_{112}$ and $\Phi_{211}$ have the same sign.  If $\Phi_{112} \Phi_{211} < 0$, therefore, we find that there are no nontrivial solutions to the bifurcation equation, while if $\Phi_{112}\Phi_{211} > 0$ then there are two nontrivial roots of $\widetilde{\mathfrak{B}}^{\pm}(\theta, 0) = 0$:
\[  \left(\theta_1 = \sqrt{\frac{\Psi_{11}\Psi_{22}}{2\Phi_{112}\Phi_{221}}},~ \theta_2 = \mp \frac{\Psi_{11}}{2\Phi_{112}}\right), \qquad\left(\theta_1 = -\sqrt{\frac{\Psi_{11}\Psi_{22}}{2\Phi_{112}\Phi_{221}}},~ \theta_2 = \mp \frac{\Psi_{11}}{2\Phi_{112}}\right) .\]

Putting these observation together with \textsc{Lemma \ref{scalinglemma}} and \textsc{Lemma \ref{regularvaluelemma}}, we obtain the following result characterizing the solution set of the bifurcation equation under various assumptions on $\Phi_{ijk}, \Theta_{ijk\ell}$.
 
 \begin{figure}[h] \setlength{\unitlength}{5cm} 
 \centering
\includegraphics[scale=0.5]{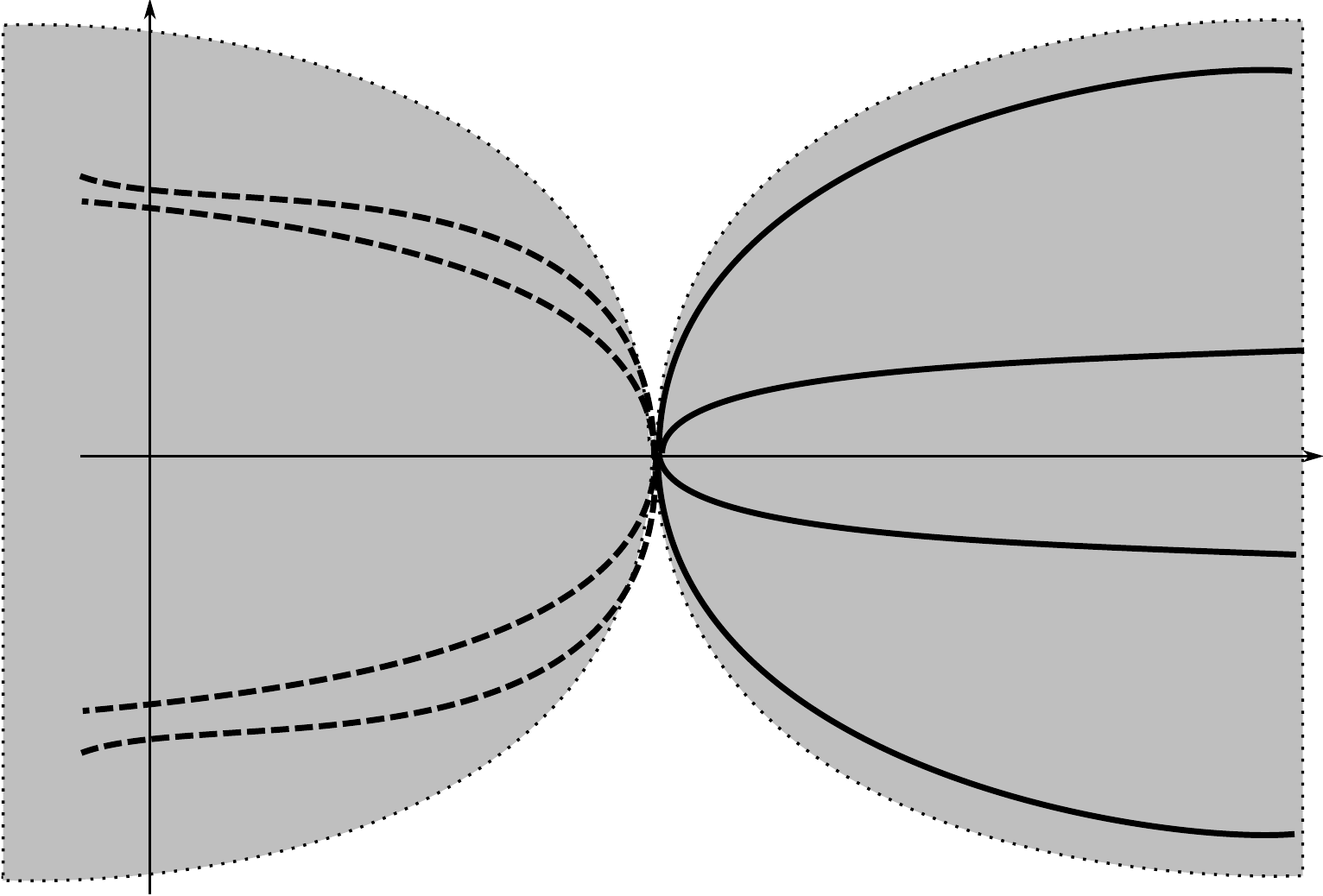}
\put(-0.65,1.05){$|\xi| = K\sqrt{\lambda - \lambda_*}$}
\put(-1.15,0.00){$|\xi| = K\sqrt{\lambda_*-\lambda}$}
\put(0.02,0.50){$\lambda$}
\put(-1.35,1.05){$\xi$}
\put(-0.15,0.70){$W^+$}
\put(-1.35,0.70){$W^-$}
\put(-1.0,0.55){$(\lambda_*, 0)$}
\caption{A typical bifurcation diagram.  In this case we have $\Phi_{112}, \Phi_{221} =0$, $\Theta_{1111}, \Theta_{2222} > 0$ and $\Theta_{1111}\Theta_{2222} < \Theta_{1122}\Theta_{2211}$.  The two solid line pitchforks on the right consist of pure solutions, while the two dashed pitchforks on the left are mixed solution curves.  The shaded regions are $W^+$ and $W^-$.}
\end{figure}
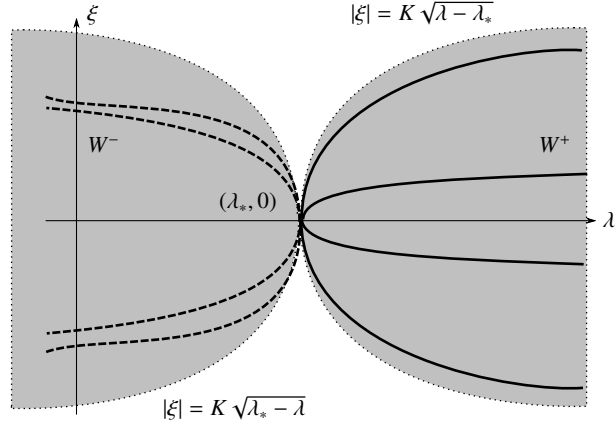

\begin{theorem} \emph{(Bifurcation structure at double eigenvalues)} Suppose that at $\lambda_* \geq -2B_{\textrm{min}}+\epsilon_0$ the null space of $\mathcal{F}_{w}(\lambda_*,0)$ is spanned by $\phi_i = M_i(p) \cos{(n_i q)}$, $i = 1, 2$,  where $0 < n_1 < n_2$.  Let $\Phi_{ij}$, $\Theta_{ijk}$ be defined as above, for $i,j,k = 1,2$.
\begin{itemize}
\item[\emph{(a)}] If $\Phi_{ij}$ and $\Theta_{ijk\ell}$ satisfy \eqref{nondegeneracycond} and \eqref{regularvalueconditions}, then there exist four $C^1$-curves, $\{ \mathcal{C}_{\mathrm{loc},i}^\prime \}_{i=1}^4$, of non-laminar solutions of the height equation bifurcating from $\mathcal{T}$ at $\lambda_*$.  In a sufficiently small neighborhood of $(Q(\lambda_*),H_*)$ in $X$, moreover, these constitute all nontrivial solutions.  Two of these curves consist of pure solutions --- one of minimal period $2\pi/n_1$, one of minimal period $2\pi/n_2$ --- and two are of mixed solutions.  \\

\item[\emph{(b)}] If instead we have that \eqref{nondegeneracyconda} holds, then there exist a $C^1$-curve of minimal period $2\pi/n_2$ solutions to the height equation bifurcating from $\mathcal{T}$ at $\lambda_*$.  If additionally we have $\Phi_{112}\Phi_{211} > 0$, then there are two $C^1$-curves of mixed-type, non-laminar solutions bifurcating from $\mathcal{T}$ at $\lambda_*$.  In a sufficiently small neighborhood of the kind described in \textsc{Lemma \ref{scalinglemma}}, these comprise all nontrivial solutions.
\end{itemize} \label{bifurcationstructuretheorem} \end{theorem}

\begin{proof2} Consider the scenario in (a).  As we have seen in the preceding lemmas, this means that the bifurcation equation is equivalent to solving \eqref{reducedbifurcationequation} near $(\theta = 0, \epsilon = 0)$.  In light of \textsc{Lemma \ref{regularvaluelemma}}, moreover, the Implicit Function theorem tells us that to each solution of $\widetilde{\mathfrak{B}}^\pm(\theta, 0) = 0$, there corresponds a curve of solutions to the bifurcation equation.  The arguments of the preceding paragraphs tell us that there are precisely four of these, two of mixed type (i.e. with $\theta_1, \theta_2 \neq 0$) and two of pure type --- one with $(\theta_1 = 0, \theta_2 \neq 0)$ and one with $(\theta_1 \neq 0, \theta_2 = 0)$.  The union of these curves gives the complete solution set of the bifurcation equation in set in a sufficiently small neighborhood of $(\lambda_*,0,0)$.   

Next let us look at the situation in (b).   Fix $C > 0$ and let $W$ be defined as before.  Working within $W$ 
we know that there are no pure solutions, but there are two of mixed type if and only if $\Phi_{112}$ and $\Phi_{211}$ share the same sign.  If this does occur, then as we argued in the previous paragraph, \textsc{Lemma \ref{regularvaluelemma}} and the Implicit Function theorem together give the existence of three curves of mixed solutions.

However, these clearly do not give the entire solution set--- Indeed, we know there must be a curve of purely $2\pi/n_2$-periodic solutions.  This follows, as in \cite{wahlen2006capgrav}, by noting that restricting the solution space $X$ to $2\pi/n_2$-periodic solutions reduces to one the dimension of the null space of $\mathcal{F}_w(\lambda_*,0)$.  Hence the existence of such a curve follows from a routine application of Crandall-Rabinowitz.  Note that this would correspond to taking $\xi_1 = 0$, and thus lies outside of $W$ no matter the choice of $C$. \end{proof2}

We are finally ready to assemble the results of the previous steps into a proof of the main theorem.  This amounts to nothing more than using the sufficient conditions derived in \textsc{Lemma \ref{quadratictermslemma}} and \textsc{Lemma \ref{cubictermslemma}} in conjunction with \textsc{Theorem \ref{bifurcationstructuretheorem}}.  

\begin{proofof}{Theorem \ref{mainresult2}.}  Recall that, by hypothesis, we have \eqref{capgravnondegeneracycondition1} and \eqref{capgravnondegeneracycondition2}.  Since we are taking $n_1 = 1$ and $n_2 \geq 3$, by \textsc{Lemma \ref{lineartermslemma}} this in turn implies \eqref{nondegeneracycond} and \eqref{regularvalueconditions} each hold.   Applying part (a) of \textsc{Theorem \ref{bifurcationstructuretheorem}}, therefore, we conclude that there exists four $C^1$-curves of solutions, $\{\mathcal{C}_{i,\textrm{loc}}^\prime\}_{i=1}^4$, of the reformulated problem bifurcating from $\mathcal{T}$ at $(Q_*, H_*)$.  Likewise, by \textsc{Lemma \ref{equivalencelemma}}, there exists four $C^1$-curves, $\{ \mathcal{C}_{i,\textrm{loc}} \}_{i=1}^4$, of solutions to the original problem \eqref{incompress}--\eqref{euler2}.  All that remains to prove are the various nodal properties of these solutions. 

Property (i) follows from the symmetries of the height equation, as we have argued before.  Now, without loss of generality, we may suppose that each $(Q,h) \in \mathcal{C}_{i,\textrm{loc}}^\prime$ can be written
\[ h(q,p) = H(p; \lambda) + \epsilon \phi_i(p) \cos{(n_i q)} + o(\epsilon), \qquad \textrm{in X}, \]
for $i = 1,2$. In other words, we are choosing the labeling so that $\mathcal{C}_{1, \textrm{loc}}^\prime$ and $\mathcal{C}_{2,\textrm{loc}}^\prime$ are the curves of pure solutions, while $\mathcal{C}_{3,\textrm{loc}}^\prime$ and $\mathcal{C}_{4,\textrm{loc}}^\prime$ are the mixed solutions.  By taking $\epsilon$ sufficiently small (or equivalently, by restricting ourselves to a sufficiently small neighborhood of the bifurcation point), we have that for $(Q,h) \in \mathcal{C}_{i,\textrm{loc}}^\prime$, $h$ has minimal period $2\pi/n_i$ in $q$. Moreover, since
\[ h_q(q,p) = - \epsilon n_i \phi_p(p) \sin{(n_i q)} + o(\epsilon), \qquad \textrm{in } X,\]
we see that, again for $\epsilon > 0$ sufficiently small, $h_q(q,p)$ has the opposite sign of $\sin{(n_i q)}$ for each $(q,p) \in \overline{R}$.  But, recall that $h_q = v/(c-u)$, which describes the slope of the vector field in the Eulerian formulation.  It follows that, if we restrict $\mathcal{C}_{i, \textrm{loc}}$ to be solutions with $\epsilon > 0$ and sufficiently small, then the wave profile, indeed all of the streamlines above the bed, are monotonic between crestline and troughline.  This implies properties (ii)--(iii).  \end{proofof}

\begin{longremark}  Elements of the mixed solution curves $\mathcal{C}_{3,\textrm{loc}}^\prime$ and $\mathcal{C}_{4,\textrm{loc}}^\prime$ can likewise be characterized
\[ h(q,p) = H(p;\lambda) + \epsilon \left( \xi_1 \phi_1(p) \cos{q} + \xi_2 \phi_2(p) \cos{(n_2 q)} \right) + o(\epsilon), \qquad \textrm{in } X,\] 
where $\xi_1, \xi_2 \in \mathbb{R}^\times$ are found from \eqref{characterizationsolutions}.  \end{longremark}

\section{Global bifurcation from simple eigenvalues} \label{globalbifurcationsection}

The next step in our program is to continue the local curves of \textsc{Theorem \ref{capgravsimplelocalbifurcation}} and  \textsc{Theorem \ref{mainresult2}} globally.  In this section we shall address the former (i.e., the simple bifurcation case), while the latter will be the subject of section \ref{analyticitysection}.

 As in previous section, we denote
 \[ X := \{ h \in C_{\textrm{per}}^{3+\alpha}(\overline{R}) : h = 0 \textrm{ on } B\}, \qquad
 Y = Y_1 \times Y_2 := C_{\textrm{per}}^{1+\alpha}(\overline{R}) \times C_{\textrm{per}}^{1+\alpha}(T),\]
 where the subscript ``per'' indicates $2\pi$-periodicity and evenness in $q$, $R$ is the rectangle $(0,2\pi)\times (p_0,0)$ and $T = [0,2\pi]\times \{p = 0\}$.  In this section, we let the nonlinear operator $\mathcal{G} = (\mathcal{G}_1, \mathcal{G}_2) : \mathbb{R}\times X \to Y$ be given by
 \begin{align}
 \mathcal{G}_1(h) & :=  (1+ h_q^2)h_{pp} + h_{qq} h_p^2 - 2 h_q h_p h_{pq} -g(h-d(h))\rho_p h_p^3 + h_p^3 \beta(-p) \label{capgravdefG1} \\
 \mathcal{G}_2(Q,h) & :=  \bigg(1+h_q^2 +h_p^2(2\sigma\kappa[h] +2g\rho h - Q)\bigg)\bigg|_{p = 0}. \label{capgravdefG2} \end{align}
Recall that $\kappa$ denotes the mean curvature: 
\[ \kappa[h] := -\frac{h_{qq}}{\left(1+h_q^2\right)^{3/2}}.\]
Note that, by construction, the laminar flow solutions $H(\cdot; Q)$ found in \textsc{Lemma \ref{laminarflowlemma}} satisfy 
\[\mathcal{G}(H(\cdot; \lambda), Q(\lambda))) = 0, \qquad \textrm{for all } \lambda \geq -2B_{\min{}}+\epsilon_0.\]

\begin{figure} \setlength{\unitlength}{5cm} 
 \centering
\includegraphics[scale=0.5]{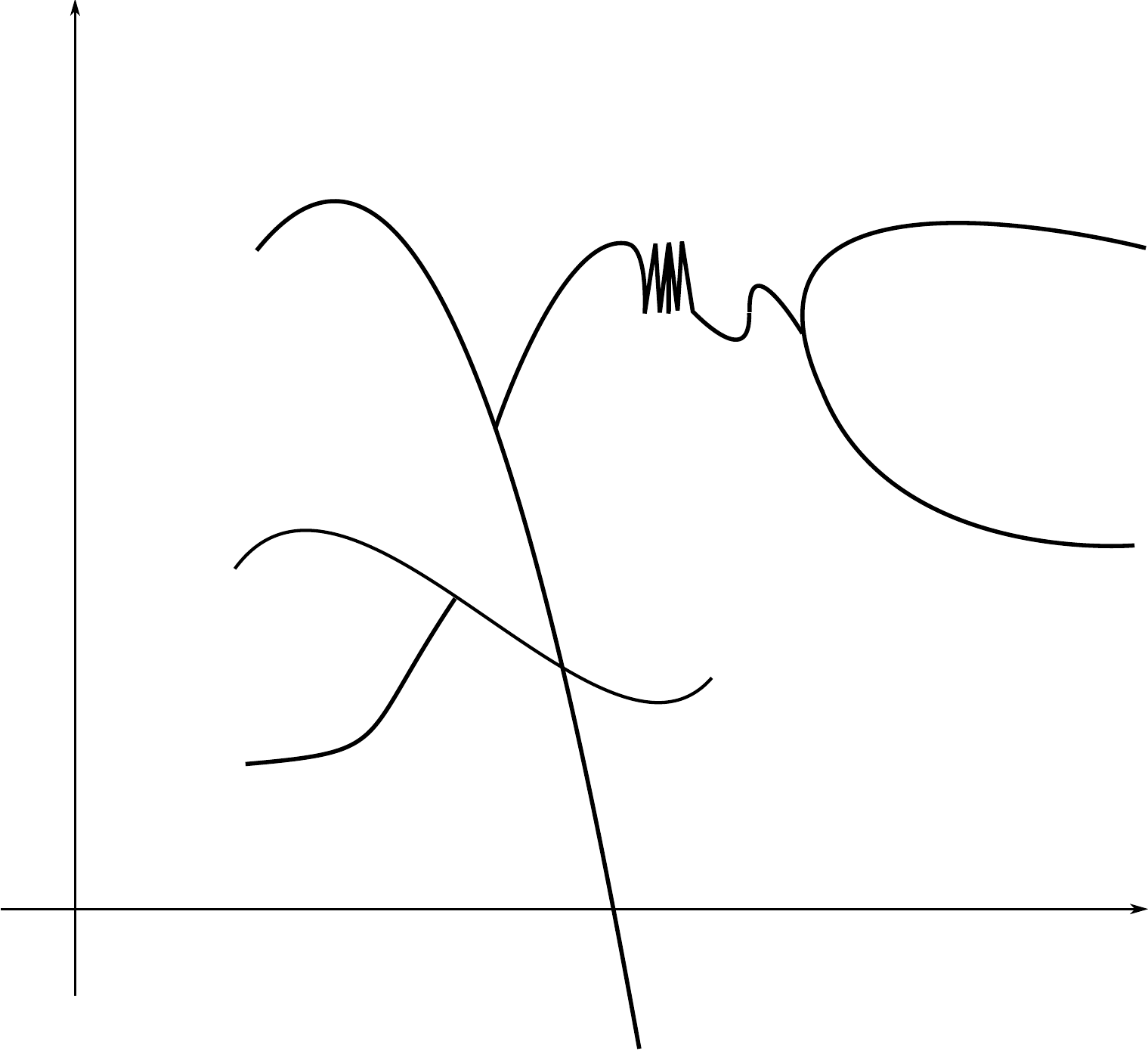}
\put(-0.42,0.70){$\mathcal{C}_{\delta}^\prime$}
\put(-0.90,0.30){$\mathcal{C}_{\mathrm{loc}}^\prime$}
\put(-0.1,0.22){$\mathcal{T}$}
\put(-1.40,1.35){$X$}
\put(-1.05,1.20){$\mathcal{O}_\delta$}
\put(-.7,0.22){$(Q^*,H^*)$}
\caption{} \label{rabinowitzdiagram}
\end{figure}

The main engine for our continuation argument will be Kielh\"ofer degree theory, a generalization of the classical theory of Leray-Schauder (cf. \cite{kielhofer1985multiple}).  The basic outline for this type of argument is by now well established (see, e.g., \cite{constantin2004exact,healey1998global,rabinowitz1971some,walsh2009stratified}).  Fixing $\delta > 0$, put 
\[ \mathcal{O}_{\delta} := \left\{ (Q,h) \in \mathbb{R}\times X : h_p > \delta \textrm{ in } \overline{R},~Q-2\sigma \kappa[h]-2g\rho h  > \delta~ \right\}.\]
Likewise, denote
\[S_\delta := \textrm{closure in } \mathbb{R}\times X \textrm{ of } \left\{(Q,h) \in \mathcal{O}_{\delta} : \mathcal{G}(Q,h) = 0, ~ h_q \nequiv 0\right\},\]
 and let $C_\delta^\prime$ be the component of $S_\delta$ that contains the point $(Q_*, H_*)$, where $Q_* := Q(\lambda_*)$, $H_* := H(\cdot; \lambda_*)$ is the point where Crandall-Rabinowitz bifurcation occurs.  We note that this is valid since the laminar flows all have zero curvature, hence $(Q^*, H^*) \in \cap_\delta \mathcal{O}_\delta$.  Thus $C_\delta^\prime$ contains the local curve $C_{\textrm{loc}}^\prime$.  A typical bifurcation diagram for this scenario is shown in Figure \ref{rabinowitzdiagram}.  Note that, the continuum $\mathcal{C}_\delta^\prime$ is connected, but need not be path-connected in general.  We shall discuss this important fact more in the next section.  
 
The main result of this section is the following.
\begin{samepage}
\begin{theorem} \emph{(Global Bifurcation from Simple Eigenvalues)} Let $\delta > 0$ be given.  One of the following alternatives must hold. 
\begin{indentenum}
 \item[\emph{(i)}] $\mathcal{C}_\delta^\prime$ is unbounded in $\mathbb{R} \times X$.
 \item[\emph{(ii)}] $\mathcal{C}_\delta^\prime$ contains another trivial point $(Q, H) \in \mathcal{T}$, with $Q \neq Q^*$.
 \item[\emph{(iii)}] $\mathcal{C}_\delta^\prime$ contains a point $(Q,h) \in \partial \mathcal{O}_{\delta}$. \end{indentenum} 
\label{capgravglobalbifurcationtheorem}
\end{theorem}
\end{samepage}
The majority of the work in proving \textsc{Theorem \ref{capgravglobalbifurcationtheorem}} is showing the admissibility of $\mathcal{G}$ in the sense of Kielh\"ofer degree.  This will, of course, hinge on being able to derive suitable \emph{a priori} estimates for the linearized operators --- a task that is difficult in two respects.  First, the addition of the nonlocal operator $d$ to the interior equation renders (at least initially) the usual elliptic regularity results inapplicable.  The method for resolving this problem, first presented \cite{walsh2009stratified}, is to ``freeze'' $d$ to a constant, apply regularity theory to the resulting elliptic problem, then argue back to estimates on the unfrozen problem.

This leads us to the second difficulty:  because of the higher order (tangential) derivatives on $T$, even the frozen elliptic problem lies outside the domain of standard Schauder theory.  Though slightly more obscure, this type of boundary condition is usually referred to as \emph{Venttsel-type}, and has received some attention in the literature (see, e.g. \cite{luo1991linear,luo1994quasilinear}).  In the interest of self-containedness, we present here one of the more basic linear result, as that is all that we shall require in this section.  This is done by fleshing out a remark in \cite{krylov1996lectures} on the solvability of such problems in the constant coefficient case.

\begin{lemma} Fix $n \geq 2$ and let $L$ be a second-order, linear, uniformly elliptic differential operator on $\mathbb{R}^n_{+}$,
\[ L := \sum_{i,j=1}^n a_{ij}(x) \partial_{x_i}\partial_{x_j} + \sum_{i=1}^n b_i(x)\partial_{x_i} + c(x), \qquad x \in \mathbb{R}^n_{+};\]
$M$ a homogeneous, second-order, linear, uniformly elliptic differential operator on $\partial\mathbb{R}_+^{n}$,
\[ M := \sum_{i,j=1}^{n-1} \mathfrak{a}_{ij}(x^\prime) \partial_{x_i} \partial_{x_j}, \qquad x^\prime \in \partial\mathbb{R}_+^{n};\]
and $N$ a first-order, linear, oblique differential operator on $\partial\mathbb{R}^n_{+}$,
\[ N := \sum_{i=1}^n \mathfrak{b}_i(x^\prime) \partial_{x_i}, \qquad x^\prime \in \partial\mathbb{R}^n_{+}, ~|\mathfrak{b}_n| \geq \beta > 0.\]  Suppose that 
\[ \|a_{ij}\|_{C^{2+\alpha}(\mathbb{R}_+^n)},~\|b_i,~c\|_{C^{\alpha}(\mathbb{R}_+^n)}, ~ \|\mathfrak{a}_{ij},~\mathfrak{b}_i \|_{C^\alpha(\partial \mathbb{R}_+^n)} < K\]
and let $\kappa_L,~\kappa_M$ denote the ellipticity constants of $L$ and $M$ respectively.  Then there exists $C > 0$ (depending on $K, \beta, n, \kappa_L$ and $\kappa_M$) such that, for all $u \in C^{2+\alpha}(\overline{\mathbb{R}_+^n})$,
\be \|u\|_{C^{2+\alpha}(\overline{\mathbb{R}^n_+})} \leq C\left( \|u\|_{C^\alpha(\overline{\mathbb{R}_+^n})} + \|Lu\|_{C^\alpha(\overline{\mathbb{R}_+^n})} +\|(M+N) u\|_{C^\alpha(\partial\mathbb{R}_+^n)} \right).\label{modschauder} \ee\label{estimateslemma} \end{lemma}

\begin{proof2} First note that, without loss of generality, we may take $b_i, c \equiv 0$.   To see this, first suppose we had the lemma for homogeneous $L$. Then, for general $L$, we can treat the lower-order terms as a right-hand side and apply \eqref{modschauder} to find,
\begin{align*}  \|u\|_{C^{2+\alpha}(\overline{\mathbb{R}^n_+})} &\leq C\bigg( \|u\|_{C^\alpha(\overline{\mathbb{R}_+^n})} +\|Lu\|_{C^{\alpha}(\overline{\mathbb{R}_+^n})} + \sum_i \|b_i\|_{C^{\alpha}(\overline{\mathbb{R}^n_+})} \|u\|_{C^{1+\alpha}(\overline{\mathbb{R}_+^n})} \\ 
&  \qquad+ \|c\|_{C^{\alpha}(\overline{\mathbb{R}_+^n})}\|u\|_{C^\alpha(\overline{\mathbb{R}_+^n})} +\|(M+N) u\|_{C^\alpha(\partial\mathbb{R}_+^n)} \bigg). \end{align*}
Next, interpolate the $C^{1+\alpha}$-norm on the right-hand side between $C^\alpha$ and $C^{2+\alpha}$ and absorb the latter term into the left-hand side to get 
\begin{align*}  \|u\|_{C^{2+\alpha}(\overline{\mathbb{R}^n_+})} &\leq \frac{1}{2} \|u\|_{C^{2+\alpha}(\overline{\mathbb{R}_+^n})} + C\left( \|u\|_{C^\alpha(\overline{\mathbb{R}_+^n})} +\|Lu\|_{C^{\alpha}(\overline{\mathbb{R}_+^n})} +\|(M+N) u\|_{C^\alpha(\partial\mathbb{R}_+^n)} \right). \\
& \leq  C\left( \|u\|_{C^\alpha(\overline{\mathbb{R}_+^n})} +\|Lu\|_{C^{\alpha}(\overline{\mathbb{R}_+^n})} +\|(M+N) u\|_{C^\alpha(\partial\mathbb{R}_+^n)} \right),  \end{align*}
where, as required, $C = C(K,\beta, n, \kappa_L, \kappa_M)$.  So, indeed, it suffices to assume that $L$ is homogeneous.
  
Now let us introduce a few notational conventions.  Let the Poisson map, $\mathcal{P}:C^{\alpha}(\mathbb{R}_+^n) \times C^{2+\alpha}(\partial \mathbb{R}_+^n) \to C^{2+\alpha}(\mathbb{R}_+^n)$, be defined by $\mathcal{P}(f,g) := u$, where $u$ is the unique solution of the Dirichlet problem
\[ (1-L)u = f,~\textrm{in } \mathbb{R}_+^n, \qquad u(x^\prime,0) = g(x^\prime),~\forall x^\prime \in \partial \mathbb{R}_n^+.\]
By standard Schauder theory this operator is well-defined, bounded and linear in each component.  We may therefore define a Dirichlet-to-Neumann operator, $\mathcal{N}$, by
\[ \mathcal{N}(f,g) := \left(\partial_{x_n} \mathcal{P}(f,g)\right)\bigg|_{\partial \mathbb{R}_+^n}.\] 
\noindent
\emph{Claim.} There exists a positive constant $C = C(K,n,\kappa_L)$ such that, for any $(f,g) \in C^{\alpha}(\overline{\mathbb{R}_+^n}) \times C^{1+\alpha}(\partial \mathbb{R}_+^n)$, 
\[ \| \mathcal{N}(f,g) \|_{C^\alpha(\overline{\mathbb{R}_+^n)}} \leq C\left( \|f\|_{C^\alpha(\overline{\mathbb{R}_+^n})} + \|g\|_{C^{1+\alpha}(\partial\mathbb{R}_+^n)}\right).\]
\noindent
\begin{proofof}{Claim.} Let $g \in C^{2+\alpha}(\partial\mathbb{R}_+^n) = C^{2+\alpha}(\mathbb{R}^{n-1})$ and $f \in C^{\alpha}(\overline{\mathbb{R}_+^n})$ be given.    By linearity of $L$, $\mathcal{N}(f,g) = \mathcal{N}(f,0)+ \mathcal{N}(0,g)$.  The first term, however, can be easily controlled since $\mathcal{P}(f,0)$ satisfies the classical Schauder estimate for vanishing boundary data:
\be \|\mathcal{P}(f,0)\|_{C^{2+\alpha}(\overline{\mathbb{R}_+^n})} \leq C\|f\|_{C^{\alpha}(\overline{\mathbb{R}_n^+})}, \label{estimateineq1} \ee 
where $C = C(\|a_{ij}\|_{C^\alpha}, \kappa_L, n) > 0$.  Therefore,
\begin{align*}
\|\mathcal{N}(f,g)\|_{C^{\alpha}(\partial\mathbb{R}_+^n)}  & \leq  \|\mathcal{N}(0,g)\|_{C^{\alpha}(\partial\mathbb{R}_+^n)} + \|\mathcal{N}(f,0)\|_{C^{\alpha}(\partial\mathbb{R}_+^n)} \\
 & \leq  \|\mathcal{N}(0,g)\|_{C^{\alpha}(\partial\mathbb{R}_+^n)} + \|\mathcal{P}(f,0)\|_{C^{1+\alpha}(\overline{\mathbb{R}_+^n})} \\ 
 & \leq  C \left(\|\mathcal{N}(0,g)\|_{C^{\alpha}(\partial\mathbb{R}_+^n)} + \|f\|_{C^{\alpha}(\overline{\mathbb{R}_+^n})} \right).\end{align*} 

Our main task is thus to estimate $\mathcal{N}(0,g)$.  Denote by $\Delta^\prime$ the Laplacian operator on $\mathbb{R}^{n-1}$.  Suppose at first that the coefficients of $L$ are constant.  The advantage of this is that it implies $1-L$ and $1-\Delta^\prime$ commute, from which we obtain the following identity
\[ \mathcal{P}\left(0,\left(1-\Delta^\prime\right)g\right) = \left(1-\Delta^\prime \right)\mathcal{P}(0,g). \]
Hence,
\begin{align} \left\|\mathcal{P}\left(0,\left(1-\Delta^\prime\right)g\right)\right\|_{C^{1+\alpha}(\overline{\mathbb{R}_+^n})} &\leq \|1-\Delta^\prime\|_{\mathcal{L}(C^{3+\alpha}(\overline{\mathbb{R}_+^n}),C^{1+\alpha}(\overline{\mathbb{R}_+^n}))}\|\mathcal{P}(0,g)\|_{C^{3+\alpha}(\overline{\mathbb{R}_+^n)}} \nonumber \\ 
& \leq  C\|g\|_{C^{3+\alpha}(\partial\mathbb{R}_+^n)}  \nonumber \\
& \leq  C\|(1-\Delta^\prime)g\|_{C^{1+\alpha}(\partial\mathbb{R}_+^n)},   \label{estimateineq2}  \end{align}
where the constant $C = C(\|a_{ij}\|_{C^\alpha},\kappa_L, n) > 0$.  

Standard elliptic theory tells us that $1-\Delta^\prime $ is a bounded linear operator mapping the space $C^{2+\alpha}(\mathbb{R}^{n-1})$ surjectively to $C^{\alpha}(\mathbb{R}^{n-1})$.  Then, replacing $g$ by $(1-\Delta^\prime)^{-1}g$ in \eqref{estimateineq2}, we get
\be \| \mathcal{P}(0,g) \|_{C^{1+\alpha}(\overline{\mathbb{R}_+^n})} \leq C \|g\|_{C^{1+\alpha}(\partial\mathbb{R}_+^n)}. \label{estimateineq3} \ee
The claim is an immediate consequence of this inequality, since
\[ \|\mathcal{N}(0,g)\|_{C^{\alpha}(\partial\mathbb{R}_+^n)} \leq  \|\mathcal{P}(0,g)\|_{C^{1+\alpha}(\overline{\mathbb{R}_+^n})}\leq C\|g\|_{C^{1+\alpha}(\partial \mathbb{R}_+^n)}.   \]
Combing this with \eqref{estimateineq1} we have the desired result in the constant coefficient case.

In order to complete the proof of the claim, therefore, we need only to reprove \eqref{estimateineq3} for variable coefficients.   Naturally, the idea will be to estimate the commutator $[1-L, 1-\Delta^\prime]$ effectively.  

Put $u := \mathcal{P}(0,(1-\Delta^\prime)g)$, $v := \mathcal{P}(0,g)$.  That is, unraveling notation, take $u$ and $v$ to satisfy:
\be (1-L)u = 0, \textrm{ in } \mathbb{R}_+^n, \qquad u = (1-\Delta^\prime)g, \textrm{ on } \partial\mathbb{R}_+^n, \label{ueqn} \ee
\be (1-L)v = 0, \textrm{ in } \mathbb{R}_+^n, \qquad v = g, \textrm{ on } \partial\mathbb{R}_+^n. \label{veqn} \ee
Put $w := (1-\Delta^\prime)v$.  Then, from \eqref{veqn}, we have that $w$ satisfies
\be (1-L)w = [1-L, 1-\Delta^\prime] v, \textrm{ in } \mathbb{R}_+^n, \qquad w = (1-\Delta^\prime)g, \textrm{ on } \partial\mathbb{R}_+^n.\label{weqn} \ee
Here we have let $[\cdot, \cdot]$ denote the commutator.  Subtracting \eqref{ueqn} from \eqref{weqn} leads to 
\[ (1-L)(w-u) = [1-L, 1-\Delta^\prime] v, \textrm{ in } \mathbb{R}_+^n, \qquad w-u = 0, \textrm{ on } \partial \mathbb{R}_+^n.  \]
Applying the Schauder estimates to this last equation we find 
\be \|w-u\|_{C^{2+\alpha}(\overline{\mathbb{R}_+^n})} \leq C\| [1-L, 1-\Delta^\prime] v \|_{C^{\alpha}(\overline{\mathbb{R}_+^n})}.\label{wminusuestimate1} \ee
Here $C = C(\|a_{ij}\|_{C^\alpha},\kappa_L,n) > 0$.  But note that $[1-L, 1-\Delta^\prime]$ is a third-order operator whose coefficients involves up to second-order derivatives of $a_{ij}$.  We therefore obtain from estimate \eqref{wminusuestimate1} that
\be  \|u - (1-\Delta^\prime) v\|_{C^{2+\alpha}(\overline{\mathbb{R}_+^n})} \leq C \|v\|_{C^{3+\alpha}(\overline{\mathbb{R}_+^n})} \leq C\|g\|_{C^{3+\alpha}(\partial\mathbb{R}_+^n)} \leq C \|(1-\Delta^\prime)g \|_{C^{1+\alpha}(\partial\mathbb{R}_+^n)}, \label{wminusuestimate2} \ee
where the constant $C = C(K,\kappa_L,n) > 0$. 

Now, rephrasing this in the earlier notation, we find
\begin{align*}  \left \| \mathcal{P}(0,(1-\Delta^\prime)g) \right\|_{C^{1+\alpha}(\overline{\mathbb{R}_+^n})} & \leq  \left \| (1-\Delta^\prime)\mathcal{P}(0,g) \right\|_{C^{1+\alpha}(\overline{\mathbb{R}_+^n})} \\
& \qquad +\left \| \mathcal{P}(0,(1-\Delta^\prime)g)-(1-\Delta^\prime)\mathcal{P}(0,g) \right\|_{C^{1+\alpha}(\overline{\mathbb{R}_+^n})}.\end{align*}
The first term on the right-hand side can be estimated as in \eqref{estimateineq2}, while the second term is bounded from above by \eqref{wminusuestimate2}.  Equivalently, if we replace $g$ by $(1-\Delta^\prime)^{-1} g$, it follows that
 \[ \left\|\mathcal{N}(0,g) \right\|_{C^{\alpha}(\partial\mathbb{R}_+^n)} \leq \left\| \mathcal{P}(0,g) \right\|_{C^{1+\alpha}(\overline{\mathbb{R}_+^n})} \leq C \|g\|_{C^{1+\alpha}(\partial\mathbb{R}_+^n},\]
 where $C = C(K, \kappa_L, n) > 0$.  This proves the claim.  \end{proofof}

Let $u \in C^{2+\alpha}(\overline{\mathbb{R}_+^n})$ be given.  For notational convenience, we let $h$ denotes the trace of $u$ on $\partial\mathbb{R}_+^n$.  The standard Schauder estimates for the Dirichlet problem are
\be \|u\|_{C^{2+\alpha}(\overline{\mathbb{R}_+^n})} \leq C \left( \|u\|_{C^\alpha(\overline{\mathbb{R}_+^n})} + \|h\|_{C^{2+\alpha}(\partial\mathbb{R}_+^n)} + \|Lu\|_{C^\alpha(\overline{\mathbb{R}_+^n})} \right),\label{standardschauderuest} \ee
where the constant $C = C(\|a_{ij}\|_{C^\alpha}, \kappa_L, n)$.  

Let $(f,g) \in C^{2+\alpha}(\overline{\mathbb{R}_+^n}) \times C^\alpha(\partial\mathbb{R}_+^n)$ be given, and suppose
\[ (1-L)u = f, \textrm{ in }\mathbb{R}_+^n, \qquad (M+N)u = g, \textrm{ on } \partial\mathbb{R}_+^n.\]
 We can express the second equation above entirely as a problem on $\partial\mathbb{R}_+^n$ by means of the operator $\mathcal{N}(f, \cdot)$:
\be \left(\sum_{i,j=1}^{n-1} \mathfrak{a}_{ij} \partial_{x_i}\partial_{x_j}  + \sum_{i=1}^{n-1} \mathfrak{b}_i \partial_{x_i} \right) h = -\mathfrak{b}_n \mathcal{N}(f,h) + g, \qquad \textrm{on } \partial \mathbb{R}_+^n.\label{Lambdaequation} \ee
Suppose that $n \geq 3$.  Applying the Schauder estimates on $\mathbb{R}^{n-1}$ to this equation yields
\be \|h\|_{C^{2+\alpha}(\partial\mathbb{R}_+^n)} \leq C\left(\|h\|_{C^\alpha(\partial\mathbb{R}_+^n)} + \|\mathcal{N}(f,h)\|_{C^\alpha(\partial\mathbb{R}_+^n)} + \|g\|_{C^\alpha(\partial\mathbb{R}_+^n)} \right),\label{hstandardschauderest} \ee
where $C = C(K, \kappa_M, \beta,  n) > 0$.  Applying the claim to the $\mathcal{N}$-term above and interpolating the resulting $C^{1+\alpha}$-norm between $C^\alpha$ and $C^{2+\alpha}$ we obtain
 \begin{align} \|h\|_{C^{2+\alpha}(\partial\mathbb{R}_+^n)} & \leq  C\left(\|h\|_{C^\alpha(\partial\mathbb{R}_+^n)} + \|f\|_{C^\alpha(\overline{\mathbb{R}_+^n})} + \|h\|_{C^{1+\alpha}(\partial\mathbb{R}_+^n)}+ \|g\|_{C^\alpha(\partial\mathbb{R}_+^n)} \right) \nonumber \\
 & \leq  \frac{1}{2}\|h\|_{C^{2+\alpha}(\partial\mathbb{R}_+^n)} + C\left(\|h\|_{C^\alpha(\partial\mathbb{R}_+^n)} + \|f\|_{C^\alpha(\overline{\mathbb{R}_+^n})} + \|g\|_{C^\alpha(\partial\mathbb{R}_+^n)} \right) \nonumber  \\
 & \leq  C\left(\|h\|_{C^\alpha(\partial\mathbb{R}_+^n)} + \|f\|_{C^\alpha(\overline{\mathbb{R}_+^n})} + \|g\|_{C^\alpha(\partial\mathbb{R}_+^n)} \right). \label{hestimate} \end{align}
 Here $C = C(K, \kappa_L, \kappa_M, \beta, n)$.  
 
 Inserting \eqref{hestimate} into \eqref{standardschauderuest} --- keeping in mind that $u = h$ on $\partial\mathbb{R}_+^n$ and $\|f\|_{C^{\alpha}(\overline{\mathbb{R}_+^n})} \leq \|u\|_{C^{\alpha}(\overline{\mathbb{R}_+^n})} + \|Lu\|_{C^{\alpha}(\overline{\mathbb{R}_+^n})}$ --- we arrive finally at the following:
\be\|u\|_{C^{2+\alpha}(\mathbb{R}^n_+)} \leq C\left(\|u\|_{C^\alpha(\overline{\mathbb{R}_+^n})} + \|Lu\|_{C^\alpha(\overline{\mathbb{R}_+^n})} +\|(M+N) u\|_{C^\alpha(\partial\mathbb{R}_+^n)}\right), \label{ufinalestimate} \ee
where $C = C(K, \kappa_L, \kappa_M, \beta, n)$.  This completes the lemma for $n > 2$.

For the case $n = 2$ we introduce an artificial variable, $t$, to move up a dimension.   Notice that the only place where the previous argument fails is in deriving \eqref{hstandardschauderest}, as Schauder theory would require the boundary to be at least two-dimensional.   Let $h, f, g$ and $\mathcal{N}(f,h)$ be as above and let their extensions $(\bar{h}, \bar{g}, \overline{\mathcal{N}(f,h)}) \in C^{2+\alpha}(\mathbb{R}^2) \times C^{\alpha}(\mathbb{R}^2) \times C^{1+\alpha}(\mathbb{R}^2)$ be defined as follows
\[ \bar{h}(t, x) := h(x), \qquad \bar{g}(t,x) := g(x), \qquad \overline{\mathcal{N}(f,h)}(t,x) := \mathcal{N}(f,h)(x), \qquad \forall (t,x) \in \mathbb{R}^2.\]
Of course, for any $k \geq 0$, the $C^{k+\alpha}(\mathbb{R}^2)$-norm of the extensions ($\bar{h}$, $\bar{g}$, $\overline{\mathcal{N}}$) coincides with the $C^{k+\alpha}(\mathbb{R})$-norm of the original original functions ($h$, $g$, $\mathcal{N}$).  

Define operators $\overline{M}$, $\overline{N}$  by
\[ \overline{M} := M + \partial_t^2, \qquad \overline{N} = N + \partial_t.\]
Then $(M+N)g = h$ iff $(\overline{M}+\overline{N})\bar{g} = \bar{h}$, so from \eqref{Lambdaequation} we obtain
\[ \left( \mathfrak{a}_{11} \partial_{x}^2  + \partial_t^2  +  \mathfrak{b}_1 \partial_{x} + \partial_t \right) \bar{h} = -\mathfrak{b}_2 \overline{\mathcal{N}(f,h)} + \bar{g}, \qquad \textrm{on } \mathbb{R}^2.\]
As we are now in the proper dimension, the usual elliptic estimates yield \eqref{hstandardschauderest} with $\bar{h}$, $\bar{g}$ and $\overline{\mathcal{N}}$ in place of $h$, $g$ and $\mathcal{N}$, respectively.  Then, using the fact that the extensions have the same relevant norms as the original functions, we may argue as before to get \eqref{hestimate} and thereby reestablish \eqref{ufinalestimate} for $n =2$.   \end{proof2}

\begin{longremark}  A standard cut-off function argument generalizes \textsc{Lemma \ref{estimateslemma}} to the case at hand, i.e. where we have $\Omega = R$, is bounded with periodic boundary conditions on the lateral sides and Dirichlet data on the top.
\label{estimatesremark}
\end{longremark}

Since we now have adequate \emph{a priori} estimates for the frozen problem, we fix any $\tau \in \mathbb{R}$ and define $\mathcal{G}^{(\tau)} : \mathbb{R} \times X \to Y$ by
\begin{align}
\mathcal{G}_1^{(\tau)}(h) & :=  (1+ h_q^2)h_{pp} + h_{qq} h_p^2 - 2h_q h_p h_{pq} - g(h-\tau) h_p^3 \rho_p + h_p^3 \beta(-p) \label{defG1tau} \\
 \mathcal{G}_2^{(\tau)}(Q,h) & :=  \mathcal{G}_2(Q,h) = \bigg(1+h_q^2 +h_p^2(2\sigma \kappa[h]+2g\rho h - Q)\bigg)\bigg|_{p = 0}. \label{defG2tau} \end{align}

We will show that, for $h \in \mathcal{O}_\delta$, $\mathcal{G}_1^{(\tau)}$ is a uniformly elliptic, quasilinear differential operator, while $\mathcal{G}_2^{(\tau)}$ is a non-degenerate, uniformly oblique Venttsel-type boundary operator.  For later reference we now compute the  Fr\'echet derivatives of $\mathcal{G}$ and $\mathcal{G}^{(\tau)}$:  
\be \begin{split}
 \mathcal{G}_{1h}(h) & =  2 h_q h_{pp} \partial_q + (1+ h_q^2) \partial_p^2 + 2h_p h_{qq} \partial_p + h_p^2 \partial_q^2 \\ 
 & \qquad - 2(h_{pq}h_q \partial_p + h_{pq}h_p \partial_q + h_p h_q \partial_p \partial_q)  \\
 & \qquad\qquad - 3g\rho_p (h-d(h)) h_p^2 \partial_p -g\rho_p h_p^3(1-d) + 3h_p^2 \beta(-p) \partial_p 
 \end{split}\label{capgravG1h} \ee
\be \mathcal{G}_{2h} (Q,h) =  \bigg(2h_q \partial_q +  2h_p(2\sigma\kappa[h] +2g\rho h- Q)\partial_p + 2g\rho h_p^2 + 2\sigma h_p^2 \kappa^\prime[h] \bigg)\bigg|_T  \label{capgravG2h} \ee
\be \begin{split}
 \mathcal{G}_{1h}^{(\tau)} (h) & =  \mathcal{G}_{1h}(h) + 3g(d(h)-\sigma) h_p^3 \rho_p \partial_p + g\rho_p h_p^3 d \\
 & =  2 h_q h_{pp} \partial_q + (1+h_q^2) \partial_p^2 + 2h_p h_{qq} \partial_p + h_p^2 \partial_q^2 \\
 & \qquad - 2(h_{pq}h_q \partial_p + h_{pq}h_p \partial_q + h_p h_q \partial_p \partial_q) \\
 & \qquad\qquad - 3g\rho_p (h-\tau) h_p^2 \partial_p -g\rho_p h_p^3 + 3h_p^2 \beta(-p) \partial_p \end{split}\label{capgravG1htau} \ee
\be \mathcal{G}_{2h}^{(\tau)}(Q,h)  =  \mathcal{G}_{2h}(Q,h), \label{capgravG2htau} \ee
where 
\[ \kappa^\prime[h] :=3 h_{qq}\left(1+h_q^2\right)^{-5/2} \partial_q - \left(1+h_q^2\right)^{-3/2} \partial_q^2.\]

Following \cite{walsh2009stratified}, we now establish the key admissibility lemmas.

\begin{lemma} \emph{(Proper Map)} Suppose $K$ is a compact subset of $Y$ and $D$ is a closed, bounded set in $\overline{\mathcal{O}_\delta}$, then $\mathcal{G}^{-1}(K) \cap D$ is compact in $\mathbb{R} \times X$. \label{capgravpropermaplemma} \end{lemma}
\begin{proof2} Let $\{(f_j, g_j)\}$ be a convergent sequence in $Y = Y_1 \times Y_2$ and assume that for each $j \geq 1$, $(f_j,g_j) = \mathcal{G}(Q_j, h_j)$ for some $(Q_j,h_j) \in \overline{\mathcal{O}_{\delta}}$ with $\{h_j\}$ bounded in $C_{\textrm{per}}^{3+\alpha}(\overline{R})$, $Q_j$ bounded in $\mathbb{R}$.  We wish to show that there exists a subsequence of $\{(Q_j,h_j)\}$ convergent in $\mathbb{R}\times X$.  
  
Denote $\theta_j := \partial_q h_j$, for $j \geq 1$.  Then, differentiating the relation between $(Q_j, h_j)$ and $(f_j, g_j)$, we find by \eqref{capgravdefG1} that, for all $j \geq 1$, 
\begin{align*}
\partial_q f_j & =  \partial_q \mathcal{G}_1(h_j) \\
& =  (1+ (\partial_q h_j)^2)\partial_p^2 \theta_j + (\partial_q(1+(\partial_q h_j)^2))\partial_p^2 h_j + (\partial_p h_j)^2 \partial_q^2 \theta_j + (\partial_q^2 h_j)\partial_q(\partial_p h_j)^2 \\
&  \qquad- 2(\partial_q h_j)(\partial_p h_j) \partial_p\partial_q \theta_j -2(\partial_q^2 h_j) (\partial_p h_j)( \partial_p\partial_q h_j) - 2(\partial_q h_j) (\partial_p\partial_q h_j)^2 \\
& \qquad\qquad  -g(\partial_p h_j)^3 \rho_p \theta_j - g(h_j - d(h_j)) \partial_q (\partial_p h_j)^3 + \partial_q(\partial_p h_j)^3 \beta(-p), \qquad \textrm{in } R.\end{align*}
Grouping terms above we may rewrite this in the form
\be \begin{split}
 (1+ (\partial_q^2 h_j)^2) \partial_p^2 \theta_j + (\partial_p h_j)^2 \partial_q^2 \theta_j 
 - 2(\partial_q h_j) (\partial_p h_j) \partial_q\partial_p \theta_j  &=  \\ 
 \partial_q f_j + F(\partial_q h_j, \partial_p \partial_q h_j, \partial_q^2 h_j, \partial_p h_j, \partial_p^2 h_j, d(h)), & \qquad \textrm{in } R, \end{split} \label{propermapReq} \ee
where $F$ is the cubic polynomial dictated by the previous expression.  Notice that the coefficients of the higher order terms above satisfy:
\[ 4(1+h_q^2)h_p^2 -4 h_q^2 h_p^2 = 4h_p^2  >  \delta^2.\] 
Hence, the right-hand side can be viewed as a uniformly elliptic operator acting on $\theta$.

We may likewise differentiate the equation on $T$ to discover that, for each $j \geq 1$,
\begin{align*}
\partial_q g_j & =  \partial_q \mathcal{G}_2 (Q_j, h_j) \\
& =  2(\partial_q h_j)\partial_q \theta_j + 2(\partial_p h_j)(2\sigma \kappa[h_j] +2g\rho h_j -Q_j) \partial_p \theta_j \\
& \qquad + 2g\rho (\partial_p h_j )^2 h_q + 2\sigma(\partial_p h_j)^2 (\partial_q \kappa[h_j]), \qquad \textrm{on } T. \end{align*}
Equivalently, 
\be \begin{split} 2\sigma(\partial_p h_j)^2 \partial_q \left(\frac{\partial_q \theta_j}{(1+h_q^2)^{3/2}}\right)+ 2(\partial_q h_j) \partial_q \theta_j + 2(\partial_p h_j)(2\sigma \kappa[h_j] +2g\rho h_j - Q_j) \partial_p \theta_j &=  \\
 \partial_q g_j + G(\partial_q h_j, \partial_p h_j), & \qquad \textrm{on } T, \end{split} \label{propermapTeq} \ee
where $G$ is the quadratic polynomial determined by the previous equation.  Observe that the boundary operator present on the left-hand side above consists of a one-dimensional elliptic operator added to a uniformly oblique operator, since
\[ |(2\sigma\kappa[h]+2g\rho h -Q)h_p | \geq \delta^2 \qquad \textrm{on } T\]
for $(Q,h) \in \overline{\mathcal{O}_{\delta}}$.
Finally, we note that 
\be \theta_j = 0 \qquad \textrm{on } B. \label{propermapBeq} \ee  

Since $\{(f_j, g_j)\}$ is convergent in $C_{\textrm{per}}^{1+\alpha}(\overline{R}) \times C_{\textrm{per}}^{1+\alpha}(T)$ we have that the sequences $\{\partial_q f_j\}$,  $\{\partial_q g_j\}$ are Cauchy in $C_{\textrm{per}}^{\alpha}(\overline{R})$ and $C_{\textrm{per}}^{\alpha}(T)$ respectively.  Now, by assumption, $\{h_j\}$ is uniformly bounded in $C_{\textrm{per}}^{3+\alpha}(\overline{R})$.   From this it follows that $F(\partial_q h_j, \partial_p \partial_q h_j, \partial_q^2 h_j, \partial_p h_j, \partial_p^2 h_j, d(h))$ is uniformly bounded in $C_{\textrm{per}}^{1+\alpha}(\overline{R})$, and $G(\partial_q h_j, \partial_p h_j)$ is uniformly bounded in $C_{\textrm{per}}^{2+\alpha}(T)$.  Each of these, viewed as a sequence in $j$, is therefore Cauchy. Then, by the compactness of the embeddings, the right-hand sides of equations \eqref{propermapReq} and \eqref{propermapTeq}  are strongly pre-compact in $C_{\textrm{per}}^\alpha(\overline{R})$ and $C_{\textrm{per}}^{\alpha}(T)$, respectively.  Possibly passing to a subsequence, we may take both to be convergent in these spaces.  

We now consider differences $\theta_j - \theta_k$,  for $j,k \geq 1$.  By \eqref{propermapReq} we have
\begin{align*}
F_{jk} & =  (1+(\partial_q^2 h_j)^2) \partial_p^2 (\theta_j-\theta_k) + (\partial_p h_j)^2 \partial_q^2 (\theta_j-\theta_k) \\
& \qquad - 2(\partial_q h_j) (\partial_p h_j) \partial_q\partial_p (\theta_j-\theta_k), \qquad \textrm{in } R, \end{align*}
where by our arguments in the previous paragraph we know $F_{jk} \to 0$ in $C^{\alpha}(\overline{R})$.  Similarly, from \eqref{propermapBeq} have that $\theta_j - \theta_k$ vanishes on the bottom and on the top \eqref{propermapTeq} tells us
 \begin{align*}
 G_{jk} & =  2\sigma(\partial_p h_j)^2 \partial_q \left(\frac{\partial_q\left(\theta_j-\theta_k\right)}{(1+h_q^2)^{3/2}}\right) + 2(\partial_q h_j) \partial_q (\theta_j-\theta_k) \\
 & \qquad + 2(\partial_p h_j)(2\sigma\kappa[h_j] - Q_j) \partial_p (\theta_j-\theta_k), \qquad \textrm{on } T.\end{align*}
Here $G_{jk} \to 0$ in $C^{\alpha}(T)$, again by the considerations of the preceding paragraph.  We now apply the mixed-boundary condition Schauder estimates of \textsc{Lemma \ref{estimateslemma}} to the differences $\theta_j - \theta_k$ to deduce that there exists some constant $C > 0$ with 
\[ C\|\theta_j-\theta_k\|_{C^{2+\alpha}(\overline{R})} \leq \|F_{jk}\|_{C^\alpha(\overline{R})} + \|G_{jk}\|_{C^{\alpha}(T)} + \|\theta_j - \theta_k\|_{C^\alpha(\overline{R})}.\]
Thus, $\theta_j - \theta_k \to 0$ in $C^{2+\alpha}$ as $j,k \to \infty$ .  We have shown, therefore, that all third derivatives of $\{h_j\}$ are Cauchy in $C^\alpha_{\textrm{per}}(\overline{R})$, except possibly for $\{\partial_p^3 h_j\}$.  To demonstrate the same holds for $\{\partial_p^3 h_j\}$, we use the PDE to express $\partial_p^2 h_j$ in terms of the derivatives of order less than or equal to two: 
 \begin{align*}
 \partial_p^2 h_j &= (1+ (\partial_q h_j)^2)^{-1} \Big( f_j -  \partial_q^2 h_j (\partial_p h)^2 + 2(\partial_q h)(\partial_p h_j)(\partial_p \partial_q h_j) \\
 & \qquad -g(h_j - d(h_j)) \rho_p (\partial_p h_j)^3 + (\partial_p h_j)^3 \beta(-p)\Big), \qquad \textrm{in } R.\end{align*}
 But we have seen that the right-hand side is Cauchy in $C^{1+\alpha}(\overline{R})$, hence $\{\partial_p^3 h_j\}$ is also Cauchy in $C_{\textrm{per}}^\alpha(\overline{R})$.  We conclude that the original sequence, $\{h_j\}$, is Cauchy in $C^{3+\alpha}$.   \end{proof2}
 
\begin{lemma} \emph{(Fredholm Map)} For each $(Q,h) \in \mathcal{O}_{\delta}$, the linearized operator $\mathcal{G}_h(Q,h)$ is a Fredholm map of index 0 from $X$ to $Y$. \label{capgravfredholmlemma} \end{lemma}

\begin{proof2} Let $\psi \in C_{\textrm{per}}^{3+\alpha}(\mathbb{R})$ be given and fix $\tau,~Q \in \mathbb{R}$, $h \in X$.  Put \[ \phi^{(i)} := \partial_q\left(\mathcal{G}_{ih}^{(\tau)}(Q,h)[\psi]\right)- \left(\partial_q \mathcal{G}^{(\tau)}_{ih}(Q,h)\right)[\psi], \qquad \textrm{for } i=1,2.\]  
Here, by $(\partial_q \mathcal{G}_{ih}^{(\tau)}(Q,h))[\psi]$ we mean differentiating the coefficients of $\mathcal{G}_{ih}^{(\tau)}$ in $q$,  then applying the resulting operator to $\psi$.  Then $\partial_q \psi$ satisfies:
\[ \left\{ \begin{array}{lll}
 \mathcal{G}_{1h}^{(\tau)}(h)  \partial_q\psi = \phi^{(1)} & & \textrm{on } R, \\
 & & \\
 \mathcal{G}_{2h}^{(\tau)}(Q,h)  \partial_q \psi =  \phi^{(2)} & & \textrm{on } T, \\ 
 & & \\
 \partial_q \psi = 0 & & \textrm{on } B, \end{array} \right.\]
 which is a uniformly elliptic PDE with a Venttsel boundary condition.  Thus there exists a constant $C > 0$, independent of $\psi$, such that
 \begin{align*}
 C \|\partial_q \psi\|_{C^{2+\alpha}(\overline{R})} & \leq   \|\partial_q \psi \|_{C^\alpha(\overline{R})} + \| \phi^{(1)} \|_{C^\alpha(\overline{R})} + \| \phi^{(2)} \|_{C^{\alpha}(T)}  \\
 & \leq  \|\psi\|_{C^{2+\alpha}(\overline{R})} + \|\partial_q (\mathcal{G}_{1h}^{(\tau)}(h) \psi) \|_{C^\alpha(\overline{R})} + \|\partial_q (\mathcal{G}_{2h}^{(\tau)}(Q,h) \psi) \|_{C^{\alpha}(T)}.  \end{align*}
On the other hand, we may express $\partial_p^2 \psi$ via the partial differential equation to arrive at an estimate for $\partial_p^3 \psi$ of the same type.  Combining these estimates we find that, for some $C > 0$ and all $\psi \in X$, 
\be C \|\psi\|_{C^{3+\alpha}(\overline{R})} \leq \|\psi \|_{C^{2+\alpha}(\overline{R})} + \|\partial_q \phi^{(1)} \|_{C^\alpha(\overline{R})} + \|\partial_q \phi^{(2)} \|_{C^\alpha(T)}. \label{capgravestimate1fredholmmap} \ee 
If we now apply \textsc{Lemma \ref{estimateslemma}} to $\psi$ directly, we can estimate the $C^{2+\alpha}$--norm of $\psi$.  Then \eqref{capgravestimate1fredholmmap} becomes
\be C \|\psi\|_X \leq \|\psi\|_{Y_1} + \|\mathcal{G}_{1h}^{(\tau)}(h) \psi\|_{Y_1} + \|\mathcal{G}_{2h}^{(\tau)}(Q,h) \psi\|_{Y_2}. \label{capgravestimate2fredholmmap} \ee

This is the key estimate, from which we can conclude  the range of $\mathcal{G}_{ih}^{(\tau)}$ is closed and the null space finite dimensional.  The task at hand, therefore, is to prove it with $\mathcal{G}$ in place of $\mathcal{G}^{(\tau)}$.  As in \cite{walsh2009stratified}, this is made possible by the following trivial but convenient fact: for any $\psi$ as above we have
\[ |d(\psi)| \leq \fint_T |\psi|dq \leq \|\psi\|_{C^\alpha(\overline{R})}. \]
We are therefore able to estimate terms involving $d$ easily in spaces of H\"older continuous functions.  In particular,
\begin{align}
\|\mathcal{G}_{1h}(h)\psi - \mathcal{G}_{1h}^{(\tau)}(h)\psi \|_{C^\alpha(\overline{R})} & =  \|3g(d(h)-\tau) h_p^3 \rho_p \psi_p + g\rho_p h_p^3 d(\psi)\|_{C^\alpha(\overline{R})} \nonumber \\
& \leq  C \|\psi\|_{C^{1+\alpha}(\overline{R})}, \label{capgravGGsigmaCalphaerror} \end{align}
where the constant $C$ above depends only on $\rho_p, \|h\|_X$ and our choice of $\sigma$.  Likewise, an identical argument gives the estimates
\[ \|\partial_i \mathcal{G}_{1h}(h)\psi - \partial_i \mathcal{G}_{1h}^{(\tau)}(h)\psi\|_{C^\alpha(\overline{R})} \leq C\|\psi\|_{C^{2+\alpha}(\overline{R})}, \qquad i = p,~q \]
where $C$ depends again on $\rho_p, \|h\|_X$ and the choice of $\tau$.  Of course we do not need to make such arguments to estimate the $\mathcal{G}_{2h}^{(\tau)}(h)$ term, as it identical to $\mathcal{G}_{2h}(h)$.

Combining these observations with \eqref{capgravestimate2fredholmmap}, we find that, for some $C > 0$ and all $\psi \in X$, 
\be C \|\psi\|_X \leq \|\psi\|_{C^{2+\alpha}(\overline{R})} + \|\mathcal{G}_{1h}(h) \psi\|_{Y_1} + \|\mathcal{G}_{2h}(Q,h) \psi\|_{Y_2}. \label{estimate2fredholmmap} \ee
Applying the arguments of \cite{walsh2009stratified}, we are able to conclude that $\mathcal{G}_h(Q,h)$ is Fredholm.  The Fredholm index has a discrete range, hence by the connectedness of $\mathcal{O}_{\delta}$ it must be constant on this set.  But $\mathcal{G}_h(Q_*, H_*)$ was shown to have a one-dimensional null space and range of codimension one.  Since $(Q_*, H_*) \in \mathcal{O}_{\delta}$ by construction, the Fredholm index is uniformly 0 along the continuum $\mathcal{C}_{\delta}^\prime$. \end{proof2}

\begin{lemma} \emph{(Spectral Properties)} \emph{(i)} $\forall \delta, \epsilon > 0$, $\exists c_1, c_2 > 0$ such that for all $(Q,h) \in \mathcal{O}_{\delta}$ with $|Q| + \|h\|_X \leq M$, for all $\psi \in X$ and for all $\mu \in \mathbb{C}$ with $|\mu| \geq c_2$ and $|\mathrm{arg~}\mu | < \pi/2 - \epsilon$, we have
\[ c_1 \|\psi\|_X \leq |\mu|^{\alpha/2}\|(A-\mu)\psi \|_{Y_1} + |\mu|^{\alpha/2} \|(B-\mu)\psi \|_{Y_2}\]
where $A = A(Q,h) = \mathcal{G}_{1h}(h)$ and $B = B(Q,h) = \mathcal{G}_{2h}(Q,h)$.  \\

\noindent
\emph{(ii)} Let  $\Sigma = \Sigma(Q,h)$ denote the spectrum of $(A,B)$.  Then $\Sigma$ consists entirely of the eigenvalues of finite multiplicity with no finite accumulation points.  Furthermore, there is a neighborhood $\mathcal{N}$ of $[0,+\infty)$ in the complex plane such that $\Sigma(\lambda, w) \cap \mathcal{N}$ is a finite set.   \label{capgravspectralpropertieslemma} \end{lemma} 

\begin{proof2}  To prove (i) we follow the argument of Agmon (cf. \cite{agmon1962eigenfunctions}).  The Venttsel boundary conditions, however, necessitate some modification of the standard argument.    All that we must do here is trivially modify the technical tricks of that argument.  Part (ii) follows immediately from part (i) in light of the previous lemmas, so we omit the details.  

Fix $\tau \in \mathbb{R}$ with $|\tau| < M$, and let $\psi \in X$, $\mu \in \mathbb{C}$ be given.  Put $\theta := \arg \mu$, and suppose for some $\epsilon > 0$, $|\theta| < \pi/2 + \epsilon$.  Consider the operators 
\[ A^{(\tau)} := \mathcal{G}_{1h}^{(\tau)}(h), ~D_1^{(\tau)} := A^{(\tau)} + e^{i\theta}\partial_t^2, \qquad \textrm{on } \mathbb{R}\times R,\]
and 
\[ D_2^{(\tau)} = D_2 := B + e^{i\theta}\partial_t^2, \qquad \textrm{on } \mathbb{R} \times T.\]
Then, for any $\tau$, $D_1^{(\tau)}$ is elliptic with constant of ellipticity independent of $\tau$, while $D_2^{(\tau)}$ is of the form of \textsc{Lemma \ref{estimateslemma}}.  Let $\zeta: \mathbb{R} \to \mathbb{R}$ be a cutoff function supported compactly in the interval $I := (-1,1)$.  Put 
\[ e(t) := e^{i|\mu|^{1/2} t} \zeta(t), \qquad \phi(t,q,p) := e(t) \psi(q,p).\]
We apply the Schauder estimates of \textsc{Lemma \ref{estimateslemma}} in $\mathbb{R}^3$ with the boundary operator $D_2^{(\tau)}$ to $\phi(t,q,p)$ to deduce
\be C\|\phi\|_{C^{2+\alpha}(I \times R)} \leq  \|\phi\|_{C^\alpha(I \times R)} + \|D_1^{(\tau)} \phi\|_{C^{\alpha}(I \times R)} + \|D_2^{(\tau)}\phi\|_{C^{\alpha}(I \times T)} \label{spectralbasicestimate} \ee
for some constant $C > 0$ independent of $\psi$.  In fact, since $\tau$ occurs as the coefficient of a first-order term in $\mathcal{G}_{1h}^{(\tau)}(Q,h)$, we have by \textsc{Lemma \ref{estimateslemma}} that $C$ can be chosen independently of $\tau$ (recall, $|\tau| < M$.)

A quick calculation readily confirms the existence of $C,~C^\prime > 0$, depending only on $\zeta$, with
\be C \mu^{\alpha/2} \leq \| e \|_{C^\alpha(\mathbb{R})} \leq C^\prime \mu^{\alpha/2}, \qquad C \mu^{(1+\alpha)/2} \leq  \| e \|_{C^{1+\alpha}(\mathbb{R})} \leq C^\prime \mu^{(1+\alpha)/2}. \label{spectralestimatee}\ee
Using this, we can unpack \eqref{spectralbasicestimate} to derive the following set of estimates: 
\begin{align}
 \|D_1^{(\tau)} \phi \|_{C^{\alpha}(I \times R)} & =  \left\| e(t) \mathcal{G}_{1h}^{(\tau)} (h) [\psi] - \mu e(t)\psi  + \left(i|\mu|^{1/2}\zeta^\prime + \zeta^{\prime\prime}\right) e^{i\theta} e^{i |\mu|^{1/2} t} \psi \right\|_{C^{\alpha}(I \times R)} \nonumber \\
 & \leq  C^\prime |\mu|^{\alpha/2} \left( \| \mathcal{G}_{1h}(h)- \mu) \psi \|_{C^{\alpha}(R)} + \|(\mathcal{G}_{1h}(h)- \mathcal{G}_{1h}^{(\tau)}(h)) \psi\|_{C^{\alpha}(R)}\right) \nonumber \\
 & \qquad  + \|  \left(i|\mu|^{1/2}\zeta^\prime + \zeta^{\prime\prime}\right) e^{i\theta} e^{i |\mu|^{1/2} t} \psi \|_{C^{\alpha}(I \times R)}  \nonumber \\
 & \leq  C^\prime |\mu|^{\alpha/2}  \left(\left\|\left(\mathcal{G}_{1h}(h) -\mu\right) \psi \right\|_{C^{\alpha}(R)}+\| \psi \|_{C^{\alpha}(R)} \right) \label{spectralAesimatetau} \\
 & \qquad + C^\prime \left(|\mu|^{(1+\alpha)/2} + |\mu|^{\alpha/2} \right) \| \psi \|_{C^{\alpha}(R)} \label{spectralAestimate} \end{align}
 for $|\mu|$ sufficiently large.  In transitioning from \eqref{spectralAesimatetau}  to \eqref{spectralAestimate}, we have made use of \eqref{capgravGGsigmaCalphaerror} taking $\tau = d(h)$  (which is valid since $d(h) < \|h\|_X < M$ and \eqref{spectralbasicestimate} holds for any $\tau$.)  A similar analysis of the boundary terms yields,
 \begin{align}
 \|D_2^{(\tau)} \phi\|_{C^{\alpha}(I \times T)} & =  \left\| e(t) \mathcal{G}_{2h}(Q,h)[\psi] - \mu e(t) \psi + \left(|\mu|^{1/2}\zeta^\prime + \zeta^{\prime\prime}\right) e^{i\theta} e^{i |\mu|^{1/2} t} \psi \right\|_{C^{\alpha}(I \times T)} \nonumber \\
 & \leq  C^\prime |\mu|^{\alpha/2} \| \left(\mathcal{G}_{2h}(Q,h) - \mu \right)\psi\|_{C^{\alpha}(T)} + C^\prime \left(|\mu|^{(1+\alpha)/2} + |\mu|^{\alpha/2} \right) \| \psi \|_{C^{\alpha}(T)}, \label{spectralBestimate} \end{align}
for $|\mu|$ is sufficiently large. 
 
 On the other hand,
 \begin{align}  [\phi ]_{2, \alpha; R \times I} & = [ \partial_t^2 e(t) \psi ]_{0,\alpha; R \times I} + [\partial_t e(t) \partial_q \psi]_{0,\alpha; R \times I} +  [\partial_t e(t) \partial_p \psi]_{0,\alpha; R \times I} \nonumber \\
 & \qquad +  [e(t) \partial_q^2 \psi]_{0,\alpha; R \times I} +  [e(t) \partial_q \partial_p \psi]_{0,\alpha; R \times I} +  [e(t) \partial_p^2 \psi]_{0,\alpha; R \times I} \nonumber \\
 & \geq C^{\prime \prime} |\mu| [\psi]_{0,\alpha; R} +  [\psi]_{2,\alpha; R}, \nonumber \\
 \| \phi \|_{C^2(I\times R)} & \geq C^{\prime\prime} |\mu| \| \psi \|_{C^0(\overline{R})} + \| \psi \|_{C^2(\overline{R})} \label{spectrallhsestimate},   \end{align}
 for $|\mu|$ sufficiently large, and for some  $C^{\prime\prime} > 0$ independent of $\mu$ and $\psi$.  Here, $[\cdot]_{k,\alpha; I \times R}$ denotes the usual H\"older semi-norm.  Combining \eqref{spectralbasicestimate}, \eqref{spectralAestimate}, \eqref{spectralBestimate}, and \eqref{spectrallhsestimate}, we find derive the following:  there exists $c_1, c_2 > 0$ such that, 
 \[ c_1 \| \psi \|_{C^{2+\alpha}(\overline{R})} \leq |\mu|^{\alpha/2} \| \left(A - \mu)\psi\right\|_{C^\alpha{(\overline{R}})} + |\mu|^{\alpha/2} \| \left(B - \mu)\psi\right\|_{C^\alpha{(\overline{R}})}, \]
 for all $\mu \in \mathbb{C}$ with $|\mu| > c_2$.  If we now repeat the same argument with $\psi_q$ and $\psi_p$ in place of $\psi$, we arrive at the inequality in (i). \end{proof2}

In order to complete the proof of \textsc{Theorem \ref{capgravglobalbifurcationtheorem}} we will make use of a generalization of Leray--Schauder degree due to Kiel\"ofer (cf. \cite{kielhofer1985multiple,kielhofer2004bifurcation}).  With that in mind, we recall the following definitions.

\begin{definition} (See \textsc{Definition II.5.2} of \cite{kielhofer2004bifurcation}.) Let $W$, $Z$ be real Banach spaces with $W$ continuously embedded in $Z$ and let $\mathcal{W}$ be an open, bounded subset of $W$.  We say that a map $G \in C^2(\overline{\mathcal{W}}; Y)$ is \emph{admissible} provided that it satisfies the following:
\begin{itemize}
\item[(i)] $G$ is a proper map.
\item[(ii)] For all $w \in \mathcal{W}$, the linearized operator $G_w(w)$ is Fredholm of index 0.  \\
\item[(iii)] There exists $\epsilon \in (0,\pi/2)$, $\alpha \in (0,1)$, and constants $c_1, c_2 > 0$ such that, for all $\psi \in X$, $w \in \mathcal{W}$, and $\mu \in \mathbb{C}$ with $|\mu| > c_2$, $|\mathrm{arg~}\mu| < \pi/2-\epsilon$, we have
\[ c_1 \|\psi\|_{W} \leq |\mu|^{\alpha/2} \| (G_w(w)-\mu) \psi\|_{Z},\]
where we have complexified $W$ and $Z$.\\
\item[(iv)] For each $w \in \mathcal{W}$, the spectrum $\Sigma$ of $G_w(w)$ consists entirely of eigenvalues of finite multiplicity with no finite accumulation points.  Furthermore, there is a neighborhood $\mathcal{N}$ of $[0,+\infty)$ in the complex plane such that $\Sigma(\lambda, w) \cap \mathcal{N}$ is a finite set. 
\end{itemize} \label{admissiblekielhofer}
\end{definition}  

The significance of admissible maps is that they permit us to define a degree.  Suppose that $G$ is such a map and let $z \in Z\setminus G(\partial \mathcal{W})$ be a regular value of $G$.  We define the \emph{(Kiel\"ofer) degree} of $G$ at $z$ with respect to $\mathcal{W}$ to be the quantity
 \[ \textrm{deg }({G}, \mathcal{W}, z) := \sum_{w \in \mathcal{G}^{-1}(\{z\})} (-1)^{\nu(w)},\]
 where $\nu(w)$ is the number of positive real eigenvalues counted by multiplicity of ${G}_{w}(w)$.  By part (iii) of \textsc{Definition \ref{admissiblekielhofer}}, we have that $\nu(w)$ is finite.  Likewise, part (i) ensures ${G}^{-1}(\{z\}) \cap \mathcal{W}$ is a finite set.  The degree is therefore well-defined at proper values.  In the usual way, this definition is extended to critical values via the Sard-Smale theorem (this is valid by since $G$ is Fredholm.)    
 
 Our interest in degree theory stems from the fact that the degree is invariant under any homotopy that respects the boundaries of $\mathcal{W}$ in the following sense.  Let $\mathcal{U} \subset [0,1] \times \mathcal{W}$ be open.  Define $\mathcal{U}_t := \{ w \in \mathcal{W} : (t,w) \in \mathcal{U}\}$ and likewise $\partial \mathcal{U}_t := \{ w \in \mathcal{W} : (t,w) \in \partial \mathcal{U}\}$.
 
 \begin{lemma} \emph{(Homotopy Invariance)} The degree is invariant under admissible homotopies. That is, suppose  $\mathcal{U} \subset [0,1] \times \mathcal{W}$ is open.  If $\mathcal{H} \in C^2(\overline{\mathcal{U}}; Y)$ is proper and, for each $t \in [0,1]$, $\mathcal{H}(t,\cdot)$ is admissible, then we say $\mathcal{H}$ is an admissible homotopy and have 
 \[ \mathrm{deg}~(\mathcal{H}(0,\cdot), \mathcal{U}_0, y) = \mathrm{deg}~(\mathcal{H}(1,\cdot), \mathcal{U}_1, y)\]
 provided $y \notin \mathcal{H}(t,\partial \mathcal{U}_t)$ for all $t \in [0,1]$. \label{homotopyinvariancelemma} \end{lemma}
 The reader is directed to \cite{kielhofer2004bifurcation} for a thorough treatment of degree theory and a proof of this result.  
 
 We are now prepared to prove the main theorem of this section.
 \begin{proofof}{Theorem \ref{capgravglobalbifurcationtheorem}}  From part (iii) of \textsc{Definition \ref{admissiblekielhofer}}, we see that in order to apply the Kielh\"ofer degree theory, we must first make a technical change.  For $h \in C_{\textrm{per}}^{3+\alpha}(\overline{R})$, put 
\[ \widetilde{h} := (h, h|_T) \in C_{\textrm{per}}^{3+\alpha}(\overline{R}) \times C_{\textrm{per}}^{3+\alpha}(T),\]
and 
\[ \widetilde{X} := \left\{ (h, h|_T) : h \in  C_{\textrm{per}}^{3+\alpha}(\overline{R}) \right\} \subset C_{\textrm{per}}^{3+\alpha}(\overline{R}) \times C_{\textrm{per}}^{3+\alpha}(T).\]
If we let $\pi$ be the bijective projection $\pi: X \to \widetilde{X}$, $h \mapsto_{\pi} (h, h|_T)$, then we can likewise define
\[ \widetilde{\mathcal{G}}(Q, \widetilde{h}) := \mathcal{G}(Q, \cdot) \circ \pi, \qquad \widetilde{\mathcal{T}} := \pi(\mathcal{T}), \qquad \widetilde{\mathcal{O}}_\delta := \pi(\mathcal{O}_\delta), \qquad \widetilde{S}_\delta := \pi(S_\delta), \qquad \widetilde{\mathcal{C}_\delta^\prime} = \pi( \mathcal{C}_\delta^\prime). \]
Let $\mathcal{W}$ be an open bounded subset of $X$ and fix $Q \in\mathbb{R}$.  By virtue of \textsc{Lemma \ref{capgravpropermaplemma}}--\textsc{Lemma \ref{capgravspectralpropertieslemma}}, we conclude that the map $\widetilde{\mathcal{G}}(Q,\cdot):\overline{\pi(\mathcal{W})} \to Y$ is admissible in the sense of \textsc{Definition \ref{admissiblekielhofer}}.  This is a trivial consequence of the fact that $\pi$ is a smooth bijection from $X$ onto $\widetilde{X}$.  
 
From our previous analysis it follows that $\widetilde{\mathcal{G}}(Q_1, \cdot) \mapsto \widetilde{\mathcal{G}}(Q_2, \cdot)$ is an admissible homotopy.  Thus, in light of \textsc{Lemma \ref{homotopyinvariancelemma}}, the degree will remain constant as we move along the continuum $\widetilde{\mathcal{C}_\delta^\prime}$.  Assuming that all three alternatives of \textsc{Theorem \ref{capgravglobalbifurcationtheorem}} fail, we use this feature to generate a contradiction.  At this stage, we have reestablished all the relevant properties of $\widetilde{\mathcal{G}}$ that were true in the constant density case.  Repeating (verbatim) the proof in \cite{constantin2004exact} --- with $\widetilde{\mathcal{G}}$ in place of $\mathcal{G}$ --- we obtain the theorem.   \end{proofof}

\begin{longremark} The motivation for choosing the Kielh\"ofer degree (rather than that of Healey and Simpson used, for example, in \cite{constantin2004exact}) lies in \textsc{Lemma \ref{capgravspectralpropertieslemma}}:  due to the presence of higher-order derivatives, the boundary operator and interior operator are equivalent, in the sense of the Agmon-type spectral estimates of \textsc{Lemma \ref{capgravspectralpropertieslemma}} (i).   This observation was first made by Wahl\'en in his analysis of the pure capillary wave case, where he also introduced the clever trick of considering $\widetilde{h}$, $\widetilde{X}$, etc. (cf. \cite{wahlenthesis}).  \label{agmonremark} \end{longremark}

\section{Analytic global bifurcation} \label{analyticitysection}

In the previous section we obtained global bifurcation results for the case when the eigenvalue at $(H_*, Q_*)$ is simple.  The basis for that method of argument, which harkens back to Rabinowitz \cite{rabinowitz1971some}, is the topological theory of degree.  In the present section, we present an alternative theory of bifurcation that instead relies on the (real) analytic properties of the nonlinear operator $\mathcal{G}$.  

This approach is originally due to Dancer (cf. \cite{dancer1973bifurcation,dancer1973globalstructure}, or \cite{buffoni2003analytic} for a survey) and has recently been applied with great success to the study of Stokes waves (cf. \cite{buffoni1998ondes,buffoni2000regularity,buffoni2000sub}).  The major insight of Dancer's theory is that the germ of the zero-set of an analytic operator has an extremely special structure.  Within the context of global continuation this means, for example, that if $\mathcal{C}$ is an analytic curve that corresponds to the zero-set of an analytic function, then there is a unique, maximal, analytic extension of $\mathcal{C}$.  Lifting this fact to infinite dimensions (and assuming some additional compactness properties) provides a very powerful set of tools for global bifurcation theory.  

With respect to the degree theoretic methods of the previous section, Dancer's analytic continuation method offers several advantages.  The first of these is that it provides a more natural way to treat the double bifurcation case.  Note that, as the multiplicity of the eigenvalue will be even, a Rabinowitz-type argument requires some additional work; in particular, special attention must be paid to the scenario where the continuations of the local branches intersect.  On the other hand, since maximal analytic extensions are unique, when working within the analytic framework, we are untroubled by such intersections.  Secondly, since we are using more than simply topological information, we obtain considerably better regularity for the global bifurcation branch.  Indeed, we will prove that the solution continuum is, in fact, path-connected and admits a locally injective, analytic reparameterization (in the sense of Definition \ref{defparameterization}.)   

Let us begin by presenting the central result of the analytic theory of global bifurcation.  This was originally proved by Dancer (cf. \cite{dancer1973bifurcation,dancer1973globalstructure}), but we shall use a slight rephrasing of that result due to Buffoni and Toland (cf. \cite{buffoni2003analytic}).  For the purposes of this exposition, we shall temporarily set aside our previous notation.  Let real Banach spaces $W,\,Z,$ open set $U \subset \mathbb{R} \times W$ and a real analytic $\mathcal{F} \colon U \to Z$ be given.  (Looking forward, we should think of $U$ as taking the place of $\mathcal{O}_\delta$, $W$ that of $X$, $Z$ that of $Y$ and $\mathcal{F}$ that of $\mathcal{G}$.)  If the following structural properties hold
\begin{indentenum}
\item[(G1)] $\mathcal{T} := \{(Q, 0) \colon Q \in \mathbb{R}\} \subset U \textrm{ and } \mathcal{F}(Q, 0) = 0, \forall Q \in \mathbb{R},$ \\

\item[(G2)] $\mathcal{F}_w(Q, w)$ is a Fredholm operator of index zero for each $(Q, w) \in U$ such that  $\mathcal{F}(Q, w) = 0$,  \\

\item[(G3)] $\exists Q_* \in \mathbb{R}$  such that {(i)} $\mathcal{N}(\mathcal{F}_w(Q_*, 0)) = \{ s \phi_* \colon s \in \mathbb{R} \}$,  and {(ii)} $\mathcal{F}_{Q w}(Q_*, 0)\phi_* \notin \mathcal{R}( \mathcal{F}_w(Q_*,0))$,
\end{indentenum}
then, by the analytic equivalent of \textsc{Theorem \ref{crandallrabinowitz}}, there exists an analytic local bifurcation curve at $(Q_*,0)$.  That is, there exists an $\epsilon > 0$ and an analytic $(\mathfrak{Q}, \mathfrak{w}) \colon (-\epsilon, \epsilon) \to \mathbb{R} \times W$ such that $\mathcal{F}(\mathfrak{Q}(s), \mathfrak{w}(s)) = 0$, for all $s \in (-\epsilon, \epsilon)$, $\mathfrak{Q}(0) = Q_*$ and $\mathfrak{w}^\prime(s) = \phi_*$.  Let us denote
\be \mathcal{C}_{\textrm{loc}}^{+} := \{ (\mathfrak{Q}(s), \mathfrak{w}(s)) \colon s \in (0, \epsilon)\}, \qquad
\mathcal{S} := \{ (Q, w) \in U \colon \mathcal{F}(Q, w) = 0\}.  \label{defCS} \ee
By taking $\epsilon$ sufficiently small, we may assume that $\mathfrak{w}^\prime \neq 0$ and $\mathcal{C}_{\mathrm{loc}}^+ \subset \mathcal{S}\setminus\mathcal{T}$.
 
\begin{theorem} \emph{(Analytic Global Simple Bifurcation)} Let $W,\,Z,\,U$ and $\mathcal{F}$ be given as above and assume that \emph{(G1), (G2)} and \emph{(G3)} hold.  We may then take $(\mathfrak{Q}, \mathfrak{w})$ as before and let  $\mathcal{C}_{\mathrm{loc}}^{+} $ and $\mathcal{S}$ be defined by \eqref{defCS}.  Suppose further that 
\be \mathfrak{Q}^\prime \nequiv 0 \textrm{ on } (-\epsilon, \epsilon), \tag{G4} \ee
\be \textrm{all bounded closed subsets of } \mathcal{S} \textrm{ are compact.} \tag{G5} \ee   
Then there exists a unique, maximal, continuous curve $\mathcal{K}$ that extends $\mathcal{C}_{\mathrm{loc}}^{+}$ in the following sense.
\begin{enumerate}
\item[\emph{(a)}] $\mathcal{K} = \{ ( \mathfrak{Q}(s), \mathfrak{w}(s)) \colon s \in [0,\infty)) \subset U$, where $(\mathfrak{Q}, \mathfrak{w}) \colon [0,\infty) \to \mathbb{R} \times W$ is continuous.  \\
\item[\emph{(b)}] $\mathcal{C}_{\mathrm{loc}}^{+} \subset \mathcal{K} \subset \mathcal{S}$. \\
\item[\emph{(c)}] The set $\{ s \geq 0 \colon \mathcal{N}( \mathcal{F}_w(\mathfrak{Q}(s), \mathfrak{w}(s))) \neq \{ 0 \} \}$ has no accumulation points. \\
\item[\emph{(d)}] At each point, $\mathcal{K}$ has a local analytic reparameterization. \\
\item[\emph{(e)}] One of the following occurs. 
\begin{enumerate}
\item[\emph{(i)}] $\| (\mathfrak{Q}(s), \mathfrak{w}(s) \| \to \infty$ as $s \to \infty$; or 
\item[\emph{(ii)}] $(\mathfrak{Q}(s), \mathfrak{w}(s))$ approaches the boundary of $U$ as $s \to \infty$; or   
\item[\emph{(iii)}] $\mathcal{K}$ is a closed loop, i.e. $\exists T > 0$ such that $(\mathfrak{Q}(s + T), \mathfrak{w}(s+T)) = (\mathfrak{Q}(s), \mathfrak{w}(s))$, for all $s \in [0,\infty)$.  This implies that $\mathcal{K} = \{ (\mathfrak{Q}(s), \mathfrak{w}(s)) \colon s \in [0,T] \}$. \\
\end{enumerate} 
\item[\emph{(f)}] If, for some $0 < s_1 < s_2$, we have 
\[(\mathfrak{Q}(s_1), \mathfrak{w}(s_1)) = (\mathfrak{Q}(s_2), \mathfrak{w}(s_2)), \textrm{ and }\mathcal{N}(\mathcal{F}_w(\mathfrak{Q}(s_1), \mathfrak{w}(s_2))) = \{ 0 \},\]
  then \emph{(e)(iii)} occurs and $s_2 - s_1$ is an integer multiple of $T$.  
\end{enumerate}
\label{tbglobalbifurcation1} \end{theorem}

In order to apply this theory to our problem we need to make several benign alterations to hypotheses (G1)--(G5).  First, we wish to replace $\mathcal{T}$ occurring in (G1) with that from the previous sections, namely the 1-parameter family of laminar flows.  This is justified since the particular form of the trivial solutions plays no role in the proof of \textsc{Theorem \ref{tbglobalbifurcation1}}.  Indeed, once the local curve $\mathcal{C}_{\mathrm{loc}}^{+}$ is found, no more mention need be made of $\mathcal{T}$.  

Secondly, we wish to consider the case where double bifurcation occurs at the local level.  Again, this is possible because the purpose of (G3) is merely to guarantee the existence of $\mathcal{C}_{\mathrm{loc}}^{+}$ and has no bearing on the global analysis.  If, as in section \ref{doubleeigenvaluebifurcationsection}, there is a \emph{set} of local bifurcation curves $\{ \mathcal{C}_{\mathrm{loc}, i} \}_{i =1}^{N}$, then there exists a corresponding set of local parametrizations $\{(\mathfrak{Q}_i, \mathfrak{w}_i)\}_{i = 1}^N$ and so we may define $\{\mathcal{C}_{\mathrm{loc}, i}^{+}\}_{i=1}^N$ in the obvious way.  Applying the argument of \textsc{Theorem \ref{tbglobalbifurcation1}} to each of these curves yields the following corollary.

\begin{corollary}  \emph{(Analytic Global Double Bifurcation)} Let $W, Z, U$ and $\mathcal{F}$ be given as before.  Suppose that the following structural properties are satisfied:
\be  \textrm{there exists a family of trivial solutions } \mathcal{T}  \subset U \cap \mathcal{S} \tag{G1'}, \ee
\be \mathcal{F}_w(Q, w) \textrm{ is a Fredholm operator of index zero for each $(Q, w) \in U \cap \mathcal{S}$}, \tag{G2'} \ee
\be \begin{array}{l} \exists (Q_*, w_*) \in \mathcal{T}, \epsilon > 0, \textrm{ and analytic functions } \{(\mathfrak{Q}_i, \mathfrak{w}_i) \colon (-\epsilon, \epsilon) \to \mathbb{R} \times W\}_{i=1}^N \\ \textrm{ such that, for all $s \in (-\epsilon, \epsilon)$, $i = 1, \ldots, N$,  } \\ (\mathfrak{Q}_i(0), \mathfrak{w}_i(0)) = (Q_*, w_*) \textrm{ and } \mathcal{F}(\mathfrak{Q}_i(s), \mathfrak{w}_i(s)) = 0,\end{array} \tag{G3'} \ee
\be \mathfrak{Q}_i^\prime \nequiv 0 \textrm{ on } (-\epsilon, \epsilon), \qquad i = 1, \ldots, N, \tag{G4'} \ee
\be \textrm{all bounded closed subset of $\mathcal{S}$ are compact.} \tag{G5'} \ee
For each $i = 1, \ldots, N$, define $\mathcal{C}_{\mathrm{loc}, i}^+ := \{ (\mathfrak{Q}_i(s), \mathfrak{w}_i(s)) \colon s \in (0, \epsilon) \}$.  Then there exists a set of continuous curves $\{\mathcal{K}_i\}_{i=1}^N$ with $\mathcal{K}_i$ the unique maximal analytic extension of $\mathcal{C}_{\mathrm{loc}, i}^+$ in the sense of \textsc{Theorem \ref{tbglobalbifurcation1}}, for $i = 1, \ldots, N$.   
\label{tbcorollary1}
\end{corollary}

This result follows directly from the proof of \textsc{Theorem \ref{tbglobalbifurcation1}}, keeping in mind our observation of the last several paragraphs.  For brevity, we therefore do not provide a proof.  Note, however, that this is not a new result: the multiple bifurcation case was considered by Dancer in \cite{dancer1973bifurcation}.  \textsc{Corollary \ref{tbcorollary1}} is merely a translation of Dancer's result to the language of Buffoni and Toland.  

With the background results in place, we now return to the notational conventions of section \ref{globalbifurcationsection} and prove our main result.  Fix $\delta > 0$ and define 
\[ \mathcal{O}_\delta := \left\{ (Q,h) \in \mathbb{R} \times X :~h_p > \delta \textrm{ in } \overline{R},~Q-\sigma\kappa[h] - 2g\rho h > \delta \right\},\]
\[ S := \{ (Q,h) \in \mathbb{R} \times X : \mathcal{G}(Q,h) = 0\}, \qquad S_\delta := \textrm{ closure in } \mathbb{R} \times X \textrm{ of } (S \cap \mathcal{O}_\delta \setminus \mathcal{T}). \] 
By the local analysis of the previous sections, we may take the non-laminar solution set near $(H^*, Q^*)$ to be comprised of a set of $C^1$-curves, $\{ \mathcal{C}_{\mathrm{loc}, i}^\prime \}_{i=1}^N$, for some $ N \geq 2$.  We now prove that these curves are, in fact, analytic.

\begin{lemma} \emph{(Analyticity)} \emph{(a)} The maps $\lambda \mapsto H(\cdot; \lambda) \in C^{3+\alpha}([p_0,0])$ and $\lambda \mapsto Q(\lambda) \in \mathbb{R}$  are analytic on the interval $(-2B_{\mathrm{min}} + \epsilon_0, \infty)$.  Furthermore, if $I$ is any subinterval of $(-2B_{\mathrm{min}}+\epsilon_0, \infty)$ that contains $\lambda_*$, then the operators $\mathcal{F}: I \times X \to Y$ and $\mathcal{G}: \mathbb{R} \times X \to Y$ (defined by \eqref{defF1}--\eqref{defF2} and \eqref{capgravdefG1}--\eqref{capgravdefG2}, respectively) are analytic. \\

\noindent \emph{(b)} There exists $\epsilon > 0$ and a set of analytic functions $\{ (\mathfrak{Q}_i(s), \mathfrak{h}_i(s)) : (-\epsilon, \epsilon) \mapsto \mathbb{R} \times X \}_{i=1}^N$ such that $\mathcal{C}_{\mathrm{loc}, i}^\prime = \{ (\mathfrak{Q}_i(s), \mathfrak{h}_i(s)) : s \in (-\epsilon, \epsilon)\}$, for $i = 1, \ldots, N$.  
\label{analyticitylemma} \end{lemma}  

\begin{proof2}  The analyticity of $\mathcal{G}$ is clear, since 
\[ \partial_Q^j \mathcal{G},~\partial_h^k \mathcal{G} \equiv 0, \qquad \forall j \geq 2,~k \geq 5.\]
Notice also that, since
\[ Q(\lambda) = \lambda + 2g\rho(0)H(0;\lambda), \qquad \mathcal{F}(\lambda, w) = \mathcal{G}(Q(\lambda), w+H(\cdot; \lambda)), \]
it suffices to prove that $\lambda \mapsto H(\cdot; \lambda)$ is analytic in order to obtain part (a). 

We shall do this by extending $\lambda \mapsto H(\cdot; \lambda)$ to a function defined on the complex half-plane $\mathbb{H} := \{ z \in \mathbb{C} : \mathrm{Re}\,z > -2B_{\mathrm{min}} + \epsilon_0 \}$, and then showing that the extended function is complex analytic.  

Mimicking the proof of \textsc{Lemma \ref{laminar2}}, we can prove that, for each $\lambda \in \mathbb{H}$ and $p \in [p_0, 0]$, there exists a unique $H = H(p; \lambda) \in \mathbb{C}$ satisfying the following implicit relationship
\be H(p; \lambda) = \int_{p_0}^p \frac{ds}{\sqrt{\lambda + G(s;\lambda)}}, \qquad p \in [p_0,0], ~\lambda \in \mathbb{H}, \label{extendedHdef} \ee
where we are taking the positive branch of the square root and 
\[ G(p;\lambda) = 2B(p) + 2g\int_{p_0}^p \rho_p(s) \left(H(s;\lambda) - H(0;\lambda)\right) ds.\]
If we denote $\lambda = \lambda_1 + i \lambda_2$ and fix $p \in [p_0,0]$, then
\begin{align*} 
\frac{d}{d{\overline{\lambda}}} H(p;\lambda) &= \frac{1}{2} \left( \frac{d}{d\lambda_1} + i\frac{d}{d\lambda_2}\right) H(p;\lambda) \\
& = -\int_{p_0}^p \left(g\int_{p_0}^r \rho_p(s) \frac{d}{d\overline{\lambda}}\left(H(s; \lambda) - H(0; \lambda)\right)ds\right) \left( \lambda + G(r; \lambda)\right)^{-3/2} \,dr. \end{align*}
The last line above was found by differentiating \eqref{extendedHdef} (with respect to $\lambda_1$ and $\lambda_2$), and making use of the definition of $G$.  

It follows that, if we denote $Y(p;\lambda) := H(p;\lambda) - H(0; \lambda)$, then 
\be \frac{d}{d\overline{\lambda}} Y(p;\lambda) = \int_{p}^0 \left( g \int_{p_0}^r \rho_p(s)\frac{d}{d\overline{\lambda}} Y(s) ds \right) \left( \lambda + G(r;\lambda)\right)^{-3/2}\,dr. \label{dYdlambdabar1} \ee
Recalling the definition of $\epsilon_0$ in \eqref{defepsilon0}, we see that, since $\lambda \in \mathbb{H}$,
\[ \left| \lambda + G(p;\lambda) \right|^{-3/2}  \leq \left| \lambda_1 + 2B_{\mathrm{min}} \right|^{-3/2} \leq (2g |p_0|^2 \|\rho_p\|_\infty )^{-1}. \]
Inserting this last expression into \eqref{dYdlambdabar1}, we find
\[ \left\| \frac{d}{d\overline{\lambda}} Y(\cdot;\lambda) \right\|_\infty \leq \frac{1}{2}  \left\| \frac{d}{d\overline{\lambda}} Y(\cdot;\lambda) \right\|_\infty,\]
and thus $dY/d\overline{\lambda} \equiv 0$.  As $H(p;\lambda) = Y(p;\lambda) - Y(p_0;\lambda)$, this implies further that $dH/d\overline{\lambda} \equiv 0$, hence $\lambda \mapsto H(\cdot;\lambda)$ is complex analytic.  

Now, (b) follows directly from the fact that $\mathcal{F}$ is analytic.    Indeed, simply substituting the Crandall-Rabinowitz and Lyapunov-Schmidt theory of the previous sections with their analytic counterparts, we obtain the analytic parameterizations $\{ (\mathfrak{Q}_i, \mathfrak{h}_i) \}_{i=1}^N$ as claimed.  \end{proof2}

Following the notation of \textsc{Theorem \ref{tbglobalbifurcation1}}, let us denote 
\[ \mathcal{C}_{\delta, i}^{\prime} := \{ (\mathfrak{Q}_i(s), \mathfrak{h}_i(s)) : s \in (0, \epsilon) \}, \qquad i = 1,\ldots, N.\]
Note that, to make our notation more concise, we have dropped the $+$ superscript. Our main result for this section is the following.

\begin{theorem} \emph{(Analytic Global Bifurcation)}  Fix $\delta > 0$.  There exists a maximal, global, continuous extension $\mathcal{K}_{\delta,i}^\prime \subset S_\delta$ of $\mathcal{C}_{\delta, i}^\prime$ (in the sense that \textsc{Theorem \ref{tbglobalbifurcation1}}\emph{(a)--(f)} holds for $\mathcal{K}_{\delta,i}^\prime$), for $i = 1, \ldots, N$.  Furthermore, for each $i$, one of the following alternatives must hold.  
\begin{indentenum}
\item[\emph{(i)}] $\| (\mathfrak{Q}_i(s), \mathfrak{h}_i(s))\|_{\mathbb{R} \times X} \to \infty$ as $s \to \infty$; or
\item[\emph{(ii)}] the closure of $\mathcal{K}_{\delta, i}^\prime$ contains a point $(Q,h) \in \partial \mathcal{O}_\delta$; or
\item[\emph{(iii)}] $\mathcal{K}_{\delta, i}^\prime$ is a closed loop.
\end{indentenum}
\label{globalbifurcationdoubleeigenvalue}
\end{theorem}

\begin{longremark} In several respects, the alternatives above are considerably stronger than their analogues in \textsc{Theorem \ref{capgravglobalbifurcationtheorem}}.  If, recalling the notation of section \ref{globalbifurcationsection}, we define $\mathcal{C}_\delta^\prime$ to be the connected component of $S_\delta$ containing $\bigcup_i \mathcal{C}_{\mathrm{loc}, i}^\prime$, then clearly $\bigcup_i \mathcal{K}_{i, \delta}^\prime \subset \mathcal{C}_\delta$.  In other words, the conclusions of \textsc{Theorem \ref{globalbifurcationdoubleeigenvalue}} applies only to a (very special) subset of $\mathcal{C}_\delta^\prime$, whereas the alternatives of \textsc{Theorem \ref{capgravglobalbifurcationtheorem}} apply to the \emph{entirety} of $\mathcal{C}_\delta^\prime$.   In that sense, the theorem above is stronger, since it identifies --- in fact, it constructs --- the precise path along which the alternative occurs.  Moreover, the fact that $\bigcup_i \mathcal{K}_{\delta, i}^\prime$ is path-connected and locally analytic is noteworthy in and of itself.  This cannot be obtained through purely degree theoretic methods, since it requires more than topological information about the solution set.  

Still, it is not possible to completely reconcile the alternatives of \textsc{Theorem \ref{capgravglobalbifurcationtheorem}} and \textsc{Theorem \ref{globalbifurcationdoubleeigenvalue}}.  For instance, it is entirely permissible for $\|(\mathfrak{Q}_i(s), \mathfrak{h}_i(s))\| \to \infty$ as $s \to \infty$, but, for some $s > 0$, we have $(\mathfrak{Q}_i(s), \mathfrak{h}_i(s)) \in \mathcal{T}$.   Moreover, it may be that $\mathcal{C}_{\delta, i}^\prime$ is a closed loop for all $i = 1,\ldots, N$, yet $\bigcup_i \mathcal{K}_{\delta, i}^\prime \cap \mathcal{T} = (Q_*, H_*)$.  This is possible, again, because $\bigcup_i \mathcal{K}_{i,\delta}^\prime$ is only a subset of $\mathcal{C}_\delta^\prime$.

Finally, let us note that a similar result can be proved for the simple bifurcation case --- we need only use \textsc{Theorem \ref{tbglobalbifurcation1}} directly, rather than \textsc{Corollary \ref{tbcorollary1}}.  
\end{longremark}

\begin{proofof}{Theorem \ref{globalbifurcationdoubleeigenvalue}.} In light of \textsc{Corollary \ref{tbcorollary1}} and \textsc{Lemma \ref{analyticitylemma}}, the structure of our argument is fairly straightforward:  we need only confirm that $\mathcal{G}$ satisfies structural properties (G1')--(G5'). 

With that in mind, let $X,\,Y, \,\mathcal{O}_\delta$ and $\mathcal{G}$ stand in for $W, \,Z, \,U$ and $\mathcal{F}$  of \textsc{Corollary \ref{tbcorollary1}}, respectively.  Then the family of laminar flow, $\mathcal{T}$ restricted to $\mathcal{O}_\delta$ satisfies (G1').  Condition (G2') was verified in \textsc{Lemma \ref{capgravfredholmlemma}}, while (G3') and (G4') follow from the analysis of section \ref{capgravlocalbifurcationsection} and section \ref{doubleeigenvaluebifurcationsection} (in particular, note that since the mixed solution curves obey the scaling $|\xi| \leq C |\lambda - \lambda_*|$, for some $C > 0$, we cannot have $\mathfrak{Q}^\prime \equiv 0$.)  Finally, (G5') is a consequence of \textsc{Lemma \ref{capgravpropermaplemma}}, since $S_\delta = \mathcal{G}^{-1}(\{0\}) \cap \mathcal{O}_\delta$.  We are therefore justified in taking $\mathcal{K}_{\delta, i}^\prime$ to be the global extension of $\mathcal{C}_{\mathrm{loc}, i}^\prime$, for $i = 1, \ldots, N$.  Alternatives (i)--(iii) follow from \textsc{Theorem \ref{tbglobalbifurcation1}}(e) applied to $\mathcal{K}_{\delta, i}^\prime$. \end{proofof}

Ideally, one would like to eliminate all but the first and second alternatives of \textsc{Theorem \ref{globalbifurcationdoubleeigenvalue}}, as this would imply that either the curve is unbounded, or the problem degenerates (i.e. loses ellipticity) along $\bigcup_{\delta > 0} \mathcal{K}_i^\prime$.  As we shall see in section \ref{finalsection}, put into the Eulerian framework both of these possibilities lead to us to conclude that there is a curve of solutions to the original problem \eqref{incompress}-\eqref{bedcond} along which either the fluid tends to stagnation at a point, or  the horizontal velocity (in the moving frame) approaches $-\infty$ somewhere within the fluid.  

Typically one excludes the undesirable alternative (iii) by means of a nodal property argument.  That is, using maximum principles one shows that the sign of certain quantities are invariant along the curve, and then argue that this is violated if $\mathcal{K}_i^\prime$ is a closed loop (or, if we are considering simple bifurcation, if $\mathcal{C}_\delta^\prime$ returns to $\mathcal{T}$) (see, for example, \cite{dancer1973globalsolution}).  Taking $\sigma = 0$, an argument of this kind was successfully carried out in the constant density case (cf. \cite{constantin2004exact}), and, to a lesser degree, in the stratified case (cf. \cite{walsh2009stratified}). Unfortunately, the introduction of surface tension--- in particular, the higher derivatives on the boundary that accompanies it--- seems to destroy the maximum principle structure of the problem.  The rather conspicuous dearth of literature on unbounded curves of capillary and capillary-gravity waves is a direct consequence of this fact.

\section{Uniform regularity} \label{uniformregularitysection}

In this section we consider the first alternatives of \textsc{Theorem \ref{capgravglobalbifurcationtheorem}} and \textsc{Theorem \ref{globalbifurcationdoubleeigenvalue}}.  As in \cite{constantin2004exact}, our goal is the establishment of uniform bounds in the H\"older norm along the continuum $\mathcal{C}_\delta^\prime$, for the simple eigenvalue case, and along the curves $\{\mathcal{K}_{\delta, i}^\prime\}_{i=1}^4$ for the higher multiplicity case.  In that paper Constantin and Strauss were able to control the $C^{3+\alpha}$-norm along the continuum in terms $\sup_{h \in \mathcal{C}_\delta^\prime} \sup_{\overline{R}} h_p$ (assuming that $Q$ is likewise bounded along the continuum).  However, the addition of the curvature term $\kappa[h]$ to the nonlinear boundary condition complicates matters: we must now allow for the possibility that the curvature may blowup (in $L^\infty$) along the continuum.
\begin{theorem} \emph{(Uniform Regularity)} Let $\delta > 0$ be given and assume that $\beta \in C^{3+\alpha}((0,|p_0|))$.  If 
\[ \sup_{(Q,h) \in \mathcal{C}_\delta^\prime} \sup_{\overline{R}} h_p < \infty, \qquad\sup_{(Q,h) \in \mathcal{C}_\delta^\prime} Q < \infty \qquad \textrm{and} \qquad \inf_{(Q,h) \in \mathcal{C}_\delta^\prime} \inf_{T} \kappa[h] > -\infty,\]
then $\sup_{(Q,h) \in \mathcal{C}_\delta^\prime}\|h\|_{C^{3+\alpha}(\overline{R})}$ is finite. \label{capgravuniformregularitytheorem} \end{theorem} 

In the proof of \textsc{Theorem \ref{capgravuniformregularitytheorem}}, the most technically complicated step is obtaining bounds on the second derivatives of $h$ along $\mathcal{C}_\delta^\prime$.  The $L^\infty$ estimates will be a consequence of $p_0$ being constant along the continuum, while the higher-order terms can be expressed (via the height equation) in terms of the lower-order terms  in the usual way. To treat the second-order terms, we shall employ a suite of \emph{a priori} estimates for nonlinear elliptic equations with Venttsel boundary conditions due to Luo and Trudinger (cf. \cite{luo1994quasilinear}).  We shall not, however, need these results in their full generality, so instead we streamline them into a single statement.  Of course, in order to apply this to our problem, we will again need to use the method of freezing the operator $\mathcal{G}$.  
 
First note that due to the periodicity in $q$, the domain $R$ can be considered as embedded on a torus, allowing us to ignore the seeming lack of smoothness at the corner points $(0,0),\,(0,p_0),\,(2\pi,0)$, and $(2\pi,p_0)$.  We consider a differential operator $F = F(x, z, \xi, r) \in C^2(R \times \mathbb{R}\times \mathbb{R}^2 \times \mathbb{S}, \mathbb{R})$ and boundary operator $G = G(x, z, \xi,r) \in C^2(R \times \mathbb{R}\times \mathbb{R}^2\times \mathbb{S}, \mathbb{R})$.  Here $\mathbb{S}$ denotes the space of $2\times 2$ real symmetric matrices, $D$ denotes the gradient operator and $D^2$ the Hessian. Both of these are taken as acting on the space of smooth real-valued functions on $\overline{R}$ which are $2\pi$-periodic in the first variable.  As usual, we say that $F$ is elliptic at $(z, \xi, r) \in  \mathbb{R} \times \mathbb{R}^2 \times \mathbb{S}$ provided that the matrix $F_r := [ \partial F/ \partial r_{ij}]_{1 \leq i, j \leq 2}$ is positive definite at that point.    Moreover, if $\Lambda_1$ and $\Lambda_2$ denote the minimum and maximum eigenvalue of $F_r$, respectively, then $F$ is said to be uniformly elliptic provided that the ratio $\Lambda_2/\Lambda_1$ is bounded.

We now specialize this framework to the problem at hand. Let $a_{ij} = a_{ij}(z,\xi)$, $b = b(x,z,\xi) \in C^2(R \times \mathbb{R} \times \mathbb{R}^2, \mathbb{R})$ and $\mathfrak{a}_{ij} \in \mathbb{S}$, $\mathfrak{b} = \mathfrak{b}(z, \xi) \in C^2(\mathbb{R} \times \mathbb{R}^2, \mathbb{R})$ be given.  We consider the case where $(F,G)$ are of the form
\[ F(z,\xi,r) = a_{ij}(z,\xi) r_{ij} + b(x,z,\xi), \qquad G(z,\xi,r) = \mathfrak{a}_{ij} r_{ij} + \mathfrak{b}(z,\xi).\]
Here $a_{ij}$ is assumed to be such that $F$ is uniformly elliptic, in the sense of the previous paragraph.  We shall take $G$ to be of Venttsel-type, which means the following:  (i) (Uniform obliqueness) There exists some $\chi > 0$ such that, at each point $(q,0)$, the inward normal derivative $\mathfrak{b}_\xi \cdot (0,1) \geq \chi > 0$, for all $(z,\xi) \in \mathbb{R}\times\mathbb{R}^2$; (ii) (Ellipticity) $\mathfrak{a}_{ij} \zeta_i \zeta_j > 0$, for all $\zeta = (\zeta_1, \zeta_2) \in \mathbb{R}^2$ such that $\zeta$ is tangent to the lower boundary of $R$.

With $F$ and $G$ given as above we wish to consider the quasilinear elliptic (Venttsel) boundary value problem 
\be \left \{ \begin{array}{ll}
F(h,Dh, D^2 h) = 0, & \textrm{in } R \\
G(h, Dh, D^2 h) = 0 & \textrm{on } p = 0 \\
h = 0 & \textrm{on } p = p_0. \end{array} \right. \label{luotrudbvp} \ee
 
\begin{theorem} \emph{(Luo-Trudinger)} Let $h \in C^5(\overline{R})$ be a solution to \eqref{luotrudbvp} that is $2\pi$-periodic in the first variable and for some constant $K > 0$ satisfies $|h|+|Dh| \leq K$ in $\overline{R}$. Suppose further that for some positive constant $M$ the functions $F = F(z, \xi, r)$ and $G = G(z,\xi)$ satisfy the following structural conditions:
\begin{eqnarray} \Lambda_2 & \leq & M \Lambda_1, \label{structure1} \\
 |F|,~|F_\xi|,~|F_z| &\leq& M \Lambda_1 (|r|+1), \label{structure2} \\
 F_{rr} &\leq& 0,  \label{structure3} \\
|G|, |G_z|, |G_\xi|, |G_{zz}|, |G_{z\xi}|, |G_{\xi\xi}| &\leq& M \chi, \label{structure4} \\
|\partial_z^2 a_{ij}|, (1+|\xi|)|\partial_z\partial_\xi a_{ij}|, (1+|\xi|^2)|\partial_\xi^2 a_{ij}| & \leq & M \Lambda_1 \label{structure5} \\
|b_x|, |b_z|, |b_\xi|(1+|\xi|) & \leq & M\Lambda_1 (1+|\xi|^2) \label{structure6} \\
|b_{xx}|, |b_{zx}|, |b_{zz}|, (1+|\xi|)|b_{z\xi} |, (1+|\xi|)|b_{x\xi}|, (1+|\xi|^2)|b_{\xi\xi}| &\leq& M \Lambda_1 (1+|\xi|^2) 
\label{structure10} \\
|\mathfrak{b}_{zz}|,(1+|\xi|)|\mathfrak{b}_{z\xi}|, (1+|\xi|) |\mathfrak{b}_{\xi\xi}| &\leq& M\chi (1+|\xi|) \label{structure7} \\
|D^3 a_{ij}|,|D^3 b| & \leq & M \Lambda_1 \label{structure8} \\
|D^3 \mathfrak{b}_| & \leq & M \chi \label{structure9}
\end{eqnarray}
for all $(z, \xi, r) \in \mathbb{R} \times \mathbb{R}^2 \times \mathbb{S}$ such that $|z|+|\xi| \leq K$.  Then there are positive constants $\mu = \mu(M)$ and $C = C(K, M)$ such that $h \in C_{\textrm{per}}^{2+\mu}(\overline{R})$ and
\[ \|h \|_{C_{\textrm{per}}^{2+\mu}(\overline{R})} \leq C.\]
\label{luotrudtheorem}
\end{theorem}

This is a (weaker) restatement of the results presented in \cite{luo1994quasilinear}.  Also, for consistency with presentation of similar theorems earlier in \cite{constantin2004exact} and \cite{walsh2009stratified}, we have slightly altered  notation.  We now give a brief sketch of the proof.  

\begin{proofof}{Theorem \ref{luotrudtheorem}.} It was already shown in \cite{walsh2009stratified} that \eqref{structure1}-\eqref{structure4} together amount to the ``natural structural conditions'' (L-T F1)--(L-T F5) presented in \cite{luo1994quasilinear}.  We will now show that, unsurprisingly, these further imply structure conditions (Luo-T A), (Luo-T A1), (Luo-T A2), (Luo-T B), (Luo-T B1) and (Luo-T B2) of \cite{luo1994quasilinear}.

First consider the conditions on $F$.  That (L-T F1) is equivalent to (Luo-T A) is obvious.   Since $F(\cdot, \cdot, \cdot, r = 0) = b(\cdot, \cdot, \cdot)$, (LT F2) immediately leads to (Luo-T A1).  Similarly, if we note that $D_r F = a_{ij}$, it can be readily verified that condition (L-T F4) implies the first line of inequalities in (Luo-T A2).  To get the second line, we appeal to (L-T F3), noting once more that $F$ evaluated at $r = 0$ is $b$.  

Next we turn to $G$.  Clearly, (Luo-T B) follows from uniform obliqueness.  Since $G(\cdot, \cdot, \cdot, r = 0) = \beta(\cdot, \cdot, \cdot)$, evaluating (L-T G2)-(L-T G3) at $r = 0$ yields the growth conditions on $\beta$ in (Luo-T B1)-(Luo-T B2).

Combining this with the new structural conditions \eqref{structure5}-\eqref{structure9}, the hypothesis of our theorem encompasses all the conditions of \cite{luo1994quasilinear}.  To arrive at the statement above, we simply apply in succession each of the results in \S2--\S4 of that paper.  \end{proofof}

Our tools fully assembled, we are ready to prove \textsc{Theorem \ref{capgravuniformregularitytheorem}}.

\begin{proofof}{Theorem \ref{capgravuniformregularitytheorem}.}
We begin by proving the $L^\infty(\overline{R})$ norm of $h$ is uniformly bounded along the continuum.  To do so we note that by \eqref{defh}
\[ 0 \leq h(q,p) = y+d \leq \eta(0)+d,\]
and thus the definition of $p_0$ in \eqref{defp0} yields
\[ \eta(0)+d \leq \frac{|p_0|}{ \inf_{\overline{D_\eta}}\left( \sqrt{\rho}(c-u)\right)} = |p_0| \sup_{\overline{R}} h_p \]
Thus assuming the hypothesis of the theorem, that is that $h_p$ is uniformly bounded in $L^\infty$ along the continuum, then the above inequalities together give the boundedness of the $L^\infty$-norm of $h$.

To bound $h_q$ we proceed as in \textsc{Theorem 6.1} of \cite{walsh2009stratified} to show that $h_q$ cannot attain an interior minimum or maximum.   Indeed, since the interior equation here is identical to the one considered in that paper, a verbatim application of their argument suffices.  As $h_q \equiv 0$ on the bottom, we need only concern ourselves with the top.  The nonlinear boundary condition then gives
\[ \inf_T h_p^2(Q-2\sigma \kappa[h]-2g\rho h) \leq 1+h_q^2 \leq \sup_T h_p^2(Q-2\sigma \kappa[h]-2g\rho h). \]
The quantity on the left-hand side is bounded uniformly from below by $\delta^3$ on $\mathcal{C}_\delta^\prime$.  Thus $\sup_{\overline{R}} h_q = \sup_{T} h_q$ can be controlled by $\sup_{\mathcal{C}_\delta^\prime}\sup_{\overline{R}} h_p$ and $\sup_{\mathcal{C}_\delta^\prime} \sup_T (-\kappa[h])$.

For the second-order terms we need to make use \textsc{Theorem \ref{luotrudtheorem}}, and so we shall once again freeze the operator $\mathcal{G}$ in the sense of section \ref{globalbifurcationsection}.  Fix $K > 0$, $\tau \in \mathbb{R}$.  Let $\zeta = \zeta(\xi)$ be a smooth cutoff function with $0 \leq \zeta \leq 1$,  $\zeta \equiv 1$ on $B_K(0) \cap \{ \xi : \xi_2 > \delta \}$ and $\textrm{supp } \zeta \subset B_{2K}(0) \cap \{ \xi : \xi_2 > \delta/2 \}$.  For each $x = (x_1, x_2) \in \mathbb{R}^2$, $z \in \mathbb{R}$, $\xi = (\xi_1, \xi_2) \in \mathbb{R}^2$, put 
\begin{align*} a_{ij}(\xi) & := \left(\begin{array}{cc} 
\xi_2^2 & - \xi_1 \xi_2 \\ -\xi_1 \xi_2 & 1+\xi_1^2 \end{array} \right) \zeta(\xi) 
+ \left(\begin{array}{cc} 
1 & 0 \\ 0 & 1 \end{array} \right)(1-\zeta(\xi)), \\
b(x, z, \xi) & := \xi_1^3\left(\beta(-x_2)-g(z-\tau)\rho_p(x_2) \right)\zeta(\xi) , \end{align*}
and
\[ \mathfrak{a}_{ij} := \left(\begin{array}{cc} 2\sigma & 0 \\ 0 & 1 \end{array} \right), \qquad \mathfrak{b}(z,\xi) := \left(1+\xi_1^2\right)^{\frac{3}{2}} \left(Q - 2g\rho(0) z -\frac{1+\xi_1^2}{\xi_2^2}\right) \zeta(\xi) +(1-\zeta(\xi))\xi_2 .\]
Suppose that for some $\tau$, $(h,Q)$ is a solution to \eqref{luotrudbvp} with $F$ and $G$ as in \textsc{Theorem \ref{luotrudtheorem}}, that is,
\begin{eqnarray*} 
F(x,z,\xi,r) &=& \left(\xi_2^2 r_{11} +(1+\xi_1^2)r_{22} -2 \xi_1\xi_2 r_{12} + \xi_1^3 \left(\beta(-x_2) -g(z-\tau) \rho_p(x_2)\right) \right) \zeta(\xi) \\
& & + r_{11} (1-\zeta(\xi)) \\
 G(z,\xi,r) &=& 2\sigma r_{11} + \left( (1+\xi_1^2)^{3/2} \left(Q-2g\rho(0)z-\frac{1+\xi_1^2}{\xi_2^2}\right)\right) \zeta(\xi) + \xi_2 (1-\zeta(\xi)).\end{eqnarray*}
It follows that each $(h,~Q) \in \mathcal{C}_\delta^\prime$ with $|h|+|Dh| < K$ satisfies the Venttsel boundary value problem \eqref{luotrudbvp} for the above choice of $F$ and $G$. The next task is verifying the structural conditions.  It is easy to see that $F$ is uniformly elliptic.  To show uniform obliqueness note that, at any point $x \in \partial R$, the inward normal is given by:
\[ \mathfrak{b}_\xi \cdot (0,1) =  \frac{2\left(1+\xi_1^2\right)^{\frac{5}{2}}}{\xi_2^3} \zeta(\xi) - (\partial_{\xi_2} \zeta(\xi))\xi_2 + 1-\zeta(\xi). \]
Inside the region where $\zeta \equiv 1$, this is clearly uniformly positive.  Also, since we may take $\zeta$ to be monotonic decreasing in the $\xi_2$-direction, the last three terms guarantee that the innward normal derivative is uniformly bounded away from zero in the entirety of $\mathbb{R}^2$.  

The final structural conditions are easy.  Since $a_{ij}$, $b$, $\mathfrak{a}_{ij}$ and $\mathfrak{b}$ are all suitably smooth, by compactness we can choose a sufficiently large $M$ so that \eqref{structure1}-\eqref{structure9} hold for $\xi \in \overline{B_K(0)} \cap \{ \xi : \xi_2 \geq \delta \}$.  (Note that here we make use of the assumption $\gamma \in C^{3}$ to get \eqref{structure8}.)

The more serious step is to guarantee that $h$ has enough regularity to begin with.  That is, we know $h \in C^{3+\alpha}(\overline{R})$, however \textsc{Theorem \ref{luotrudtheorem}} requires it to be of class $C^5$.  In order to conclude that $\mathcal{C}_\delta^\prime$ lies in $C^5(\overline{R})$ we shall apply standard linear elliptic theory.   Setting aside the question of $\beta$ for the moment,  the fact that $F$ is quasilinear tells us that, if $h \in C^{k+\alpha}(\overline{R})$ with $k > 2$, then the coefficients of the problem are at worst $C^{k-1+\alpha}(\overline{R})$.    Thus, if we assume $\beta \in C^{3+\alpha}((0,|p_0|))$, we have $a_{ij},~b \in C^{2+\alpha}$, which implies by the estimates of \textsc{Lemma \ref{estimateslemma}} for the linear Venttsel problem that $h \in C^{4+\alpha}(\overline{R})$.  This means, in turn, that the coefficients are actually $C^{3+\alpha}$.  Turning the crank once more, we conclude $h \in C^{5+\alpha}(\overline{R})$. \textsc{Theorem \ref{luotrudtheorem}} is therefore applicable, and thus the second order derivatives are bounded.

Thus $(F,G)$ satisfy the structural requirements of \textsc{Theorem \ref{luotrudtheorem}}.   We conclude that there exists a constants $C = C(K, M, \tau)$ and $\mu = \mu(K,M, \tau)$ so that, for any solution $(h,Q)$ of \eqref{luotrudbvp} with $|h| + |Dh| < K$,
\[ \|h\|_{C^{2+\mu}(\overline{R})} \leq C. \]
Note that every $h \in \mathcal{C}_\delta^\prime$ is such a solution for $\tau = d(h)$.  In order to complete the estimates of the second-order terms, therefore, requires removing the dependency on $\tau$ from $C$.  Toward that end we recall that $\tau = d(h) \leq \sup_{R} |h|$. Thus we may select $K$ independently of $\tau$.  Moreover, since in the definition of $F$, $\tau$ does not occur in any of the second-order terms $a_{ij}$, it does not affect $\Lambda_1$ and $\Lambda_2$.  Indeed, the characteristic polynomial satisfied by the eigenvalues is
\[ \Lambda_i^2 - (1+\xi_1^2+\xi_2^2)\Lambda_i - (1+\xi_1^2)\xi_2^2 - 4 \xi_1^2 \xi_2^2 = 0, \qquad i =1,2\]
which does not depend on $\tau$.  This is an obvious consequence of the fact that the effects of density variation in the height equation are consigned to the lower-order terms.  The term in $F$ that does contain $\tau$ is $b$, and it can be easily estimated in terms of $K$.  Thus $C$ and $\mu$ can be chosen independently of $\tau$, which gives the second-derivative bounds.

Finally we establish third derivative bounds.  Denoting $\theta := h_q$, we differentiate the height equation in $q$ to find,
\be \left \{ \begin{array}{ll}
(1+h_q^2) \theta_{pp} -2h_p h_q \theta_{pq} + h_p^2 \theta_{qq} & \\
\qquad = C(h, h_p, h_q, h_{pp}, h_{pq}, h_{qq}), & \textrm{in } R, \\
& \\
2\sigma h_p^2 \left(\theta_q(1+h_q^2)^{-3/2}\right)_q + 2h_q \theta_q + (2g\rho h + 2\sigma \kappa[h] - Q)h_p \theta_p & \\ 
\qquad = D(h_q, h_p) & \textrm{on } T, \\
& \\
\theta = 0 & \textrm{on } B. \end{array} \right. \label{unifregularitythetaeq} \ee
Thus $\theta$ solves a second-order elliptic PDE with $C_{\textrm{per}}^{1+\mu}(\overline{R})$ coefficients and right-hand sides $C \in C_{\textrm{per}}^\mu(\overline{R})$, and $D \in C_{\textrm{per}}^{\mu}(T)$.  Since we have already proved $h$ is bounded uniformly in $C_{\textrm{per}}^{2+\mu}(\overline{R})$, it follows that the right-hand side of \eqref{unifregularitythetaeq} is bounded uniformly in $C_{\textrm{per}}^\mu(\overline{R})$ along the continuum.  As we have seen in section \ref{globalbifurcationsection}, so long as $h_p$ is bounded uniformly away from 0 along the continuum, problem \eqref{unifregularitythetaeq} above will be uniformly elliptic with a nondegenerate Venttsel-type boundary condition.  We may therefore apply the linear theory for such problems (e.g. \textsc{Lemma \ref{estimateslemma}}) to obtain uniform \emph{a priori} estimates of $\theta$ in $C_{\textrm{per}}^{2+\mu}(\overline{R})$ over all of $\mathcal{C}_\delta^\prime$.

It remains only to bound $h_{ppp}$ in $C_{\textrm{per}}^\mu(\overline{R})$.  By means of the height equation, we may express $h_{pp}$ in terms of the lower order derivatives of $h$:
\[ h_{pp} = -(1+h_q^2)^{-1} \left( h_p^3 \gamma(-p) + h_{qq}h_p^2 - 2h_q h_p h_{pq} \right) \]
But the right-hand side above is in $C^{1+\mu}_{\textrm{per}}(\overline{R})$ by the arguments of the previous paragraphs.  Thus $h_{ppp}$ in $C^{\mu}_{\textrm{per}}(\overline{R})$ as desired.  

To transition back to the original H\"older exponent, $\alpha$, we merely note that if $h \in C^{3+\mu}_{\textrm{per}}(\overline{R})$, then it is certainly in $C^{2+\alpha}_{\textrm{per}}(\overline{R})$.  The arguments we have used to derive the third derivative bounds in no way relied on the particular value of $\mu$, so running through them with $\alpha$ instead we find $h \in C^{3+\alpha}_{\textrm{per}}(\overline{R})$.  \end{proofof}

\section{Proof of the main results} \label{finalsection}

With the uniform regularity results we established in the previous section, we are now prepared to begin the task of proving the main theorems.  \textsc{Theorem \ref{capgravglobalbifurcationtheorem}} gave us three possibilities for the continuation of the local bifurcation curve(s).   In order to prove \textsc{Theorem \ref{mainresult1}}, we need merely to make a careful study of these alternatives in light of \textsc{Theorem \ref{capgravuniformregularitytheorem}}.

\begin{proofof}{Theorem \ref{mainresult1}.}  Fix $\delta > 0$.  If $\mathcal{C}_\delta^\prime$ is unbounded in $\mathbb{R}\times X$, then at least one of the following must occur:
\begin{indentenum}
\item[(i)] there exists a sequence $(Q_n, h_n) \in \mathcal{C}_\delta^\prime$ with $\lim_{n\to\infty} Q_n = \infty$; or
\item[(ii)] there exists a sequence $(Q_n, h_n) \in \mathcal{C}_\delta^\prime$ with $\lim_{n\to\infty} \min_{T} \kappa[h_n]  = -\infty$; or
\item[(iii)] there exists a sequence $(Q_n, h_n) \in \mathcal{C}_\delta^\prime$ with $\lim_{n\to\infty} \max_{\overline{R}} \partial_p h_n = \infty$,
\end{indentenum}
where the second and third possibilities are a consequence of \textsc{Theorem \ref{capgravuniformregularitytheorem}}.  If $\mathcal{C}_\delta^\prime$ contains a point of $\partial\mathcal{O}_\delta$, then either
\begin{indentenum}
\item[(iv)] there exists a $(Q, h) \in \mathcal{C}_\delta^\prime$ with $\partial_p h = \delta$ somewhere in $\overline{R}$, or
\item[(v)] there exists a $(Q, h) \in \mathcal{C}_\delta^\prime$ with $2g\rho h + 2\sigma \kappa[h] = Q-\delta$ somewhere on $T$.  
\end{indentenum} 
Thus, \textsc{Theorem \ref{capgravglobalbifurcationtheorem}} tell us that one of the five alternatives above occurs, or else the bifurcation curve intersects $\mathcal{T}$ at a point other than $(Q^*, H^*)$.

Consider the first alternative (i).  By the definition of $p_0$ in \eqref{defp0}, we can estimate
\[ \inf_{y \in [-d, \eta(x)]} \sqrt{\rho(x,y)}\Big(c-u(x,y)\Big) \leq \frac{p_0}{\eta(x)+d} \leq \sup_{y \in [-d, \eta(x)]} \sqrt{\rho(x,y)}\Big(c-u(x,y)\Big), \]
and thus, appealing to \eqref{defQ}, we deduce
\be
\begin{split} Q &= \rho(0, \eta(0))\Big(c-u(0,\eta(0))\Big)^2 + 2g\rho(0, \eta(0))(\eta(0)+d) \\
& \leq \sup_{\overline{D_\eta}} \rho (c-u)^2 + \frac{2g \rho(0, \eta(0)) p_0}{\inf_{\overline{D_\eta}} \sqrt{\rho}(c-u)}. \end{split} \label{boundQ} \ee
If for some $\delta > 0$ the first alternative holds, we conclude from \eqref{boundQ} that for the corresponding sequence $(u_n, v_n, \rho_n,\eta_n)$ of solutions to \eqref{euler2}-\eqref{boundcond}, either $\sup_{\overline{D_{\eta_n}}} \sqrt{\rho_n} (c-u_n) \to \infty$, or $\inf_{\overline{D_{\eta_n}}} \sqrt{\rho_n}(c-u_n) \to 0$.  

Similarly, if the third alternative (iii) holds, then, simply by rearranging the change of variables equations in \eqref{uvhphq}, we have
\[ \partial_p h_n = \frac{1}{\sqrt{\rho_n} (c-u_n)},\]
whence $\max_{\overline{R}} \partial_p h_n \to \infty$ if and only if $\inf_{\overline{D_{\eta_n}}} \sqrt{\rho_n}(c-u_n) \to 0$.  

Now consider the second alternative (ii) wherein, along some sequence, we have that a point of infinite (negative) curvature is developing.  Then, by the dynamic boundary condition \eqref{presscond}, we see that 
\[ P_n - P_{\mathrm{atm}} = -\sigma \eta_n^{\prime\prime} \left(1+(\eta_n^\prime)^2\right)^{-3/2},\]
and so $\max_{\eta_n} (P_n - P_{\mathrm{atm}}) \to \infty$.  In other words, as the free surface develops a corner point, the pressure at that points (relative to the atmospheric pressure) goes to positive infinity. This, we anticipate, leads to a blow-up in the energy on the free surface.  Indeed, from the definition of $Q$ in the Eulerian framework, we obtain
\begin{align*}
Q_n &= 2\left(E_n|_{\eta_n} - P_{\mathrm{atm}} + g\rho_n|_{\eta_n} d_n\right) \\
& = 2\left(P_n - P_{\mathrm{atm}}\right) + \left[ \rho_n \left((u_n-c)^2+v_n^2\right) + 2g\rho_n \left(\eta_n+d_n\right)\right]\bigg|_{\eta_n} \\
& > 2\left(P_n - P_{\mathrm{atm}}\right).\end{align*}
Thus alternative (ii) reduces to alternative (i).  

Next, suppose that for some decreasing sequence $\delta_n \to 0$ the fourth (iv) alternative holds.  Again, appealing to \eqref{uvhphq}, we see immediately that the corresponding sequence in the Eulerian framework, $(u_n, v_n, \rho_n,\eta_n)$, must satisfy $\sup_{\overline{D_{\eta_n}}} \sqrt{\rho_n} (c-u_n) \to \infty$ as $n \to \infty$.  

Finally, if we have a sequence $\delta_n \to 0$ for which the fifth alternative (v) holds, then there exists a sequence $\{(Q_n, h_n)\}$ with $(Q_n, h_n) \in \mathcal{C}_{\delta_n}^\prime$ for each $n \geq 1$ and $\sup_T (Q_n-2g\rho h_n 2\sigma \kappa[h_n]) \to 0$.  If we evaluate the boundary condition on $T$ for this sequence, we find
\[ \max_{T} \left(\frac{1}{(\partial_p h_n) ^2}\right) \leq \max_{T} \left(\frac{1+(\partial_q h_n)^2}{(\partial_p h_n)^2}\right) = \max_T \left( Q_n -2g\rho_n h_n-2\sigma \kappa[h_n] \right). \]
We infer, therefore, that $\max_{\overline{R}} \partial_p h_n \to \infty$, and hence the fifth alternative implies the third.

By construction, the family of continua $\mathcal{C}_\delta^\prime$ indexed by $\delta > 0$ is increasing as $\delta$ decreases.  We may therefore define $\mathcal{C}^\prime := \bigcup_{\delta > 0} \mathcal{C}_\delta^\prime$ to be the maximal continuum.  Then by the considerations of the \S2, and in particular \textsc{Lemma \ref{equivalencelemma}}, there exists a connected set $\mathcal{C}$ of solutions to \eqref{incompress}-\eqref{bedcond} corresponding to $\mathcal{C}^\prime$.  The arguments of the preceding paragraphs, moreover, imply that along some sequence in $\mathcal{C}$, either (i) $\sup_{\overline{D_{\eta_n}}} u_n \to c$, or  (ii) $\inf_{\overline{D_{\eta_n}}} u_n \to -\infty$, or (iii) $\mathcal{C}$ contains more than one laminar flow solution.  This proves part (a).  

Part (b) follows from part (a), taking into account the remarks following \textsc{Theorem \ref{globalbifurcationdoubleeigenvalue}}.  That is, for each $\delta > 0$, we have a path-connected set $\mathcal{K}_\delta^\prime \subset \mathcal{C}_{\delta}^\prime$ where, along some path one of alternative (i)--(v) occurs or $\mathcal{K}_\delta^\prime$ is a closed loop.  Arguing as before, it is clear that, if $\mathcal{K}^\prime := \bigcup_{\delta > 0} \mathcal{K}_\delta^\prime \subset \mathcal{C}^\prime$, then there is a corresponding path-connected set $\mathcal{K} \subset \mathcal{C}$ of solutions to \eqref{incompress}-\eqref{bedcond}.  It follows that either $\mathcal{K}$ is a closed loop, or along some path in $\mathcal{K}$ we have either (i) $\sup_{\overline{D_{\eta}}} u \to c$, or  (ii) $\inf_{\overline{D_{\eta}}} u \to -\infty$.   

Finally, (c) is an automatic consequence of the definition of the space $X$.  The proof is complete.\end{proofof}

\begin{proofof}{Theorem \ref{mainresult3}.}
Fix $\delta > 0$ and $i \in \{1,2,3,4\}$.  In view of \textsc{Theorem \ref{globalbifurcationdoubleeigenvalue}} and \textsc{Theorem \ref{capgravuniformregularitytheorem}}, the following list of alternatives is exhaustive.  
\begin{indentenum}
\item[(i)] $\mathfrak{Q}_i(s) \to \infty$, as $s \to \infty$; or
\item[(ii)] $\min_T \kappa[\mathfrak{h}_i(s)] \to -\infty$, as $s \to \infty$; or
\item[(iii)] $\max_{\overline{R}} \partial_p \mathfrak{h}_i(s) \to \infty$, as $s \to \infty$; or
\item[(iv)] $\min_{\overline{R}} \partial_p \mathfrak{h}_i(s) \to \delta$, as $s \to \infty$; or
\item[(v)] $\min_T \left(\mathfrak{Q}_i(s) -2g\rho \mathfrak{h}_i(s) - 2\sigma \kappa[\mathfrak{h}_i(s)]\right) \to \delta$, as $s \to \infty$; or
\item[(vi)] $\mathcal{K}_{\delta,i}$ is a closed loop.
\end{indentenum}
Here $(\mathfrak{Q}_i(\cdot), \mathfrak{h}_i(\cdot))$ is a global parameterization of $\mathcal{K}_{\delta,i}$ in the sense of \textsc{Theorem \ref{tbglobalbifurcation1}} (by abuse of notation we have suppressed the dependency of the parameterization on $\delta$).

In light of \textsc{Lemma \ref{equivalencelemma}}, to $\mathcal{K}_{\delta, i}^\prime$ there is associated a path-connected set $\mathcal{K}_{\delta,i}$ of solutions to the problem in its Eulerian formulation \eqref{incompress}-\eqref{boundcond}. By the same arguments as in the proof of \textsc{Theorem \ref{mainresult1}}, we deduce that the first three alternatives (i)--(iii) imply that one of the following occurs:  
\[ \sup_{\overline{D_{\eta_n}}} u_n \to c \qquad \textrm{or} \qquad \inf_{\overline{D_{\eta_n}}} u_n \to -\infty,\]
for some sequence $\{(u_n, v_n, \rho_n, \eta_n, Q_n)\} \subset \mathcal{K}_{\delta,i}$.  The same holds true, again by the considerations of \textsc{Theorem \ref{mainresult1}}, if either alternative (iv) or alternative (v) hold for a decreasing sequence $\delta_n$ with $\delta_n \to 0$.  If we take $\mathcal{K}_i := \bigcup_{\delta > 0} \mathcal{K}_{\delta,i}$ and $\mathcal{K} := \bigcup_i \mathcal{K}_i$, then the theorem is proved.  
\end{proofof}

\begin{longremark}  Naively, one might hope that since each $\mathcal{K}_{\delta, i}^\prime$ has an analytic parameterization, the same might be true of the associated curves $\mathcal{K}_{\delta, i}$.   Unfortunately, this does not seem to be true in general.  The difficulty lies in \textsc{Lemma \ref{equivalencelemma}}, where the key step in transitioning from the height equation formulation back to the Eulerian is to reconstruct $\psi$ from $(h,Q)$.  In \cite{constantin2004exact,walsh2009stratified}, this is accomplished as follows.  Let $(\mathfrak{Q}, \mathfrak{h})$ be an analytic parameterization of $\mathcal{K}_{\delta, i}^\prime$ and, fixing $s \in \mathbb{R}^+$, let $h = \mathfrak{h}(s)$.  In view of \eqref{uvhphq}, for each $x_0 \in \mathbb{R}$, we define $\psi(x_0,\cdot; s)$ to be the unique solution of the first-order ODE: 
\be \left\{ \begin{array}{l} \psi_y(x_0,y;s) = -h_p^{-1}\left(x_0, -\psi(x_0,y;s)\right) \\
\psi(x_0, h(x_0,0)-d(h);s) = 0. \end{array} \right. \label{equivalenceode} \ee 
Following the arguments of \cite{constantin2004exact,walsh2009stratified}, we can show that this $\psi$ indeed corresponds to the pseudo-stream function for a solution of \eqref{incompress}--\eqref{boundcond}, thus obtaining \textsc{Lemma \ref{equivalencelemma}} for fixed $s$.   Varying $s$ gives a curve of pseudo-stream functions and the path-connected set $\mathcal{K}_{\delta, i}$. But, since $\mathfrak{h}(s)$ need not be analytic in $p$, we \emph{cannot} conclude from \eqref{equivalenceode} that $\psi(\cdot,\cdot;s)$ is analytic in $s$. \label{eulernotanalyticremark} \end{longremark}

\section{Examples} \label{examplessection}

In this section, we consider the implications of our results in a a few of the most important special cases for the given functions $\rho$ and $\beta$.  

\subsection{Pure capillary waves}  Examining \eqref{euler2}--\eqref{boundcond}, it is apparent that the contribution of surface tension comes in the form of a term involving second-order spatial derivatives, whereas the gravitational constant appears as the coefficient of  first-order derivative terms.   When considering flows in which the wave length is small with respect to $\sigma$, therefore, we expect the motion to be driven largely by the effects of capillarity.  In such circumstances, it is often appropriate to discount the presence of gravity entirely, both in the interior of the fluid and on the free surface.   A common example is that of a pebble being thrown into a quiescent body of water; the ripples that result are typically small enough in width as to be accurately modeled as pure capillary waves.  

 Strictly speaking, of course, pure capillary waves lie outside of the regime we are considering in this paper, as we have taken $g$ to be a strictly positive constant.  Nonetheless, our arguments can be easily altered to include this scenario. Indeed, defining the pseudostream function and proceeding as before, one arrives at a height equation formulation analogous to \eqref{heighteq}:
\be \left \{ \begin{array}{lll}
(1+h_q^2)h_{pp} + h_{qq}h_p^2 - 2h_q h_p h_{pq}  = -h_p^3 \beta(-p) & p_0 < p < 0, \\
1+h_q^2 + h_p^2( 2\sigma\kappa[h]- Q)  = 0 & p = 0, \\
h = 0 & p = p_0. \end{array} \right. \label{purecapheighteq} \ee
As one would expect --- and, indeed, \eqref{purecapheighteq} confirms --- ignoring gravitational effects completely negates the significance of stratification.  That is, for every solution of \eqref{purecapheighteq} and any choice of streamline density function, there corresponds a solution to \eqref{incompress}--\eqref{boundcond} (cf. \cite{yih1965dynamics}).  

The system of equations \eqref{purecapheighteq} was first considered by Wahl\'en in \cite{wahlen2006capillary}, where he used it to give a systematic study of small amplitude solutions.  In his thesis, Wahl\'en expanded upon these results, proving a global bifurcation theorem in the same vein as \textsc{Theorem \ref{capgravglobalbifurcationtheorem}} (cf. \cite{wahlenthesis}).  As an easy byproduct of our work in sections \ref{analyticitysection}, \ref{uniformregularitysection} and \ref{finalsection}, we are now able to significantly extend upon this result.  

Let us first begin by presenting the key results for the local theory developed in \cite{wahlen2006capillary}.  In keeping with our previous notation, we have written these in terms of the Bernoulli function $\beta$, but it is worth noting that, since $\rho_p \equiv 0$ for pure capillary waves, $\beta$ reduces to the familiar vorticity function $\gamma$ of \cite{constantin2004exact,constantin2007rotational,wahlen2006capillary,wahlen2006capgrav,wahlen2007onrotational}.

The first step is to establish the existence of a curve of laminar flow solutions, i.e. those where $h_q \equiv 0$.  In light of the boundary condition on the top, this task is somewhat simpler than for capillary-gravity waves.  In fact, there is no need to introduce an additional parametrizing variable.

\begin{lemma} \emph{(Wahl\'en, \cite{wahlen2006capillary})} For each $Q$ with $0 \leq -2B_{\textrm{min}} < Q$, the corresponding laminar flow solution $H = H(p; Q)$ of \eqref{heighteq} is given by
\[ H(p;Q) := \int_{p_0}^p \frac{ds}{\sqrt{Q+2B(s)}}, \qquad p_0 \leq p \leq 0.\]
\label{purecapillarylaminarflow} \end{lemma}

The existence of small amplitude solutions is proved by a bifurcation argument as was carried out in section \ref{capgravlocalbifurcationsection}, with the end result being the following.  

\begin{theorem} \emph{(Wahl\'en, \cite{wahlen2006capillary})} Let the wave speed $c >0$, the relative mass flux $p_0 < 0$, coefficient of surface tension $\sigma$ and the Bernoulli function $\beta \in C^{1+\alpha}([0,|p_0|])$ and $\alpha \in (0,1)$ be given. Then, if
\be   \sup_{p_1 \in M} \left( \frac{\sigma p_1^2 - \int_{p_1}^0 (p-p_1)^2 (2B(p) - 2B_{\textrm{min}})^{1/2} dp}{ \int_{p_1}^0 ( 2B(p) - 2B_{\textrm{min}} )^{3/2} dp } \right)^{1/2} > 1, \label{purecapillarylbc1} \ee 
where
\[ M := \left\{ p_1 \in [p_0, 0] : \sigma p_1^2 > \int_{p_1}^0 (p-p_1)^2 (2B(p) - 2B_{\textrm{min}}) dp \right\},\]
there exists a $C^1$-curve $\mathcal{D}_{\textrm{loc}} \subset \mathscr{S}$ of small amplitude traveling wave solutions $(u, v, \eta,Q)$ of \eqref{incompress}-\eqref{bedcond} with period $2\pi$, speed $c$ and relative mass flux $p_0$ , satisfying $u <c$ throughout the fluid. The curve $\mathcal{D}_{\textrm{loc}}$, moreover, exhibits the structural properties of \textsc{Theorem \ref{capgravsimplelocalbifurcation}} \emph{(i)}--\emph{(iii)}.  \label{purecapillarylocal} \end{theorem}
Note that the statement above differs from the theorem presented in \cite{wahlen2006capillary} in two ways. First, we have specialized to the case where $L = 2\pi$.  This is done simply as a matter of taste and to stress similarities with \cite{constantin2004exact} and \cite{walsh2009stratified}.  More importantly, we have increased the required regularity of the solutions $(u,v,\eta) \in \mathcal{D}_{\textrm{loc}}$ and the vorticity function by one (that this is valid is straightforward to show, by the remarks following the statement of \eqref{mainresult1}.)  The motivation for working with more regular solutions is merely to allow us to later employ \textsc{Theorem \ref{capgravuniformregularitytheorem}}.  

Also, let us draw attention to the fact that, for pure capillary waves, the bifurcation from $\mathcal{T}$ is always simple.  This is to be expected as sending $g$ to zero in \eqref{defsigmac}, we find $\sigma_c = 0$ so that $\Sigma_1 = \mathbb{R}^+$ and $\Sigma_2\cup \Sigma_3 = \emptyset$.  Interestingly, then, double bifurcation is truly a capillary-gravity wave phenomenon:  it is absent in both pure capillary wave and gravity wave regimes.

Now, observe that the arguments of section \ref{globalbifurcationsection} does not in any way rely on assuming $g$ to be strictly positive.  In fact, by taking $g = 0$ and replacing the local theory of section \ref{capgravlocalbifurcationsection} with that above, we recreate exactly the global bifurcation theorem for pure capillary waves proved by Wahl\'en in \cite{wahlenthesis}.  Let $\delta > 0$ be given and define $\mathcal{D}_\delta^\prime$, $\mathcal{O}_\delta$ and $S_\delta$ as in section \ref{globalbifurcationsection}, but with $\mathcal{D}_{\textrm{loc}}$ in place of $\mathcal{C}_{\textrm{loc}}$. 
\begin{theorem} \emph{(Wahl\'en, \cite{wahlenthesis})} Let $\delta > 0$ be given and let $\mathcal{D}_\delta$ be as described above.  One of the following alternatives must hold: 
\begin{indentenum}
 \item[\emph{(i)}] $\mathcal{D}_\delta^\prime$ is unbounded in $\mathbb{R} \times X$.
 \item[\emph{(ii)}] $\mathcal{D}_\delta^\prime$ intersects $\mathcal{T}$ at more than one point.
 \item[\emph{(iii)}] $\mathcal{D}_\delta^\prime$ contains a point $(Q,h) \in \partial \mathcal{O}_{\delta}$. \end{indentenum} 
\label{purecapglobalbifurcationtheorem}
\end{theorem}

Fortunately, our arguments in the remaining sections are similarly indifferent to the strict positivity of $g$.  For instance, taking $g = 0$ and running through the proofs of section \ref{uniformregularitysection} verbatim yields a theorem identical to \textsc{Theorem \ref{capgravuniformregularitytheorem}} but for pure capillary waves.  One way of interpreting this robustness is to note that, from a technical standpoint, the difficulty introduced by stratification and capillarity are neatly confined to the interior and boundary, respectively.  They must therefore be treated more-or-less independently, and hence we have already done the heavy lifting necessary to obtain the corresponding results for the pure capillary case.  

Combining these results (and likewise generalizing the proof of \textsc{Theorem \ref{mainresult1}}) yields the following.

\begin{theorem} Fix a wave speed $c > 0$, wavelength $L > 0$, relative mass flux $p_0 < 0$ and coefficient of surface tension $\sigma > 0$.  Fix any $\alpha \in (0,1)$, and let the function $\beta \in C^{3+\alpha}([0,|p_0|])$ be given such that condition \eqref{purecapillarylbc1}.  Let $\rho \in C^{1+\alpha}([p_0, 0]; \mathbb{R}^+)$ be any choice of streamline density function. 

There exists a connected set $\mathcal{D} \subset \mathscr{S}$ of solutions $(u,v,\rho, \eta,Q)$ of the traveling, stratified, pure capillary wave problem.  The set $\mathcal{D}$, moreover, exhibits all of the structural properties described in \textsc{Theorem \ref{mainresult1}} \emph{(a)}--\emph{(c)}.  
\label{purecapmaintheorem} 
\end{theorem}

This theorem represents a substantial strengthening of \textsc{Theorem \ref{purecapglobalbifurcationtheorem}}, elucidating as it does the meaning of alternatives (i) and (iii) in the Eulerian framework, as well as proving the existence of a path-connected subset.

\subsection{Irrotational homogeneous capillary-gravity waves}  This simplest case to consider with $\sigma > 0$ is when $\beta \equiv 0$ and $\rho \equiv \rho_0$, a positive constant.  Since $\rho_p \equiv 0$, the laminar problem simplifies greatly and we are able to solve it explicitly.  In particular, for $\lambda \geq 0$, $p_0 < p < 0$, we have
\[ H(p;\lambda) = \frac{p-p_0}{\sqrt{\lambda}}, \qquad Q(\lambda) = \lambda + \frac{2g \rho_0 |p_0|}{\sqrt{\lambda}}.\]
Thus $\lambda_0 = (g\rho_0|p_0|)^{2/3}$ and $a = \sqrt{\lambda}$.  Furthermore, since we are working in the constant density regime, $\lambda_c$ must coincide with $\lambda_0$, whence 
\[ \sigma_c = \frac{1}{3}(g\rho_0 )^2 |p_0|^3 \lambda_0^{-5/2} = \frac{1}{3} (g \rho_0)^{1/3} |p_0|^{4/3}.\]

Arguing as we have in the previous sections, we see that there exists a local curve of non-laminar solutions of period $2\pi/n$ bifurcating from $\mathcal{T}$ at $(H(\cdot; \lambda_*), \lambda_*)$, for some $n \in \mathbb{N}^{\times}$, $\lambda_* > 0$, provided there exists $M = M(p)$ satisfying
\[ \left\{\begin{array}{l} a^3 M^{\prime\prime} = n^2 a M, \qquad p_0 < p < 0, \\
M(p_0) = 0, \\ \lambda_*^{3/2} M^\prime(0) = (g\rho_0+n^2\sigma)M(0). \end{array} \right.\]
From the first equation above we infer that any such $M$ must take the form
\[ M(p) = \sinh{\left(\frac{n(p-p_0)}{\sqrt{\lambda}}\right)}, \qquad p_0 < p < 0.\] 
Imposing the boundary conditions on the top, it follows that bifurcation will occur for $(n,\lambda_*)$ satisfying 	
the following relation 
\be \lambda_*= \frac{n^2 \sigma+g\rho_0}{n} \tanh{\left(\frac{n|p_0|}{\sqrt{\lambda_*}}\right)}. \label{lambdanirrotationalcase} \ee
Note that, for fixed $n$, there is always a unique choice of  $\lambda_*$ for which \eqref{lambdanirrotationalcase} will hold.  This corresponds to the fact that, for constant density and zero vorticity, the Capillary-Gravity Local Bifurcation Condition is satisfied automatically.  However, it may very well be the case that $(n_1, \lambda_*)$ and $(n_2, \lambda_*)$ each satisfy \eqref{lambdanirrotationalcase} for some $\lambda_* > 0$ and distinct $n_1, n_2 \in \mathbb{N}^{\times}$.  

For definiteness, take $n_1 = 1$ and let $\lambda_*$ be chosen so that $(1, \lambda_*)$ satisfies relation \eqref{lambdanirrotationalcase}.  If there is no $n_2 >1$ such that \eqref{lambdanirrotationalcase}  holds for $(n_2, \lambda_*)$, then $\sigma \in \Sigma_1$ and we are in the simple bifurcation setting described by \textsc{Theorem \ref{capgravsimplelocalbifurcation}}.  In view of \textsc{Lemma \ref{eigenvaluelemma}} we know that simple bifurcation will occur whenever $\sigma \geq \sigma_c$.   This fact can be seen directly from \eqref{lambdanirrotationalcase} since, for any such $\sigma$, we must have that $\lambda$ is strictly increasing as a function of $n$.  

On the other hand, there is a set $\Sigma_2$ of $\sigma$ lying to the left of $\sigma_c$ for which there \emph{will} be an $n_2 > 1$ such that $(n_2, \lambda_*)$ satisfies relation \eqref{lambdanirrotationalcase}.   In this case, double bifurcation occurs at $(H_*, Q(\lambda_*))$ and we are in the setting of \textsc{Theorem \ref{mainresult2}}.  As we might expect, the set $\Sigma_2$ is highly non-generic.  Indeed, the following result shows that it has measure zero.  

\begin{proposition} For homogeneous, irrotational, capillary-gravity waves, $\Sigma_2$ is countable.  In fact, it consists entirely of isolated points.  \label{genericproposition}
\end{proposition}
\begin{proof2}
Let $\sigma_* \in \Sigma_2$ be given.  Since \eqref{lb} is always satisfied for irrotational waves, for each $\sigma \in \mathbb{R}^+$, there exists a unique $\lambda_* > \epsilon_0$ for which 
\[ -(a^3 \phi_1^\prime)^\prime + g\rho^\prime \phi_1 = \mu a \phi_1, ~\textrm{for } p_0 < p < 0, \qquad \phi(p_0) = 0,~\lambda_*^{3/2} \phi^\prime_1(0) = (-\mu \sigma + g\rho(0)) \phi_1(0).\]
has a (nontrivial) solution $\phi_1$, where $\mu = \mu_1 = -1$.   By continuous dependence on parameters, we may view $\lambda_*$ as a smooth function of $\sigma$.  

As $\sigma_*$ is a point of double bifurcation, there exists precisely one other nonpositive integer, $\mu_2 = -n_2^2$, for which the above eigenvalue problem has nontrivial solution with $\mu = \mu_2$.  Clearly $n_2 \neq 1$, since this would entail simple bifurcation.  Let us first consider the case where $n_2 \geq 2$.

Again appealing to continuous dependence, we may consider the negative, real eigenvalues $\mu_1$ and $\mu_2$ to be smooth functions of $(\lambda, \sigma)$ defined in a sufficiently small neighborhood $\mathcal{U}_1 \times \mathcal{U}_2$ of $(\lambda_*,\sigma_*)$.  Unravelling definitions, we see that $\sigma \in \mathcal{U}_2$ is an element of $\Sigma_2$ if and only if $\sqrt{-\mu_2(\lambda_*(\sigma),\sigma)}$ is a positive integer.  It follows that to show $\Sigma_2$ is countable it suffices to confirm
\[ \frac{\partial (\mu_2)}{\partial \sigma} \bigg|_{(\lambda_*,\sigma_*)} \neq 0.\]

We now prove this claim.  Evaluating  \eqref{lambdanirrotationalcase} for $n = 1$, we find
\be \lambda_* = (\sigma + g\rho_0) \tanh{\left( \frac{|p_0|}{\sqrt{\lambda_*}} \right)}.  \label{generic1} \ee 
Differentiating the relation in $\sigma$ and rearranging terms, we are able to compute
\be \begin{split}
\frac{d \lambda_*}{d\sigma} &= \tanh{\left(\frac{|p_0|}{\sqrt{\lambda_*}}\right)} \left(1 + (\sigma+g\rho_0) \frac{|p_0|}{2\lambda_*^{3/2}} \sech^2{\left(\frac{|p_0|}{\sqrt{\lambda_*}}\right)} \right)^{-1}\\
& = \lambda_* \left(\sigma + g\rho_0 + (\sigma+g\rho_0)^2 \frac{|p_0|}{2\lambda_*^{3/2}} \sech^2{\left(\frac{|p_0|}{\sqrt{\lambda_*}}\right)} \right)^{-1}, \end{split} \label{dlambdastardsigma} \ee
which is strictly positive.  (In fact, we already knew this in the general case by \textsc{Lemma \ref{zeroeigenvaluelemma}}.)

Let us denote $n(\sigma) := \sqrt{-\mu_2(\lambda_*(\sigma),\sigma)}$, for each $\sigma \in \mathcal{U}_2$ (in particular, observe that $n(\sigma_*) = n_2$.)  From \eqref{lambdanirrotationalcase}, therefore, we have
\be n \lambda_* = (n^2 \sigma + g\rho_0) \tanh{\left(\frac{n|p_0|}{\sqrt{\lambda_*}}\right)}. \label{genericn} \ee
Differentiating \eqref{genericn} by $\sigma$, we obtain
\begin{align*} \frac{dn}{d\sigma} \lambda_* + n \frac{d\lambda_*}{d\sigma} &= \left(2n \frac{dn}{d\sigma} \sigma + n^2\right) \tanh{\left(\frac{n|p_0|}{\sqrt{\lambda_*}}\right)} \\
& \qquad + (n^2 \sigma + g \rho_0) \sech^2{\left(\frac{n|p_0|}{\sqrt{\lambda_*}}\right)} \left( \frac{|p_0|}{\sqrt{\lambda_*}} \frac{dn}{d\sigma} - \frac{n|p_0|}{2\lambda_*^{3/2}} \frac{d\lambda_*}{d\sigma} \right).\end{align*}
Upon combining like terms, this simplifies to
\[ C_1 \frac{d n}{d \sigma} = C_2 \label{dndsigma}, \]
where 
\be 
\begin{split}
C_1 &:= \lambda_*-2n\sigma \tanh{\left(\frac{n|p_0|}{\sqrt{\lambda_*}}\right)} - (n^2 \sigma + g \rho_0) \frac{|p_0|}{\sqrt{\lambda_*}} \sech^2{\left(\frac{n|p_0|}{\sqrt{\lambda_*}}\right)}, \\
C_2 &:= n^2 \tanh{\left(\frac{n|p_0|}{\sqrt{\lambda_*}}\right)} - n \frac{d\lambda_*}{d \sigma} - (n^2 \sigma + g\rho_0) \frac{n|p_0|}{2 \lambda_*^{3/2}} \sech^2{\left(\frac{n|p_0|}{\sqrt{\lambda_*}}\right)} \frac{d\lambda_*}{d\sigma}.  \end{split}
\label{genericdefC1C2}
\ee
All that remains is to prove $C_2 \neq 0$.  


Using \eqref{dlambdastardsigma}, we may write
\be \begin{split}
C_2 &= n \left[ n \tanh{\left(\frac{n|p_0|}{\sqrt{\lambda_*}}\right)} - \frac{d\lambda_*}{d \sigma}\left(1 + (n^2 \sigma + g\rho_0) \frac{|p_0|}{2 \lambda_*^{3/2}} \sech^2{\left(\frac{n|p_0|}{\sqrt{\lambda_*}}\right)} \right) \right] \\
& = n \tanh{\left(\frac{|p_0|}{\sqrt{\lambda_*}}\right)} \left[ \frac{\displaystyle n \tanh{\left(\frac{n|p_0|}{\sqrt{\lambda_*}}\right)}}{\displaystyle\tanh{\left(\frac{|p_0|}{\sqrt{\lambda_*}}\right)}} - \frac{\displaystyle 1 + (n^2 \sigma + g\rho_0) \frac{|p_0|}{2 \lambda_*^{3/2}} \sech^2{\left(\frac{\displaystyle n|p_0|}{\displaystyle\sqrt{\lambda_*}}\right)}}{\displaystyle 1 + ( \sigma + g\rho_0) \frac{|p_0|}{2 \lambda_*^{3/2}} \sech^2{\left(\frac{\displaystyle |p_0|}{\displaystyle\sqrt{\lambda_*}}\right)}}\right] \\
& =: n \tanh{\left(\frac{|p_0|}{\sqrt{\lambda_*}}\right)} \left[ C_3 - C_4 \right].
\end{split} \label{genericC2computation1} \ee
Returning to \eqref{generic1} and \eqref{genericn}, we see
\[ (n^2 \sigma + g\rho_0) \sech^2{\left(\frac{n|p_0|}{\sqrt{\lambda_*}}\right)} = n\lambda_* \sech{\left(\frac{n|p_0|}{\sqrt{\lambda_*}}\right)} \csch{\left(\frac{n|p_0|}{\sqrt{\lambda_*}}\right)}, \]
\[  ( \sigma + g\rho_0) \sech^2{\left(\frac{|p_0|}{\sqrt{\lambda_*}}\right)} = \lambda_* \sech{\left(\frac{|p_0|}{\sqrt{\lambda_*}}\right)} \csch{\left(\frac{|p_0|}{\sqrt{\lambda_*}}\right)}.\]
Thus,  $C_4$ can be estimated
\be
\begin{split} 
C_4 & = \frac{\displaystyle \sinh{\left(\frac{|p_0|}{\sqrt{\lambda_*}} \right)}\cosh{\left(\frac{|p_0|}{\sqrt{\lambda_*}} \right)}}{\displaystyle \sinh{\left(\frac{|p_0|}{\sqrt{\lambda_*}} \right)}\cosh{\left(\frac{|p_0|}{\sqrt{\lambda_*}} \right)} + \frac{|p_0|}{2\lambda_*^{3/2}}}+ n \frac{\displaystyle \sinh{\left(\frac{|p_0|}{\sqrt{\lambda_*}} \right)} \cosh{\left(\frac{|p_0|}{\sqrt{\lambda_*}} \right)} }{\displaystyle \sinh{\left(\frac{n|p_0|}{\sqrt{\lambda_*}} \right)} \cosh{\left(\frac{n|p_0|}{\sqrt{\lambda_*}} \right)} } \\
& \leq 1+C_3 \sinh^2{\left(\frac{|p_0|}{\sqrt{\lambda_*}} \right)}\csch^2{\left(\frac{n|p_0|}{\sqrt{\lambda_*}} \right)}. \end{split} \label{genericC4computation}
\ee
Observe, however, that the map 
\[ n \mapsto \sinh^2{\left(\frac{|p_0|}{\sqrt{\lambda_*}} \right)}\csch^2{\left(\frac{n|p_0|}{\sqrt{\lambda_*}} \right)} \]
has a global upper bound  of $1/4$ for $n \geq 2$ and any values of $|p_0|, \lambda_* > 0$.  Also, since $\tanh{(n|p_0|/\sqrt{\lambda_*})} > \tanh{(|p_0|/\sqrt{\lambda_*})}$, we have that $C_3 > n \geq 2$.  From \eqref{genericC4computation} it then follows that $C_4 < 3C_3/4$, whence \eqref{genericC2computation1} implies $C_2 > 0$. Altogether, then, we have proved that $dn/d\sigma < 0$.  It follows that the points $\sigma \in \Sigma_2$ where $n_2 \geq 2$ are isolated. \end{proof2}


Let us investigate further the case when $\sigma \in \Sigma_2$. For simplicity, let us suppose $n_2 > 2$ and fix $\lambda = \lambda_*$.  Note that, by evaluating the right-hand side of \eqref{lambdanirrotationalcase} for $n = 1$ and $n = 2$, we find that nondeneracy condition \eqref{capgravnondegeneracycondition1} is equivalent to 
\[ \frac{4 \sigma + g\rho_0}{2\sigma + 2g\rho_0} \tanh{\left(\frac{2|p_0|}{\sqrt{\lambda_*}}\right)} \mathrm{\,cotanh\,}{\left(\frac{|p_0|}{\sqrt{\lambda_*}}\right)} \neq 1.\] 
Now, recalling our earlier notation, we compute, for $i = 1,2$, 
\be \begin{split}
\frac{1}{\pi} \Psi_{ii} &= -\frac{ n_i^2 |p_0| }{2\sqrt{\lambda}} -g\rho_0 \lambda^{-1} \sinh^{2}{\left(\frac{n_i |p_0|}{\sqrt{\lambda}}\right)} - \sigma n_i^2\lambda^{-1} \\
& = -\frac{n_i^2 |p_0|}{2 \sqrt{\lambda}} -\frac{n_i}{2} \sinh{\left(\frac{2n_i|p_0|}{\sqrt{\lambda}}\right)},
\end{split} \label{Psiiexample} \ee

\be \begin{split}
\frac{1}{\pi}\Theta_{iiii} & = \frac{n_i^2 \sqrt{\lambda}}{4}  \left( \frac{7|p_0|}{4} - \frac{\sqrt{\lambda}}{4 n_i} \sinh{\left(\frac{2n_i|p_0|}{\sqrt{\lambda}}\right)} - \frac{5 \sqrt{\lambda}}{16 n_i} \sinh{\left(\frac{4n_i|p_0|}{\sqrt{\lambda}}\right)} \right) \\
& \qquad + \frac{3 n_i^2}{4 \sqrt{\lambda}} \sinh{\left(\frac{n_i|p_0|}{\sqrt{\lambda}}\right)} \cosh^3{\left(\frac{n_i|p_0|}{\sqrt{\lambda}}\right)},
\end{split}
\label{Thetaiiiiexample} \ee


\be
\begin{split}
\frac{1}{\pi} \Theta_{iijj} & =  \frac{1}{2} \sqrt{\lambda} \left( \frac{1}{2} n_i^2 n_j^2 |p_0| - \frac{1}{4}n_i n_j^2 \sqrt{\lambda} \sinh{\left(\frac{2n_i |p_0|}{\sqrt{\lambda}}\right)} \right) \\
& \qquad + n_i n_j^3 \sqrt{\lambda} \left( \frac{\displaystyle \sinh{\left(\frac{2(n_i-n_j)|p_0|}{\sqrt{\lambda}}\right)}}{8(n_j - n_i)} + \frac{\sinh{\left(\frac{\displaystyle 2(n_i+n_j)|p_0|}{\displaystyle\sqrt{\lambda}}\right)}}{8(n_i+n_j)} \right) \\
& \qquad\qquad + \frac{n_i n_j^2}{2\sqrt{\lambda}}  \sinh{\left(\frac{n_i |p_0|}{\sqrt{\lambda}}\right)} \cosh{\left(\frac{n_i |p_0|}{\sqrt{\lambda}}\right)} \cosh^2{\left(\frac{n_j |p_0|}{\sqrt{\lambda}}\right)}. 
\end{split}
\label{Thetaiijjexample} \ee
Observe that the quantity in \eqref{Psiiexample} is strictly negative, as expected.  The signs and relative sizes of the $\Theta_{iiii}$ and $\Theta_{iijj}$ are more difficult to ascertain in general, but inserting \eqref{Psiiexample}--\eqref{Thetaiijjexample} into \eqref{capgravnondegeneracycondition2} yields an explicitly verifiable condition under which \textsc{Theorem \ref{mainresult2}} holds.  We omit the details.  

To obtain the dispersion relation for these solutions, let $U_*$ denote the horizontal velocity corresponding to the laminar flow $H_*$.  Then, by the definition of $p_0$, 
\[ |p_0| = \int_{-d}^{0} \sqrt{\rho_0}\left(c-U_*\right) dy = \int_{-d}^0 \left(\partial_p H_* \right)^{-1} dy =  d \sqrt{\lambda_*}.\]
Writing \eqref{lambdanirrotationalcase} in terms of the depth and the horizontal velocity, therefore, we obtain
\[ c-U_* = \sqrt{\frac{n^2 \sigma + g\rho_0}{n} \tanh{\left(nd\right)}}.\]

\begin{acknowledgements}
The author wishes to thank the many people who aided in the development of this work.  I am grateful to  W. Strauss for his invaluable support at every stage of this research, and to J. Mallet--Paret for considerable guidance on the material of section \ref{doubleeigenvaluebifurcationsection}.  Likewise, I am   deeply indebted to J.F. Toland, who suggested the approach of section \ref{analyticitysection} and offered much encouragement.  Finally, many thanks are due to A. Constantin, H. Dong, V. Hur and E. Wahl\'en for enlightening comments and discussions at various points along the way.    
\end{acknowledgements}

\bibliographystyle{spmpsci}
\bibliography{paperdraft5}   

\begin{thebibliography}{10}
\providecommand{\url}[1]{{#1}}
\providecommand{\urlprefix}{URL }
\expandafter\ifx\csname urlstyle\endcsname\relax
  \providecommand{\doi}[1]{DOI~\discretionary{}{}{}#1}\else
  \providecommand{\doi}{DOI~\discretionary{}{}{}\begingroup
  \urlstyle{rm}\Url}\fi

\bibitem{agmon1962eigenfunctions}
Agmon, S.: On the eigenfunctions and on the eigenvalues of general elliptic
  boundary value problems.
\newblock Comm. Pure and Appl. Math. \textbf{15}(2), 119--147 (1962)

\bibitem{buffoni1998ondes}
Buffoni, B., Dancer, E., Toland, J.: Sur les ondes de stokes et une conjecture
  de levi-civita.
\newblock C. R. Acad. Sci. Paris S{\'e}r. I Math. \textbf{326}, 1265--1268
  (1998)

\bibitem{buffoni2000regularity}
Buffoni, B., Dancer, E., Toland, J.: The regularity and local bifurcation of
  stokes waves.
\newblock Arch. Rational Mech. Anal. \textbf{152}(3), 207--240 (2000)

\bibitem{buffoni2000sub}
Buffoni, B., Dancer, E., Toland, J.: The sub-harmonic bifurcation of stokes
  waves.
\newblock Arch. Rational Mech. Anal. \textbf{152}(3), 241--271 (2000)

\bibitem{buffoni2003analytic}
Buffoni, B., Toland, J.: Analytic theory of global bifurcation: an
  introduction.
\newblock Princeton University Press (2003)

\bibitem{constantin2004exact}
Constantin, A., Strauss, W.: Exact steady periodic water waves with vorticity.
\newblock Comm. Pure Appl. Math. \textbf{57}(4), 481--527 (2004)

\bibitem{constantin2007rotational}
Constantin, A., Strauss, W.: Rotational steady water waves near stagnation.
\newblock Philos. Trans. Roy. Soc. London Ser. A \textbf{365}(1858), 2227--2239
  (2007)

\bibitem{crandall1971bifurcation}
Crandall, M.G., Rabinowitz, P.H.: Bifurcation from simple eigenvalues.
\newblock J. Func. Anal. \textbf{8}, 321--340 (1971)

\bibitem{dancer1973bifurcation}
Dancer, E.: Bifurcation theory for analytic operators.
\newblock Proc. London Math. Soc. \textbf{26}, 359--384 (1973)

\bibitem{dancer1973globalsolution}
Dancer, E.: Global solution branches for positive mappings.
\newblock Arch. Rational Mech. Anal. \textbf{52}(2), 181--192 (1973)

\bibitem{dancer1973globalstructure}
Dancer, E.: Global structure of the solutions of nonlinear real analytic
  eigenvalue problems.
\newblock Proc. London Math. Soc \textbf{27}, 747--765 (1973)

\bibitem{dubreil1934determination}
Dubreil-Jacotin, M.: Sur la determination rigoureuse des ondes permanentes
  periodiques d'ampleur finie.
\newblock J. Math. Pures Appl. \textbf{13}(3), 217--291 (1934)

\bibitem{dubreil1937theoremes}
Dubreil-Jacotin, M.: Sur les theoremes d'existence relatifs aux ondes
  permanentes periodiques a deux dimensions dans les liquides heterogenes.
\newblock J. Math. Pures Appl. \textbf{16}(9), 43--67 (1937)

\bibitem{healey1998global}
Healey, T., Simpson, H.: Global continuation in nonlinear elasticity.
\newblock Arch. Rational Mech. Anal. \textbf{143}(1), 1--28 (1998)

\bibitem{jones1986symmetry}
Jones, M., Toland, J.: Symmetry and the bifurcation of capillary-gravity waves.
\newblock Arch. Rational Mech. Anal. \textbf{96}(1), 29--53 (1986)

\bibitem{kielhofer1985multiple}
Kielh{{\"o}}fer, H.: Multiple eigenvalue bifurcation for {F}redholm operators.
\newblock J. Reine Angew. Math. \textbf{358}, 104--124 (1985)

\bibitem{kielhofer2004bifurcation}
Kielh{{\"o}}fer, H.: Bifurcation theory, \emph{Applied Mathematical Sciences},
  vol. 156.
\newblock Springer-Verlag, New York (2004).
\newblock An introduction with applications to PDEs

\bibitem{kinsman}
Kinsman, B.: Wind Waves.
\newblock Prentice Hall, New Jersey (1965)

\bibitem{krylov1996lectures}
Krylov, N.: Lectures on elliptic and parabolic equations in H{\"o}lder spaces.
\newblock American Mathematical Society (1996)

\bibitem{levicivita1925determination}
Levi-Civita, T.: D{\'e}termination rigoureuse de ondes permanentes d'ampleur
  finie.
\newblock Ann. Math. \textbf{93}, 264--314 (1925)

\bibitem{luo1991linear}
Luo, Y., Trudinger, N.: Linear second order elliptic equations with venttsel
  boundary conditions.
\newblock In: Proc. R. Soc. Edinb., Sect. A, vol. 118, pp. 193--207 (1991)

\bibitem{luo1994quasilinear}
Luo, Y., Trudinger, N.: Quasilinear second order elliptic equations with
  venttsel boundary conditions.
\newblock Potential Anal. \textbf{3}(2), 219--243 (1994)

\bibitem{mei1989applied}
Mei, C.: The applied dynamics of ocean surface waves.
\newblock World Scientific Pub Co Inc (1989)

\bibitem{nekrasov1951exact}
Nekrasov, A.I.: The exact theory of steady waves on the surface of a heavy
  fluid.
\newblock Izdat. Akad. Nauk SSSR, Moscow (1951)

\bibitem{okamoto1990on}
Okamoto, H.: On the problem of water waves of permanent configuration.
\newblock Nonlinear Anal. \textbf{14}(6), 469--481 (1990)

\bibitem{okamoto1996resonance}
Okamoto, H., Sh{\=o}ji, M.: The resonance of modes in the problem of
  two-dimensional capillary-gravity waves.
\newblock Physica D: Nonlinear Phenomena \textbf{95}(3-4), 336--350 (1996)

\bibitem{rabinowitz1971some}
Rabinowitz, P.: Some global results for nonlinear eigenvalue problems.
\newblock J. Funct. Anal \textbf{7}, 487--513 (1971)

\bibitem{schwartz1979numerical}
Schwartz, L., Vanden-Broeck, L.: Numerical solution of the exact equations for
  capillary gravity waves.
\newblock J. Fluid Mech. \textbf{95}, 119--139 (1979)

\bibitem{shoji1989new}
Sh{\=o}ji, M.: New bifurcation diagrams in the problem of permanent progressive
  waves.
\newblock J. Fac. Sci. Univ. Tokyo Sect. IA Math. \textbf{36}, 571--613 (1989)

\bibitem{toland1985bifurcation}
Toland, J., Jones, M.: The bifurcation and secondary bifurcation of
  capillary-gravity waves.
\newblock Proc. Roy. Soc. London Ser. A \textbf{399}(1817), 391--417 (1985)

\bibitem{turner2002traveling}
Turner, R.E.L.: Traveling waves in natural systems.
\newblock In: Variational and topological methods in the study of nonlinear
  phenomena ({P}isa, 2000), \emph{Progr. Nonlinear Differential Equations
  Appl.}, vol.~49, pp. 115--131. Birkh{\"a}user Boston, Boston, MA (2002)

\bibitem{wahlen2006capgrav}
Wahl{{\'e}}n, E.: Steady periodic capillary-gravity waves with vorticity.
\newblock SIAM J. Math. Anal. \textbf{38}(3), 921--943 (electronic) (2006)

\bibitem{wahlen2006capillary}
Wahl{{\'e}}n, E.: Steady periodic capillary waves with vorticity.
\newblock Ark. Mat. \textbf{44}(2), 367--387 (2006)

\bibitem{wahlen2007onrotational}
Wahl{{\'e}}n, E.: On rotational water waves with surface tension.
\newblock Philos. Trans. R. Soc. Lond. Ser. A Math. Phys. Eng. Sci.
  \textbf{365}(1858), 2215--2225 (2007)

\bibitem{wahlenthesis}
Wahl{{\'e}}n, E.: On some nonlinear aspects of wave motion.
\newblock Ph.D. thesis, Lund University (2008)

\bibitem{walsh2009stratified}
Walsh, S.: Stratified and steady periodic water waves.
\newblock SIAM J. Math. Anal. \textbf{41}(3), 1054--1105 (2009)

\bibitem{wilton1915ripples}
Wilton, J.: On ripples.
\newblock Phil. Mag. \textbf{29}, 688--700 (1915)

\bibitem{yih1965dynamics}
Yih, C.S.: Dynamics of nonhomogeneous fluids.
\newblock The Macmillan Co., New York (1965)

\end{thebibliography}

\end{document}